\documentclass{article}
\usepackage{amsmath}
\usepackage{amssymb}
\usepackage{amsthm}
\begin{document}
\def\e#1\e{\begin{equation}#1\end{equation}}
\def\eq#1{{\rm(\ref{#1})}}
\theoremstyle{plain}
\newtheorem{thm}{Theorem}[section]
\newtheorem{lem}[thm]{Lemma}
\newtheorem{prop}[thm]{Proposition}
\newtheorem{cor}[thm]{Corollary}
\theoremstyle{definition}
\newtheorem{dfn}[thm]{Definition}
\def\dim{\mathop{\rm dim}}
\def\Re{\mathop{\rm Re}}
\def\Im{\mathop{\rm Im}}
\def\diam{\mathop{\rm diam}}
\def\vol{\mathop{\rm vol}}
\def\U{\mathbin{\rm U}}
\def\SU{\mathop{\rm SU}}
\def\ge{\geqslant} 
\def\le{\leqslant} 
\def\R{\mathbin{\mathbb R}}
\def\Z{\mathbin{\mathbb Z}}
\def\C{\mathbin{\mathbb C}}
\def\al{\alpha}
\def\be{\beta}
\def\la{\lambda}
\def\ga{\gamma}
\def\de{\delta}
\def\ep{\epsilon}
\def\ka{\kappa}
\def\th{\theta}
\def\vp{\varphi}
\def\si{\sigma}
\def\ze{\zeta}
\def\De{\Delta}
\def\La{\Lambda}
\def\Om{\Omega}
\def\Ga{\Gamma}
\def\om{\omega}
\def\d{{\rm d}}
\def\pd{\partial}
\def\db{{\bar\partial}}
\def\ts{\textstyle}
\def\w{\wedge}
\def\sm{\setminus}
\def\bigot{\bigotimes}
\def\iy{\infty}
\def\ra{\rightarrow}
\def\hookra{\hookrightarrow}
\def\t{\times}
\def\ha{{\textstyle\frac{1}{2}}}
\def\ti{\tilde}
\def\ovB{\,\overline{\!B}}
\def\ms#1{\vert#1\vert^2}
\def\bms#1{\bigl\vert#1\bigr\vert^2}
\def\md#1{\vert #1 \vert}
\def\bmd#1{\big\vert #1 \big\vert}
\def\nm#1{\Vert #1 \Vert}
\def\cnm#1#2{\Vert #1 \Vert_{C^{#2}}} 
\def\lnm#1#2{\Vert #1 \Vert_{L^{#2}}} 
\def\snm#1#2#3{\Vert #1 \Vert_{L^{#2}_{#3}}} 
\def\bcnm#1#2{\bigl\Vert #1 \bigr\Vert_{C^{#2}}} 
\def\blnm#1#2{\bigl\Vert #1 \bigr\Vert_{L^{#2}}} 
\title{$\U(1)$-invariant special Lagrangian 3-folds. II.\\
Existence of singular solutions}
\author{Dominic Joyce \\ Lincoln College, Oxford}
\date{}
\maketitle

\section{Introduction}
\label{us1}

Special Lagrangian submanifolds (SL $m$-folds) are a distinguished
class of real $m$-dimensional minimal submanifolds in $\C^m$, which are
calibrated with respect to the $m$-form $\Re(\d z_1\w\cdots\w\d z_m)$.
They can also be defined in Calabi--Yau manifolds, are important in
String Theory, and are expected to play a r\^ole in the eventual
explanation of Mirror Symmetry between Calabi--Yau 3-folds.

This is the second in a suite of three papers \cite{Joyc3,Joyc4}
studying special Lagrangian 3-folds $N$ in $\C^3$ invariant under
the $\U(1)$-action
\e
{\rm e}^{i\th}:(z_1,z_2,z_3)\mapsto
({\rm e}^{i\th}z_1,{\rm e}^{-i\th}z_2,z_3)
\quad\text{for ${\rm e}^{i\th}\in\U(1)$.}
\label{us1eq1}
\e
These three papers and \cite{Joyc5} are surveyed in
\cite{Joyc6}. Locally we can write $N$ as
\e
\begin{split}
N=\Bigl\{(z_1&,z_2,z_3)\in\C^3:
\Im(z_3)=u\bigl(\Re(z_3),\Im(z_1z_2)\bigr),\\
&\Re(z_1z_2)=v\bigl(\Re(z_3),\Im(z_1z_2)\bigr),
\quad\ms{z_1}-\ms{z_2}=2a\Bigr\},
\end{split}
\label{us1eq2}
\e
where $a\in\R$ and $u,v:\R^2\ra\R$ are continuous functions. It
was shown in \cite{Joyc3} that when $a\ne 0$, $N$ is an SL 3-fold
in $\C^3$ if and only if $u,v$ satisfy
\e
\frac{\pd u}{\pd x}=\frac{\pd v}{\pd y}\quad\text{and}\quad
\frac{\pd v}{\pd x}=-2\bigl(v^2+y^2+a^2\bigr)^{1/2}\frac{\pd u}{\pd y},
\label{us1eq3}
\e
and then $u,v$ are smooth and $N$ is nonsingular.

The goal of this paper and its sequel \cite{Joyc4} is to study what
happens when $a=0$. In this case, at points $(x,0)$ with $v(x,0)=0$
the factor $-2(v^2+y^2+a^2)^{1/2}$ in \eq{us1eq3} becomes zero, and
then \eq{us1eq3} is no longer elliptic. Because of this, when $a=0$
the appropriate thing to do is to consider {\it weak solutions} of
\eq{us1eq3}, which may have {\it singular points} $(x,0)$ with
$v(x,0)=0$. At such a point $u,v$ may not be differentiable, and
$\bigl(0,0,x+iu(x,0)\bigr)$ is a singular point of the SL 3-fold
$N$ in~$\C^3$.

This paper will be concerned largely with technical analytic
issues, to do with the existence, uniqueness and regularity of
weak, singular solutions of \eq{us1eq3} in the case $a=0$. The
sequel \cite{Joyc4} will describe the singularities of solutions
of \eq{us1eq3} with $a=0$, prove that under mild conditions the
singularities are isolated and have a unique {\it multiplicity}
and {\it type}, and show that for each $k\ge 1$ singularities with
multiplicity $k$ exist and occur in codimension $k$, in some sense.

In \cite{Joyc4} we also construct large families of {\it special
Lagrangian fibrations} on open subsets of $\C^3$. These fibrations
are used in \cite{Joyc5} as local models to study special Lagrangian
fibrations of (almost) Calabi--Yau 3-folds, and to draw some
conclusions about the {\it SYZ Conjecture} \cite{SYZ}. For a brief
summary of the results of all four papers, see~\cite{Joyc6}.

In \cite{Joyc3} we showed that if $S$ is a domain in $\R^2$
and $u,v\in C^1(S)$ satisfy \eq{us1eq3}, then $v$ satisfies
\e
\frac{\pd}{\pd x}\Bigl[\bigl(v^2+y^2+a^2\bigr)^{-1/2}
\frac{\pd v}{\pd x}\Bigr]+2\,\frac{\pd^2v}{\pd y^2}=0
\label{us1eq4}
\e
in $S^\circ$, and also there exists $f\in C^2(S)$ with $\frac{\pd
f}{\pd y}=u$ and $\frac{\pd f}{\pd x}=v$, unique up to addition
of a constant, satisfying
\e
\Bigl(\Bigl(\frac{\pd f}{\pd x}\Bigr)^2+y^2+a^2
\Bigr)^{-1/2}\frac{\pd^2f}{\pd x^2}+2\,\frac{\pd^2f}{\pd y^2}=0.
\label{us1eq5}
\e
Both \eq{us1eq4} and \eq{us1eq5} are {\it second-order quasilinear
elliptic equations}. We then proved existence and uniqueness of
solutions of the Dirichlet problems for \eq{us1eq4} and \eq{us1eq5}
on (strictly convex) domains when~$a\ne 0$.

The main results of this paper are Theorems \ref{us6thm1},
\ref{us6thm2}, \ref{us7thm1} and \ref{us7thm2} below, which
prove existence and uniqueness of {\it singular} solutions to
the Dirichlet problems for \eq{us1eq4} and \eq{us1eq5} when
$a=0$. They give detailed results on the {\it regularity} of the
solutions --- basically $u,v$ are $C^0$, have weak derivatives
$\frac{\pd u}{\pd x},\frac{\pd u}{\pd y},\frac{\pd v}{\pd x},
\frac{\pd v}{\pd y}$ in $L^p$ for $p$ in various ranges, and
$u,v$ are real analytic away from their singular points.
Also, $f,u,v$ vary continuously with the boundary data.
These singular solutions give many examples of {\it singular
$\U(1)$-invariant SL\/ $3$-folds} in~$\C^3$.

A fundamental question about compact SL 3-folds $N$ in Calabi--Yau
3-folds $M$ is: {\it how stable are they under large deformations\/}?
Here we mean both deformations of $N$ in a fixed $M$, and what
happens to $N$ as we deform $M$. The deformation theory of compact
SL 3-folds under {\it small}\/ deformations is already well
understood, and is described in \cite[\S 9]{Joyc2}. But to extend
this understanding to large deformations, one needs to take into
account singular behaviour.

One possible moral of this paper is that {\it compact SL\/ $3$-folds
are pretty stable under large deformations}. That is, we have shown
existence and uniqueness for (possibly singular) $\U(1)$-invariant
SL 3-folds in $\C^3$ satisfying certain boundary conditions. This
existence and uniqueness is {\it entirely unaffected} by
singularities that develop in the SL 3-folds, which is quite
surprising, as one might have expected that when singularities
develop the existence and uniqueness properties would break down.

This is encouraging, as both the author's programme for constructing
invariants of Calabi--Yau 3-folds in \cite{Joyc1} by counting special
Lagrangian homology 3-spheres, and proving some version of the SYZ
Conjecture \cite{SYZ} in anything other than a fairly weak, limiting
form, will require strong stability properties of compact SL 3-folds
under large deformations; so these papers may be taken as a small
piece of evidence that these two projects may eventually be successful.

In \S\ref{us6} and \S\ref{us7} we shall prove {\it existence and
uniqueness of the Dirichlet problems} for \eq{us1eq4} and
\eq{us1eq5} respectively when $a=0$. That is, given a suitable
domain $S$ in $\R^2$ and boundary data $\phi\in C^{k+2,\al}(\pd S)$,
we construct unique $v\in C^0(S)$ with $v\vert_{\pd S}=\phi$
satisfying \eq{us1eq4} weakly, or unique $f\in C^1(S)$ with
$f\vert_{\pd S}=\phi$ satisfying \eq{us1eq5} with weak
second derivatives.

The basic method is this. For each $a\in(0,1]$ we let
$v_a$ or $f_a$ in $C^{k+2,\al}(S)$ be the unique solution
of \eq{us1eq4} or \eq{us1eq5} with $v_a\vert_{\pd S}=\phi$
or $f_a\vert_{\pd S}=\phi$. Then we aim to prove that 
$v_a\ra v$ in $C^0(S)$ or $f_a\ra f$ in $C^1(S)$ as $a\ra 0_+$
for some unique $v,f$ which are (weak, singular) solutions
of the Dirichlet problems for~$a=0$.

To show that these limits $v,f$ exist, the main issue is to
prove {\it a priori estimates} of $v_a,f_a$ and their derivatives
that are {\it uniform in} $a$. That is, we need bounds such as
$\cnm{v_a}0\le C$ for all $a\in(0,1]$, with $C$ independent of $a$.
Getting such estimates is difficult, since equations \eq{us1eq4}
and \eq{us1eq5} {\it really are\/} singular when $a=0$, so many
norms of $v_a,f_a$ such as $\cnm{\pd v_a}0$ can diverge to
infinity as $a\ra 0_+$, and uniform a priori estimates of these
norms {\it do not exist}.

Although there are many results on a priori estimates for
nonlinear elliptic equations in the literature, I could not
find any that told me what I needed to know about \eq{us1eq4}
and \eq{us1eq5}, so I was forced to invent my own method.
It gives a priori estimates for the first derivatives of
bounded solutions of {\it nonlinear Cauchy-Riemann equations\/}
such as \eq{us1eq3}. The underlying idea is geometrical,
and comes from complex analysis.

We shall explain \cite[Th.~6.9]{Joyc3}, which is the key tool.
Suppose $S,T$ are domains in $\R^2$, and $(u,v):T\ra S^\circ$,
$(\hat u,\hat v):S\ra T$ satisfy \eq{us1eq3} for some $a\ne 0$
with $(u,v)(x_0,y_0)=(u_0,v_0)$, $(\hat u,\hat v)(u_0,v_0)=
(x_0,y_0)$ for $(u_0,v_0)\in S^\circ$, $(x_0,y_0)\in T^\circ$.
Then the {\it graphs} $\Ga,\hat\Ga$ of $(u,v)$ and $(\hat u,
\hat v)$ are 2-submanifolds of $S\t T$ intersecting
at~$(u_0,v_0,x_0,y_0)$.

As \eq{us1eq3} is a {\it nonlinear Cauchy-Riemann equation},
it turns out there is an {\it almost complex structure} on
$S\t T$, making $\Ga,\hat\Ga$ into {\it pseudo-holomorphic
curves}. Therefore each intersection point in $\Ga\cap\hat\Ga$
has a positive integer {\it multiplicity}. Now by an argument
like those used to count zeroes of holomorphic functions in complex
analysis, by considering {\it winding numbers\/} along $\pd S$
we find that the {\it total multiplicity} of $\Ga\cap\hat\Ga$
is 1. Thus, $(u_0,\ldots,y_0)$ has multiplicity 1, so the tangent
spaces of $\Ga,\hat\Ga$ at $(u_0,\ldots,y_0)$ are {\it distinct}.

So suppose $(\hat u,\hat v):S\ra T$ satisfies \eq{us1eq3} with
$(\hat u,\hat v)(u_0,v_0)=(x_0,y_0)$. Then there {\it cannot
exist\/} any solution $(u,v):T\ra S^\circ$ of \eq{us1eq3}
with $(u,v)(x_0,y_0)=(u_0,v_0)$ and first derivatives
$\pd u,\pd v$ at $(x_0,y_0)$ taking prescribed values,
that is, those necessary to make $T_{(u_0,\ldots,y_0)}\Ga,
T_{(u_0,\ldots,y_0)}\hat\Ga$ coincide. In this way we
translate {\it existence results\/} for $(\hat u,\hat v):S\ra T$
satisfying \eq{us1eq3} with prescribed values and derivatives
at $(u_0,v_0)$ into {\it nonexistence results\/} for
$(u,v):T\ra S^\circ$ satisfying \eq{us1eq3} with
prescribed values and derivatives at~$(x_0,y_0)$.

Sections \ref{us4} and \ref{us5} use this to prove
{\it a priori estimates} for $\frac{\pd u}{\pd x},
\frac{\pd u}{\pd y},\frac{\pd v}{\pd x}$ and
$\frac{\pd v}{\pd y}$ when $(u,v)$ are bounded
solutions of \eq{us1eq3}. In \S\ref{us4} we
construct two families of solutions $(\hat u,\hat v)$
to \eq{us1eq3}, with $\hat v\approx \al+\be y+\ga x$
or $\hat v\approx \al+\be y+\ga xy$ for $\ga$ small.
Roughly speaking, these examples $(\hat u,\hat v):S\ra T$
fill out all possible values and derivatives at $(u_0,v_0)$
with $\pd\hat v$ small. Then \cite[Th.~6.9]{Joyc3} shows
that for solutions $(u,v):T\ra S^\circ$ of \eq{us1eq3},
all values and derivatives with $\pd v$ large are excluded.
So we can give a priori bounds for~$\pd u,\pd v$.

Here is another interesting analytic feature. Because of
the involvement of $y$ in \eq{us1eq3}--\eq{us1eq5}, $x$
{\it and\/ $y$ derivatives behave differently}. Roughly
speaking, we find in \S\ref{us6} that if $(u,v)$ are weak
solutions of \eq{us1eq3} when $a=0$ then using the material
of \S\ref{us5}, we can show that $\frac{\pd u}{\pd x},
\frac{\pd v}{\pd y}$ lie in $L^p$ for $p\in[1,\frac{5}{2})$,
and $\frac{\pd u}{\pd y}$ lies in $L^q$ for $q\in[1,2)$,
and $\frac{\pd v}{\pd x}$ lies in $L^r$ for~$r\in[1,\iy)$.

As $L^p_1\hookra C^0$ for $p>2$ by the Sobolev embedding
theorem, $v$ is continuous. But our $L^q$ estimate of
$\frac{\pd u}{\pd y}$ is too weak to show $u$ is continuous
in this way. Because of this, we prove a {\it nonstandard
Sobolev embedding theorem}, Theorem \ref{us2thm}. It allows
us to use a stronger $L^p$ norm of $\frac{\pd u}{\pd x}$ to
compensate for the weaker $L^q$ norm of $\frac{\pd u}{\pd y}$,
and still prove $u$ is continuous. Proving $u,v$ are continuous
is important geometrically: without continuity, $N$ in
\eq{us1eq2} would not be {\it locally closed}, and one singular
point of $u,v$ would correspond to many in~$N$.

Readers may wonder why we study both equations \eq{us1eq4}
and \eq{us1eq5}, rather than just one. The answer is that
for different applications, each may be preferable. For
instance, in \S\ref{us4} we construct our solutions with
the Dirichlet problem for \eq{us1eq5}, but estimate them
using elliptic regularity for \eq{us1eq4}. In \cite{Joyc4}
both Dirichlet problems are used, usually \eq{us1eq5}, and
properties of \eq{us1eq4} such as $v<v'$ on $\pd S$ implies
$v<v'$ on $S$ crop up continually. In \cite{Joyc5} the
Dirichlet problem for \eq{us1eq5} is best.
\medskip

\noindent{\it Acknowledgements.} I would like to thank Mark Gross,
Rafe Mazzeo and Rick Schoen for helpful conversations, and the
referee for useful suggestions. I was supported by an EPSRC
Advanced Fellowship whilst writing this paper.

\section{Background material from analysis}
\label{us2}

In \S\ref{us21} we briefly summarize some background material
we will need for later analytic results. Our principal reference
is Gilbarg and Trudinger \cite{GiTr}. Section \ref{us22} proves
a Sobolev embedding type result for functions $u$ on subsets
of $\R^2$, in terms of bounds for $\lnm{\frac{\pd u}{\pd x}}{p}$
and $\lnm{\frac{\pd u}{\pd y}}{q}$ with~$p\ne q$.

\subsection{Domains, function spaces and operators}
\label{us21}

First we define {\it domains} in $\R^2$, and {\it Banach
spaces of functions} upon them.

\begin{dfn} A {\it domain} in $\R^2$ is a compact subset
$S\subset\R^2$ which is topologically a disc with smooth
(or sometimes piecewise smooth) {\it boundary} $\pd S$.
The interior is $S^\circ=S\sm\pd S$. A convex domain $S$
is {\it strictly convex} if the curvature of $\pd S$ is
strictly positive.

For each $k\ge 0$, write $C^k(S)$ for the space of continuous
functions $f:S\ra\R$ with $k$ continuous derivatives, with norm
$\cnm{f}k=\sum_{j=0}^k\sup_S\bmd{\pd^jf}$. For $\al\in(0,1]$,
the {\it H\"older space} $C^{k,\al}(S)$ is the
subset of $f\in C^k(S)$ for which
\begin{equation*}
[\pd^k f]_\al=\sup_{x\ne y\in S}
\frac{\bmd{\pd^kf(x)-\pd^kf(y)}}{\md{x-y}^\al}
\end{equation*}
is finite, with norm $\cnm{f}{k,\al}=\cnm{f}k+[\pd^kf]_\al$.
Set~$C^\iy(S)=\bigcap_{k=0}^\iy C^k(S)$.

For $q\ge 1$, the {\it Lebesgue space} $L^q(S)$ is the set
of locally integrable functions $f$ on $S$ for which the norm
$\lnm{f}q=\bigl(\int_S\md{f}^q\d{\bf x}\bigr){}^{1/q}$ is finite.
For $k\ge 0$, the {\it Sobolev space} $L^q_k(S)$ is the set of
$f\in L^q(S)$ which are $k$ times weakly differentiable with
$\md{\nabla^jf}\in L^q(S)$ for~$j\le k$.
\label{us2def1}
\end{dfn}

Next we discuss {\it differential operators} on domains.

\begin{dfn} Let $S$ be a domain in $\R^2$. A {\it second-order 
linear differential operator} $P:C^2(S)\ra C^0(S)$ may be written
\e
\bigl(Pu\bigr)(x)=
\sum_{i,j=1}^2a^{ij}(x)\frac{\pd^2u}{\pd x_i\pd x_j}(x)
+\sum_{i=1}^2b^i(x)\frac{\pd u}{\pd x_i}(x)+c(x)u(x),
\label{us2eq1}
\e
for $u\in C^2(S)$. Here $a^{ij},b^i,c\in C^0(S)$ are
the {\it coefficients} of $P$, with $a^{ij}=a^{ji}$.

A {\it second-order quasilinear operator} $Q:C^2(S)\ra C^0(S)$
may be written
\e
\bigl(Qu\bigr)(x)=
\sum_{i,j=1}^2a^{ij}(x,u,\pd u)\frac{\pd^2u}{\pd x_i\pd x_j}(x)
+b(x,u,\pd u),
\label{us2eq2}
\e
for $u\in C^2(S)$. Here $a^{ij},b\in C^0\bigl(S\t\R\t(\R^2)^*\bigr)$
are the {\it coefficients} of $Q$, with $a^{ij}=a^{ji}$. We call $P$
or $Q$ {\it elliptic} if the symmetric $2\t 2$ matrix $(a^{ij})$ is
positive definite at every point. A second-order quasilinear operator
$Q:C^2(S)\ra C^0(S)$ is in {\it divergence form} if it is written as
\e
\bigl(Qu\bigr)(x)=\sum_{j=1}^2\frac{\pd}{\pd x_j}
\bigl(a^j(x,u,\pd u)\bigr)+b(x,u,\pd u),
\label{us2eq3}
\e
for functions $a^j\in C^1\bigr(S\t\R\t(\R^2)^*\bigr)$
and~$b\in C^0\bigr(S\t\R\t(\R^n)^*\bigr)$.
\label{us2def2}
\end{dfn}

Let $Q$ be a quasilinear operator as in \eq{us2eq2} or \eq{us2eq3}.
We shall consider three different senses in which $Qu=0$ can hold:
\begin{itemize}
\setlength{\itemsep}{0pt}
\setlength{\parsep}{0pt}
\item We say $Qu=0$ {\it holds} if $u\in C^2(S)$ and $Qu=0$
in~$C^0(S)$.
\item We say $Qu=0$ {\it holds with weak derivatives} if $u$
is twice {\it weakly differentiable}, so that $\pd u$, $\pd^2u$
exist almost everywhere, and $Qu=0$ holds almost everywhere.
\item For $Q$ in {\it divergence form} \eq{us2eq3}, we say
$Qu\!=\!0$ {\it holds weakly} if $u\!\in\!L^1_1(S)$~and
\e
-\sum_{j=1}^2\int_S\frac{\pd\psi}{\pd x_j}\cdot a^j(x,u,\pd u)
\d{\bf x}+\int_S\psi\cdot b(x,u,\pd u)\d{\bf x}=0
\label{us2eq4}
\e
for all $\psi\in C^1(S)$ supported in $S^\circ$. We get \eq{us2eq4}
by multiplying \eq{us2eq3} by $\psi$ and integrating by parts.
Note that \eq{us2eq4} makes sense even if $u$ is only once weakly
differentiable.
\end{itemize}
Clearly, the first sense implies the second implies the third.
But if $Q$ is {\it elliptic}, under suitable assumptions on
$a^j,b$ one can show that a weak solution $u$ is a solution,
so all three senses coincide. See for instance~\cite[\S 8]{GiTr}.

\subsection{Continuity of functions with $L^p,L^q$ derivatives}
\label{us22}

In \S\ref{us6} and \S\ref{us7} we will need the following
result, to prove that $u$ is continuous in weak solutions
$u,v$ of \eq{us1eq3} when~$a=0$.

\begin{thm} Let\/ $S$ be a domain in $\R^2$, and\/ $p,q\in(1,\iy)$
with\/ $p^{-1}+q^{-1}<1$. Then there exist continuous functions
$G,H:S\t S\ra[0,\iy)$ depending only on $S,p$ and\/ $q$ satisfying
\e
\begin{gathered}
G(w',x',w,x)=G(w,x,w',x'),\quad H(w',x',w,x)=H(w,x,w',x'),\\
\text{and}\quad
G(w,x,w,x)=H(w,x,w,x)=0\quad\text{for all\/ $(w,x),(w',x')\in S$,}
\end{gathered}
\label{us2eq5}
\e
such that whenever $u\in C^1(S)$ then for all\/ $(w,x),(w',x')\in S$
we have
\e
\bmd{u(w,x)-u(w',x')}\le\blnm{\ts\frac{\pd u}{\pd x}}{p}
G(w,x,w',x')+\blnm{\ts\frac{\pd u}{\pd y}}{q}H(w,x,w',x').
\label{us2eq6}
\e
\label{us2thm}
\end{thm}

Before we prove the theorem, we explain its significance. The
{\it Sobolev Embedding Theorem} \cite[\S 2.3]{Aubi} implies
that if $S$ is a domain in $\R^2$ then $L^p_1(S)\hookra C^0(S)$
is a continuous inclusion when $p>2$. That is, functions in a
bounded subset of $L^p_1(S)$ are all uniformly continuous. Theorem
\ref{us2thm} is essentially equivalent to this when~$p=q$.

However, when $p\ne q$ and $p\le 2$ or $q\le 2$, our result is more
general than this consequence of the Sobolev Embedding Theorem. The
basic idea is that by taking a stronger norm of $\frac{\pd u}{\pd x}$
we can make do with a weaker norm of $\frac{\pd u}{\pd y}$, or
vice versa, and still prove that $u$ is continuous.

In fact one can prove a stronger result: the functions $G,H$ of
Theorem \ref{us2thm} actually satisfy $G(w,x,w',x')\le
C\bigl(\md{w-w'}^\ga+\md{x-x'}^\de\bigr)$ and $H(w,x,w',x')\le
C'\bigl(\md{w-w'}^\ep+\md{x-x'}^\ze\bigr)$ for $\ga,\de>0$ depending
on $p,\be$ and $\ep,\ze>0$ depending on $q,\be$ and $C,C'>0$. One can
then use this to prove that $u$ is {\it uniformly H\"older
continuous}, rather than just uniformly continuous. But we will
not need this.

We will prove the theorem explicitly in the case that $S$ is a
rectangle $R=[k,l]\t[m,n]$ in $\R^2$, and then explain briefly how
to extend the proof to general $S$. We begin with a convolution
formula for functions on~$R$.

\begin{prop} Let\/ $R$ be the closed rectangle $[k,l]\t[m,n]$ in $\R^2$,
with\/ $k<l$ and\/ $m<n$. Fix $\be>0$, and define $E:R\t R\ra[0,1)$
and\/ $F:R\t R\ra(0,\iy]$~by
\begin{gather}
E(w,x,y,z)\!=\!\begin{cases}
\max\bigl((y\!-\!w)^\be(l\!-\!w)^{-\be},(z\!-\!x)(n\!-\!x)^{-1}\bigr),
& y\!\ge\!w,\; z\!\ge\!x, \\
\max\bigl((y\!-\!w)^\be(l\!-\!w)^{-\be},(x\!-\!z)(x\!-\!m)^{-1}\bigr),
& y\!\ge\!w,\; z\!<\!x, \\
\max\bigl((w\!-\!y)^\be(w\!-\!k)^{-\be},(z\!-\!x)(n\!-\!x)^{-1}\bigr),
& y\!<\!w,\; z\!\ge\!x, \\
\max\bigl((w\!-\!y)^\be(w\!-\!k)^{-\be},(x\!-\!z)(x\!-\!m)^{-1}\bigr),
& y\!<\!w,\; z\!<\!x,\end{cases}
\nonumber\\
\text{and}\/\quad
F(w,x,y,z)=-\,\frac{E(w,x,y,z)^{-1-1/\be}-1}{(\be+1)(l-k)(m-n)}.
\label{us2eq7}
\end{gather}
Then for all\/ $u\in C^1(R)$ and $(w,x)\in R$ we have
\e
\begin{split}
u(w,x)=\frac{1}{(l-k)(n-m)}&\int_m^n\!\int_k^l\!u(y,z)\d y\,\d z\\
+&\int_m^n\!\int_k^l\!\frac{\pd u}{\pd x}(y,z)(y-w)F(w,x,y,z)\d y\,\d z\\
+\be&\int_m^n\!\int_k^l\!\frac{\pd u}{\pd y}(y,z)(z-x)F(w,x,y,z)\d y\,\d z.
\end{split}
\label{us2eq8}
\e
\label{us2prop1}
\end{prop}

\begin{proof} Let $(w,x)$ and $(r,s)$ lie in $R$. Then
\begin{align*}
u(w,x)&=u(r,s)-\int_0^1\frac{\d}{\d t}\Bigl(
u\bigl(w+t^{1/\be}(r-w),x+t(s-x)\bigr)\Bigr)\d t\\
&=u(r,s)-\int_0^1\Bigl(\be^{-1}(r-w)t^{1/\be-1}
\frac{\pd u}{\pd x}\bigl(w+t^{1/\be}(r-w),x+t(s-x)\bigr)\\
&\qquad\qquad\qquad\qquad\qquad\;
+(s-x)\frac{\pd u}{\pd y}\bigl(w+t^{1/\be}(r-w),x+t(s-x)\bigr)\Bigr)\d t
\end{align*}
Note that if $(w,x),(r,s)\in R$ then $\bigl(w+t^{1/\be}(r-w),x+t(s-x)
\bigr)\in R$ for all $t\in[0,1]$, as $R$ is a rectangle. Integrate this
equation over $(r,s)\in R$, regarding $(w,x)$ as fixed. We get
\begin{align*}
&(l-k)(n-m)u(w,x)=\int_m^n\!\int_k^l\!u(r,s)\d r\,\d s\\
&-\int_m^n\!\int_k^l\!\int_0^1\Bigl(\be^{-1}t^{1/\be-1}(r-w)
\frac{\pd u}{\pd x}\bigl(w+t^{1/\be}(r-w),x+t(s-x)\bigr)\\
&\qquad\qquad\qquad\qquad\quad\;\;
+(s-x)\frac{\pd u}{\pd y}\bigl(w+t^{1/\be}(r-w),x+t(s-x)\bigr)\Bigr)
\d t\,\d r\,\d s.
\end{align*}
Change variables from $(r,s,t)$ to $(y,z,t)=\bigl(w+t^{1/\be}(r-w),
x+t(s-x),t\bigr)$ in the triple integral. Then $\d t\,\d r\,\d s=
t^{-1-1/\be}\d t\,\d y\,\d z$, and $t^{1/\be-1}(r-w)=t^{-1}(y-w)$,
and $(s-x)=t^{-1}(z-x)$. Furthermore, if $(y,z)\in R$ and $t\in(0,1]$
then the condition for $(r,s)=\bigl(w+t^{-1/\be}(y-w),x+t^{-1}(z-x)\bigr)$
to lie in $R$ is $E(w,x,y,z)<t\le 1$, because this is how we defined $E$.
Hence
\begin{align*}
&(l-k)(n-m)u(w,x)=\int_m^n\!\int_k^l\!u(y,z)\d y\,\d z\\
&-\int_m^n\!\int_k^l\!\Bigl((y-w)\frac{\pd u}{\pd x}(y,z)+\be(z-x)
\frac{\pd u}{\pd y}(y,z)\Bigr)\int_{E(w,x,y,z)}^1\be^{-1}t^{-2-1/\be}
\d t\,\d y\,\d z.
\end{align*}
Equation \eq{us2eq8} then follows by dividing by $(l-k)(n-m)$, doing
the $t$ integral explicitly, and substituting in~\eq{us2eq7}.
\end{proof}

Now $F$ is given entirely explicitly in \eq{us2eq7}, so the
following is an exercise in Lebesgue integration, which we leave
to the reader.

\begin{prop} In the situation of Proposition \ref{us2prop1}, for
$(w,x)\in R$ the function $(y,z)\mapsto(y-w)F(w,x,y,z)$ lies in $L^s(R)$
if\/ $1\le s<1+\be^{-1}$, and then the map $R\ra L^s(R)$ taking $(w,x)$
to the function $(y,z)\mapsto(y-w)F(w,x,y,z)$ is continuous. Similarly, 
the function $(y,z)\mapsto(z-x)F(w,x,y,z)$ lies in $L^t(R)$ if\/
$1\le t<1+\be$, and then the corresponding map $R\ra L^t(R)$ is
continuous.
\label{us2prop2}
\end{prop}

We can now define the functions $G,H$ in Theorem~\ref{us2thm}.

\begin{dfn} Let $R$ be the closed rectangle $[k,l]\t[m,n]$ in $\R^2$,
with $k<l$ and $m<n$. Suppose $p,q\in(1,\iy)$ with $p^{-1}+q^{-1}<1$.
This is equivalent to $1/(q-1)<p-1$. Choose $\be$ with $1/(q-1)<\be<p-1$.
Define $s=p/(p-1)$ and $t=q/(q-1)$. Then
\begin{equation*}
\frac{1}{p}+\frac{1}{s}=1,\quad
\frac{1}{q}+\frac{1}{t}=1,\quad
1<s<1+\be^{-1}\quad\text{and}\quad
1<t<1+\be.
\end{equation*}
Let $F$ be as in Proposition \ref{us2prop1}, with this value of
$\be$. Define functions $G,H:R\t R\ra[0,\iy)$ by
\begin{align*}
G(w,x,w',x')&\!=\!\Bigl(\int_m^n\!\!\int_k^l\!\bmd{(y\!-\!w)F(w,x,y,z)
\!-\!(y\!-\!w')F(w',x',y,z)}^s\d y\,\d z\Bigr)^{1/s},\\
H(w,x,w',x')&\!=\!\be\Bigl(\int_m^n\!\!\int_k^l\!\bmd{(z\!-\!x)F(w,x,y,z)
\!-\!(z\!-\!x')F(w',x',y,z)}^t\d y\,\d z\Bigr)^{1/t}.
\end{align*}
Then $G$ is well-defined and continuous, as the functions
\begin{equation*}
(y,z)\mapsto(y-w)F(w,x,y,z) \quad\text{and}\quad
(y,z)\mapsto(y-w')F(w',x',y,z)
\end{equation*}
lie in $L^s(R)$ by Proposition \ref{us2prop2} and depend continuously
on $(w,x)$ and $(w',x')$, and $G(w,x,w',x')$ is the $L^s$ norm of their
difference. Similarly, $H$ is well-defined and continuous.
\label{us2def3}
\end{dfn}

\begin{prop} Theorem \ref{us2thm} holds when $S$ is a
rectangle~$[k,l]\t[m,n]$.
\label{us2prop3}
\end{prop}

\begin{proof} Let $G,H$ be as in Definition \ref{us2def3}. Then
$G,H$ are continuous, and \eq{us2eq5} is immediate from the
definition. Subtracting equation \eq{us2eq8} with values
$(w,x)$ and $(w',x')$ gives
\begin{align*}
u(w&,x)-u(w,x)=\\
&\int_m^n\int_k^l\frac{\pd u}{\pd x}(y,z)
\bigl((y-w)F(w,x,y,z)-(y-w')F(w',x',y,z)\bigr)\d y\,\d z\\
+\be&\int_m^n\int_k^l\frac{\pd u}{\pd y}(y,z)
\bigl((z-x)F(w,x,y,z)-(z-x')F(w',x',y,z)\bigr)\d y\,\d z.
\end{align*}
Equation \eq{us2eq6} follows by using H\"older's inequality
to estimate the first integral in terms of the $L^p$ norm of
$\frac{\pd u}{\pd x}$ and the $L^s$ norm of the other factor,
and the second integral in terms of the $L^q$ norm of
$\frac{\pd u}{\pd y}$ and the $L^t$ norm of the other factor.
\end{proof}

To extend Theorem \ref{us2thm} to general domains $S$ in $\R^2$, we
mimic the proof of Proposition \ref{us2prop1} to derive an analogue
of \eq{us2eq8} that holds on $S$, not just on a rectangle $R$. The
main problem in doing this is that if $(w,x),(r,s)\in S$ then
$\bigl(w+t^{1/\be}(r-w),x+t(s-x)\bigr)$ may not lie in $S$ for
all $t\in[0,1]$. So we choose a more general family of paths
joining $(w,x)$ to all points $(r,s)$ in $S$. If these have the
same qualitative behaviour near $(w,x)$, scaling like $t^{1/\be}$
in the $x$ coordinate and $t$ in the $y$ coordinate, then the
analogue of Proposition \ref{us2prop2} should hold, and the rest
of the proof follows with little change.

\section{$\U(1)$-invariant special Lagrangian 3-folds}
\label{us3}

For introductions to special Lagrangian geometry, see Harvey and
Lawson \cite[\S III]{HaLa} and the author \cite{Joyc2}. Here is
the definition of special Lagrangian submanifolds in $\C^m$,
taken from~\cite[\S III]{HaLa}.

\begin{dfn} Let $\C^m$ have complex coordinates $(z_1,\dots,z_m)$, 
and define a metric $g$, a real 2-form $\om$ and a complex $m$-form 
$\Om$ on $\C^m$ by
\e
\begin{split}
g=\ms{\d z_1}+\cdots+\ms{\d z_m},\quad
\om&=\frac{i}{2}(\d z_1\w\d\bar z_1+\cdots+\d z_m\w\d\bar z_m),\\
\text{and}\quad\Om&=\d z_1\w\cdots\w\d z_m.
\end{split}
\label{us3eq1}
\e
Then $\Re\Om$ and $\Im\Om$ are real $m$-forms on $\C^m$. Let
$L$ be an oriented real submanifold of $\C^m$ of real dimension 
$m$. We say that $L$ is a {\it special Lagrangian submanifold} 
of $\C^m,$ or {\it SL\/ $m$-fold}\/ for short, if $L$ is calibrated 
with respect to~$\Re\Om$. Equivalently \cite[Cor.~III.1.11]{HaLa},
$L$ is special Lagrangian (with some orientation) if
$\om\vert_L\equiv 0$ and\/~$\Im\Om\vert_L\equiv 0$.
\label{us3def1}
\end{dfn}

We now recall a few fundamental results from \cite{Joyc3} that
will be used very often later. This paper is not designed to
be read independently of \cite{Joyc3}, and many other results
from \cite{Joyc3} will be cited when they are needed. Readers
are referred to \cite{Joyc3} for proofs, discussion and motivation.
The following result \cite[Prop.~4.1]{Joyc3} is the starting
point for everything in \cite{Joyc3,Joyc4} and this paper.

\begin{prop} Let\/ $S$ be a domain in $\R^2$ or $S=\R^2$, let\/
$u,v:S\ra\R$ be continuous, and\/ $a\in\R$. Define
\e
\begin{split}
N=\bigl\{(z_1,z_2,z_3)\in\C^3:\,& z_1z_2=v(x,y)+iy,\quad z_3=x+iu(x,y),\\
&\ms{z_1}-\ms{z_2}=2a,\quad (x,y)\in S\bigr\}.
\end{split}
\label{us3eq2}
\e
Then 
\begin{itemize}
\setlength{\itemsep}{0pt}
\setlength{\parsep}{0pt}
\item[{\rm(a)}] If\/ $a=0$, then $N$ is a (possibly singular)
special Lagrangian $3$-fold in $\C^3$, with boundary over
$\pd S$, if\/ $u,v$ are differentiable and satisfy
\e
\frac{\pd u}{\pd x}=\frac{\pd v}{\pd y}
\quad\text{and}\quad
\frac{\pd v}{\pd x}=-2\bigl(v^2+y^2\bigr)^{1/2}\frac{\pd u}{\pd y},
\label{us3eq3}
\e
except at points $(x,0)$ in $S$ with\/ $v(x,0)=0$, where $u,v$ 
need not be differentiable. The singular points of\/ $N$ are those
of the form $(0,0,z_3)$, where $z_3=x+iu(x,0)$ for $x\in\R$ 
with\/~$v(x,0)=0$.
\item[{\rm(b)}] If\/ $a\ne 0$, then $N$ is a nonsingular SL\/
$3$-fold in $\C^3$, with boundary over $\pd S$, if and only if\/
$u,v$ are differentiable on all of\/ $S$ and satisfy
\e
\frac{\pd u}{\pd x}=\frac{\pd v}{\pd y}\quad\text{and}\quad
\frac{\pd v}{\pd x}=-2\bigl(v^2+y^2+a^2\bigr)^{1/2}\frac{\pd u}{\pd y}.
\label{us3eq4}
\e
\end{itemize}
\label{us3prop1}
\end{prop}

In \cite[Prop.~7.1]{Joyc3} we show that solutions $u,v\in C^1(S)$ 
of \eq{us3eq4} are derived from a {\it potential function}
$f\in C^2(S)$ satisfying a {\it second-order quasilinear elliptic
equation}.

\begin{prop} Let\/ $S$ be a domain in $\R^2$ and\/ $u,v\in C^1(S)$
satisfy \eq{us3eq4} for $a\ne 0$. Then there exists $f\in C^2(S)$
with\/ $\frac{\pd f}{\pd y}=u$, $\frac{\pd f}{\pd x}=v$ and
\e
P(f)=\Bigl(\Bigl(\frac{\pd f}{\pd x}\Bigr)^2+y^2+a^2
\Bigr)^{-1/2}\frac{\pd^2f}{\pd x^2}+2\,\frac{\pd^2f}{\pd y^2}=0.
\label{us3eq5}
\e
This $f$ is unique up to addition of a constant, $f\mapsto f+c$.
Conversely, all solutions of\/ \eq{us3eq5} yield solutions 
of\/~\eq{us3eq4}. 
\label{us3prop2}
\end{prop}

Equation \eq{us3eq5} may also be written in {\it divergence form} as
\begin{align}
P(f)&=\frac{\pd}{\pd x}\Bigl[A\Bigl(a,y,\frac{\pd f}{\pd x}\Bigr)\Bigr]
+2\,\frac{\pd^2f}{\pd y^2}=0,
\label{us3eq6}\\
\text{where}\quad
A(a,y,v)&=\int_0^v\bigl(w^2+y^2+a^2\bigr)^{-1/2}\,\d w,
\label{us3eq7}
\end{align}
so that $\frac{\pd A}{\pd v}=(v^2+y^2+a^2)^{-1/2}$. Note that
$A$ is undefined when~$a=y=0$.

In \cite[Prop.~8.1]{Joyc3} we show that if $u,v$ satisfy \eq{us3eq4}
then $v$ satisfies a {\it second-order quasilinear elliptic equation},
and conversely, any solution $v$ of this equation extends to a solution 
$u,v$ of~\eq{us3eq4}.

\begin{prop} Let\/ $S$ be a domain in $\R^2$ and\/ $u,v\in C^2(S)$
satisfy \eq{us3eq4} for $a\ne 0$. Then
\e
Q(v)=\frac{\pd}{\pd x}\Bigl[\bigl(v^2+y^2+a^2\bigr)^{-1/2}
\frac{\pd v}{\pd x}\Bigr]+2\,\frac{\pd^2v}{\pd y^2}=0.
\label{us3eq8}
\e
Conversely, if\/ $v\in C^2(S)$ satisfies \eq{us3eq8} then 
there exists $u\in C^2(S)$, unique up to addition of a 
constant\/ $u\mapsto u+c$, such that $u,v$ satisfy~\eq{us3eq4}.
\label{us3prop3}
\end{prop}

Defining $A$ as in \eq{us3eq7}, equation \eq{us3eq8} may
also be written
\e
\frac{\pd^2}{\pd x^2}\bigl(A(a,y,v)\bigr)+2\,\frac{\pd^2v}{\pd y^2}=0.
\label{us3eq9}
\e
In \cite[Th.~7.6]{Joyc3} and \cite[Th.~8.8]{Joyc3} we prove
existence and uniqueness of solutions to the Dirichlet problems for
\eq{us3eq5} and \eq{us3eq8} in (strictly convex) domains in~$\R^2$.

\begin{thm} Let\/ $S$ be a strictly convex domain in $\R^2$, 
and let\/ $a\ne 0$, $k\ge 0$ and\/ $\al\in(0,1)$. Then for 
each\/ $\phi\in C^{k+2,\al}(\pd S)$ there exists a unique 
solution $f$ of\/ \eq{us3eq5} in $C^{k+2,\al}(S)$ with\/ 
$f\vert_{\pd S}=\phi$. This $f$ is real analytic in $S^\circ$, 
and satisfies $\cnm{f}{1}\le C\cnm{\phi}{2}$, for some $C>0$
depending only on~$S$.
\label{us3thm1}
\end{thm}

\begin{thm} Let\/ $S$ be a domain in $\R^2$. Then whenever $a\ne 0$, 
$k\ge 0$, $\al\in(0,1)$ and\/ $\phi\in C^{k+2,\al}(\pd S)$ there 
exists a unique solution $v\in C^{k+2,\al}(S)$ of\/ \eq{us3eq8} 
with\/ $v\vert_{\pd S}=\phi$. Fix a basepoint\/ $(x_0,y_0)\in S$.
Then there exists a unique $u\in C^{k+2,\al}(S)$ with\/ $u(x_0,y_0)=0$
such that $u,v$ satisfy \eq{us3eq4}. Furthermore, $u,v$ are real
analytic in~$S^\circ$.
\label{us3thm2}
\end{thm}

Combined with Propositions \ref{us3prop1} and \ref{us3prop2}, these
give existence and uniqueness results for $\U(1)$-invariant special
Lagrangian 3-folds in $\C^3$ satisfying certain boundary conditions.

\section{Two families of model solutions of \eq{us3eq4}}
\label{us4}

In this section we construct two families of solutions $u,v$
and $u',v'$ of \eq{us3eq4} on the unit disc $D$ in $\R^2$, and
make detailed analytic estimates of $u,v$, $u',v'$ and their
derivatives. The reason for doing this is that in \S\ref{us5}
we will use these examples and \cite[Th.~6.9]{Joyc3} to prove
a priori estimates of the first derivatives of solutions $u,v$
of \eq{us3eq4} satisfying a bound $u^2+v^2<L^2$. These estimates
are the key technical tool we will need to extend Theorems
\ref{us3thm1} and \ref{us3thm2} to the case $a=0$, and prove
other important facts about singular solutions of~\eq{us3eq3}.

The method we shall use is to start with the exact solutions
$u=\be x$, $v=\al+\be y$ of \eq{us3eq4}, and add on a small
perturbation. This perturbation is the sum of a known, exact
solution of the linearization of \eq{us3eq4} at $u=\be x$,
$v=\al+\be y$, and an `error term'. Most of the hard work
below is in estimating this error term, and showing that in
some circumstances it is small, so that the explicit approximate
solution of \eq{us3eq4} is close to an exact solution.

\subsection{The main results}
\label{us41}

Here are the two main results of this section.

\begin{thm} Let\/ $K,L,M,N>0$ be given. Then there exist\/ $A,B,C>0$
depending only on $K,L,M,N$ such that the following is true.

Suppose $a,x_0,y_0,u_0,v_0,p_0$ and\/ $q_0$ are real numbers satisfying
\begin{equation}
\begin{aligned}{}
&\md{p_0}\!\le\!A\,\max\bigl((v_0^2\!+\!a^2)^{1/2}\md{y_0}^{1/2},
(v_0^2\!+\!a^2)^{3/4}\bigr),&\;\>
\md{y_0}\md{q_0}\!\le\!B(v_0^2\!+\!a^2)^{1/2},&\\
&\md{q_0}\le C, \quad a\ne 0, \quad \md{a}\le K, 
\quad x_0^2+y_0^2\le L^2,& \quad\text{and}\quad u_0^2+v_0^2\le M^2.&
\end{aligned}
\label{us4eq1}
\end{equation}
Define $D_L=\bigl\{(x,y)\in\R^2:x^2+y^2\le L^2\bigr\}$. Then 
there exist $u,v\in C^\iy(D_L)$ satisfying \eq{us3eq4},
$(u-u_0)^2+(v-v_0)^2<N^2$ on $D_L$, and
\e
u(x_0,y_0)=u_0,\;\> v(x_0,y_0)=v_0,\;\>
\frac{\pd v}{\pd x}(x_0,y_0)=p_0\;\>\text{and\/}\;\>
\frac{\pd v}{\pd y}(x_0,y_0)=q_0.
\label{us4eq2}
\e
\label{us4thm1}
\end{thm}

\begin{thm} Let\/ $K,L,M,N>0$. Then there exist\/ $A,B>0$ depending
only on $K,L,M,N$ such that the following is true.

Suppose $a,x_0,y_0,u_0,v_0,p_0$ and\/ $q_0$ are real numbers satisfying
\e
\begin{aligned}
\md{p_0}&\le A\,\md{y_0}^{5/2},&\quad
\md{q_0}&\le B,&\quad a&\ne 0, \\
\md{a}&\le K, &\quad x_0^2+y_0^2&\le L^2, 
\quad&\text{and}\quad u_0^2+v_0^2&\le M^2.
\end{aligned}
\label{us4eq3}
\e
Define $D_L=\bigl\{(x,y)\in\R^2:x^2+y^2\le L^2\bigr\}$. Then 
there exist $u',v'\in C^\iy(D_L)$ satisfying \eq{us3eq4},
$(u'-u_0)^2+(v'-v_0)^2<N^2$ on $D_L$, and
\e
u'(x_0,y_0)=u_0,\;\> v'(x_0,y_0)=v_0,\;\>
\frac{\pd v'}{\pd x}(x_0,y_0)=p_0\;\>\text{and\/}\;\>
\frac{\pd v'}{\pd y}(x_0,y_0)=q_0.
\label{us4eq4}
\e
\label{us4thm2}
\end{thm}

We can interpret Theorems \ref{us4thm1} and \ref{us4thm2} like this:
given $a,x_0,y_0,u_0,v_0,L$ and $N$ with $a\ne 0$ and $x_0^2+y_0^2\le
L^2$, we wish to know what are the possible values of $p_0=\frac{\pd
v}{\pd x}(x_0,y_0)$ and $q_0=\frac{\pd v}{\pd y}(x_0,y_0)$ for
solutions $u,v$ of \eq{us3eq4} on $D_L$ with $(u,v)(x_0,y_0)=(u_0,v_0)$
and $(u-u_0)^2+(v-v_0)^2<N^2$. The theorems give ranges of values of
$p_0$ and $q_0$ for which such solutions are guaranteed to exist, in
terms of upper bounds $K$ for $\md{a}$ and $M$ for $\bmd{(u_0,v_0)}$.
In \S\ref{us5} we will combine the theorems with \cite[Th.~6.9]{Joyc3}
to derive a priori estimates of the first derivatives of solutions
$u,v$ of~\eq{us3eq4}.

The only difference between Theorems \ref{us4thm1} and \ref{us4thm2}
is the conditions on $p_0$ and $q_0$ in \eq{us4eq1} and \eq{us4eq3}.
Roughly speaking, when $y_0$ is small but $v_0^2+a^2$ is not small
Theorem \ref{us4thm1} is a stronger result than Theorem \ref{us4thm2}
--- that is, the requirements on $p_0$ and $q_0$ in \eq{us4eq3} are
unnecessarily strict. And when $v_0$ and $a$ are small but $y_0$ is
not small, Theorem \ref{us4thm2} is stronger --- that is, the
requirements on $p_0$ and $q_0$ in \eq{us4eq1} are unnecessarily
strict.

In \S\ref{us5} we will use Theorems \ref{us4thm1} and \ref{us4thm2}
to prove {\it a priori estimates} for $\frac{\pd u}{\pd x},\frac{\pd
u}{\pd y},\frac{\pd v}{\pd x}$ and $\frac{\pd v}{\pd y}$ when $u,v$
satisfy \eq{us3eq4} and a $C^0$ bound. Only by using {\it both\/}
Theorems \ref{us4thm1} and \ref{us4thm2} together will we be able
to make these a priori estimates powerful enough for the applications
in \S\ref{us6}--\S\ref{us7}, in particular, proving the {\it
continuity} of singular solutions $u,v$ of~\eq{us3eq3}.

Theorems \ref{us4thm1} and \ref{us4thm2} are the only results of
this section that will be used later, so readers not interested
in the proofs should skip on to \S\ref{us5}. A sketch of the
proofs is given after Definition~\ref{us4def2}. 

\subsection{Definition of the solutions}
\label{us42}

We define the first of two families of solutions $u,v$ of
\eq{us3eq4}. In \S\ref{us47} we will use rescaled versions
of these $u,v$ to prove Theorem~\ref{us4thm1}.

\begin{dfn} Let $a,s,\al,\be$ and $\ga$ be real numbers satisfying
\e
a\ne 0, \quad 0<s\le 1,\quad \al^2+a^2=s^2,\quad \md{\be}\le 1 
\quad\text{and}\quad \md{\ga}\le\ts\frac{1}{20}s.
\label{us4eq5}
\e
Let $D$ be the closed unit disc in $\R^2$. Define functions
$F,G\in C^\iy(\R)$ and $g\in C^\iy(D)$~by 
\begin{gather}
F(y)=-\ha\int_0^y\bigl((\al+\be w)^2+w^2+a^2\bigr)^{-1/2}\d w,
\quad G(y)=\int_0^yF(w)\d w
\label{us4eq6}\\
\text{and}\quad
g(x,y)=\al x+\be xy+\ga\bigl(\ha x^2+G(y)\bigr).
\label{us4eq7}
\end{gather}
Let $f\in C^\iy(D)$ be the unique solution of \eq{us3eq5} on $D$ with
$f\vert_{\pd D}=g\vert_{\pd D}$. This exists by Theorem \ref{us3thm1}.
Define $\phi=f-g$ and $\psi=\frac{\pd\phi}{\pd x}$. Then $\phi,\psi\in 
C^\iy(D)$ and $\phi\vert_{\pd D}=0$. Define $u=\frac{\pd f}{\pd y}$ and 
$v=\frac{\pd f}{\pd x}$. Then $u,v$ satisfy \eq{us3eq4} and $v$
satisfies \eq{us3eq8}. Also, as $f=g+\phi$, from \eq{us4eq7} we get
\e
u=\be x+\ga F(y)+\frac{\pd\phi}{\pd y} \;\>\text{and}\;\>
v=\al+\be y+\ga x+\frac{\pd\phi}{\pd x}=\al+\be y+\ga x+\psi.
\label{us4eq8}
\e
\label{us4def1}
\end{dfn}

As $v=\al+\be y+\ga x+\psi$ satisfies \eq{us3eq8}, we get
\e
\begin{split}
&\bigl((\al+\be y+\ga x+\psi)^2+y^2+a^2\bigr)^{-1/2}
\frac{\pd^2\psi}{\pd x^2}+2\,\frac{\pd^2\psi}{\pd y^2}=\\
&\ts\bigl((\al+\be y+\ga x+\psi)^2+y^2+a^2\bigr)^{-3/2}
\bigl(\al+\be y+\ga x+\psi\bigr)\bigl(\frac{\pd\psi}{\pd x}+\ga\bigr)^2.
\end{split}
\label{us4eq9}
\e
This is a second-order quasilinear elliptic equation on $\psi$. We will 
use elliptic regularity results for it to estimate the derivatives 
of~$\psi$.

The function $g$ constructed above is the sum of an exact solution
$\al x+\be xy$ of \eq{us3eq5}, and a multiple $\ga$ of an exact
solution $\ha x^2+G(y)$ of the {\it linearization} of \eq{us3eq5}
at $\al x+\be xy$. Thus, when $\ga$ is small we expect $g$ to be
an approximate solution of \eq{us3eq5}, with error roughly of
order~$\ga^2$.

Here is our second family of solutions $u',v'$ of \eq{us3eq4}.
In \S\ref{us47} we will use rescaled versions of these $u',v'$
to prove Theorem~\ref{us4thm2}.

\begin{dfn} Let $a,s,\al,\be$ and $\ga$ be real numbers satisfying
\e
\ts a\ne 0, \quad 0<s\le 1,\quad \al^2+a^2=s^2,
\quad \md{\be}\le \frac{1}{40} 
\quad\text{and}\quad \md{\ga}\le\frac{1}{40}.
\label{us4eq10}
\e
Let $D$ be the closed unit disc in $\R^2$. Define functions
$F',G'\in C^\iy(\R)$ and $g'\in C^\iy(D)$~by 
\begin{gather}
F'(y)=-\ha\int_0^yw\bigl((\al+\be w)^2+w^2+a^2\bigr)^{-1/2}\d w,
\quad G'(y)=\int_0^yF'(w)\d w
\nonumber\\
\text{and}\quad
g'(x,y)=\al x+\be xy+\ga\bigl(\ha x^2y+G'(y)\bigr).
\label{us4eq11}
\end{gather}
Let $f'\in C^\iy(D)$ be the unique solution of \eq{us3eq5} on 
$D$ with $f'\vert_{\pd D}=g'\vert_{\pd D}$. This exists by
Theorem~\ref{us3thm1}.

Define $\phi'=f'-g'$ and $\psi'=\frac{\pd\phi'}{\pd x}$. Then
$\phi',\psi'\in C^\iy(D)$ and $\phi'\vert_{\pd D}=0$. Define
$u'=\frac{\pd f'}{\pd y}$ and $v'=\frac{\pd f'}{\pd x}$. Then
$u',v'$ satisfy \eq{us3eq4} and $v'$ satisfies \eq{us3eq8}.
Also, as $f'=g'+\phi'$, from \eq{us4eq11} we get
\begin{equation*}
u'=\be x+\ga F'(y)+\frac{\pd\phi'}{\pd y} \;\>\text{and}\;\>
v'=\al+\be y+\ga xy+\frac{\pd\phi'}{\pd x}=\al+\be y+\ga xy+\psi'.
\end{equation*}
\label{us4def2}
\end{dfn}

Again, we start with an approximate solution $g'(x,y)$ of
\eq{us3eq5} which is the sum of the exact solution $\al x+\be xy$
of \eq{us3eq5} and an exact solution of the linearization of
\eq{us3eq5} at $\al x+\be xy$, but this time it is a {\it
different\/} solution of the linearization to that used in
Definition~\ref{us4def1}.

Here is how the rest of the section is laid out. In
\S\ref{us43} we derive {\it a priori $C^0$ estimates\/}
for $\phi$ and $\phi'$. We do this by noting that
$f,f'$ satisfy \eq{us3eq5}, and writing down explicit
{\it super-} and {\it subsolutions\/} $f_\pm$, $f'_\pm$ of
\eq{us3eq5} with $f_-=f=f_+$ and $f'_-=f'=f'_+$ on $\pd D$,
so that $f_-\le f\le f_+$ and $f'_-\le f'\le f'_+$ on $D$.
Similarly, \S\ref{us44} derives a priori $C^0$ estimates
for $\psi,\psi'$, by noting that $v,v'$ satisfy \eq{us3eq8},
and writing down explicit super- and subsolutions
$v_\pm,v'_\pm$ of~\eq{us3eq8}.

Section \ref{us45} extends these to a priori estimates of {\it
higher derivatives\/} of $\psi,\psi'$, using standard interior
elliptic regularity results and a scaling argument. Section
\ref{us46} shows that if $a,v_0,p_0,q_0$ satisfy some
inequalities then in Definition \ref{us4def1} we can choose
$\al,\be,\ga$ such that $v(x_0,y_0)=v_0$, $\frac{\pd v}{\pd x}
(x_0,y_0)=p_0$ and $\frac{\pd v}{\pd y}(x_0,y_0)=q_0$, and
similarly for $u',v'$. Finally, \S\ref{us47} completes the proofs.

\subsection{A priori $C^0$ estimates for $\phi$ and $\phi'$}
\label{us43}

Tedious manipulation of inequalities and completing the square proves:

\begin{lem} In Definition \ref{us4def1} or \ref{us4def2},
whenever $\md{w}\le\frac{1}{10}s$ and\/ $\md{y}\le 1$ we have
\e
\ts\frac{1}{4}(y^2+s^2)\le (\al+\be y+w)^2+y^2+a^2\le 4(y^2+s^2).
\label{us4eq12}
\e
\label{us4lem}
\end{lem}

Here are explicit {\it super}- and {\it subsolutions} $f_\pm$ of
\eq{us3eq5} for $f$ in~$D$.

\begin{prop} In the situation of Definition \ref{us4def1}, define
functions $f_\pm\in C^\iy(D)$ by $f_\pm=g\pm 8s^{-1}\ga^2(1\!-\!
x^2\!-\!y^2)$. Then $P(f_+)\le 0$ and\/ $P(f_-)\ge 0$, where
$P$ is defined in~\eq{us3eq5}.
\label{us4prop1}
\end{prop}

\begin{proof} Computation using \eq{us3eq5} and \eq{us4eq7} shows that
$P(f_\pm)=\mp A+B$, where
\begin{align*}
A&=8s^{-1}\ga^2\Bigl[2\bigl((\al+\be y+\ga x\mp 8s^{-1}\ga^2x)^2
+y^2+a^2\bigr)^{-1/2}+4\Bigr] \quad\text{and}\\
B&=\ga\Bigl[\bigl((\al\!+\!\be y\!+\!\ga x\!\mp\!8s^{-1}\ga^2x)^2\!+
\!y^2\!+\!a^2\bigr)^{-1/2}\!-\!\bigl((\al\!+\!\be y)^2\!+
\!y^2\!+\!a^2\bigr)^{-1/2}\Bigr].
\end{align*}
We show that $A\ge\md{B}$, which gives $P(f_+)\le 0$ and~$P(f_-)\ge 0$.

Define $w=\ga x\mp 8s^{-1}\ga^2x$. Then as $\md{x}\le 1$ and $\md{\ga}
\le\frac{1}{20}s$ we see that $\md{w}\le 2\md{\ga}$. Applying the Mean 
Value Theorem to $h(z)=((\al+\be y+z)^2+y^2+a^2)^{-1/2}$ 
between $z=0$ and $w$ implies that
\begin{equation*}
B=-\,\ga w(\al+\be y+z)\bigl((\al+\be y+z)^2+y^2+a^2\bigr)^{-3/2},
\end{equation*}
for some $z$ between 0 and $w$. As $\md{\al+\be y+z}\le((\al+\be y+z)^2
+y^2+a^2)^{1/2}$ and $\md{z}\le\md{w}\le 2\md{\ga}\le\frac{1}{10}s$, 
this gives
\begin{equation*}
\md{B}\le\frac{\md{\ga}\md{w}}{(\al+\be y+z)^2+y^2+a^2}
\le \frac{2\ga^2}{\frac{1}{4}(y^2+s^2)}\le 8s^{-2}\ga^2,
\end{equation*}
using the first inequality of Lemma \ref{us4lem} with $z$ in place 
of $w$. Also, the second inequality of Lemma \ref{us4lem} yields
\begin{equation*}
A\ge 8s^{-1}\ga^2\bigl((y^2+s^2)^{-1/2}+4\bigr)\ge 8s^{-2}\ga^2.
\end{equation*}
The last two equations give $A\ge\md{B}$, and the proof is complete.
\end{proof}

This shows that $f_+$ is a {\it supersolution} and $f_-$ a
{\it subsolution} of \eq{us3eq5}. As $1-x^2-y^2=0$ on $\pd D$, we
have $f_\pm=g=f$ on $\pd D$. Therefore \cite[Prop.~7.5]{Joyc3} shows
that $f_-\le f\le f_+$ on $D$. Subtracting $g$ from each side gives:

\begin{cor} In Definition \ref{us4def1}, $\md{\phi}\le 8s^{-1}
\ga^2(1-x^2-y^2)$ on $D$. Hence 
\e
\cnm{\phi}{0}\le 8s^{-1}\ga^2, \quad\text{and}\quad
\ts\bmd{\frac{\pd\phi}{\pd x}},\bmd{\frac{\pd\phi}{\pd y}}\le 
16s^{-1}\ga^2 \quad\text{on $\pd D$.}
\label{us4eq13}
\e
\label{us4cor1}
\end{cor}

The first inequality is the a priori $C^0$ estimate for $\phi$ that 
we want. The analogues of Proposition \ref{us4prop1} and Corollary
\ref{us4cor1} for Definition \ref{us4def2} are:

\begin{prop} Define functions $f'_\pm\in C^\iy(D)$ by $f'_\pm=g'\pm\ga^2
(1-x^2-y^2)$. Then $P(f'_+)\le 0$ and\/ $P(f'_-)\ge 0$, where $P$ is 
defined in~\eq{us3eq5}.
\label{us4prop2}
\end{prop}

\begin{proof} Computation using \eq{us3eq5} and \eq{us4eq11} shows 
that $P(f'_\pm)=\mp A+B$, where
\begin{align*}
A&=\ga^2\Bigl[2\bigl((\al+\be y+\ga xy\mp\ga^2x)^2
+y^2+a^2\bigr)^{-1/2}+4\Bigr] \quad\text{and}\\
B&=\ga y\Bigl[\bigl((\al\!+\!\be y\!+\!\ga xy\!\mp\ga^2x)^2\!+
\!y^2\!+\!a^2\bigr)^{-1/2}\!-\!\bigl((\al\!+\!\be y)^2\!+
\!y^2\!+\!a^2\bigr)^{-1/2}\Bigr].
\end{align*}
Calculations similar to the proof of Proposition \ref{us4prop1}
using \eq{us4eq10} and \eq{us4eq12} show that $A\ge\md{B}$, which
gives $P(f'_+)\le 0$ and~$P(f'_-)\ge 0$.
\end{proof}

As for Corollary \ref{us4cor1}, we deduce an a priori $C^0$
estimate for~$\phi'$:

\begin{cor} In Definition \ref{us4def2}, $\md{\phi'}\le\ga^2
(1-x^2-y^2)$ on $D$. Hence 
\begin{equation*}
\cnm{\phi'}{0}\le\ga^2, \quad\text{and}\quad
\ts\bmd{\frac{\pd\phi'}{\pd x}},\bmd{\frac{\pd\phi'}{\pd y}}\le 
2\ga^2 \quad\text{on $\pd D$.}
\end{equation*}
\label{us4cor2}
\end{cor}

\subsection{A priori $C^0$ estimates for $\psi$ and $\psi'$}
\label{us44}

For our next three results we work in the situation of 
Definition \ref{us4def1}. As $\psi=\frac{\pd\phi}{\pd x}$,
the second inequality of \eq{us4eq13} shows that $\md{\psi}
\le 16s^{-1}\ga^2$ on $\pd D$. We shall use this to construct
{\it super}- and {\it subsolutions} $v_\pm$ for $v$, and hence
derive an {\it a priori $C^0$ estimate} for~$\psi$.

\begin{prop} Define functions $v_\pm\in C^\iy(D)$ by
\e
v_\pm=\al+\be y+\ga x\pm s^{-1}\ga^2(20-4x^2-4y^2). 
\label{us4eq16}
\e
Then $Q(v_+)\le 0$ and\/ $Q(v_-)\ge 0$, where $Q$ is defined in~\eq{us3eq8}.
\label{us4prop3}
\end{prop}

\begin{proof} From \eq{us3eq8} we find that $Q(v_\pm)=\mp A-B$, where
\begin{align}
A&=8s^{-1}\ga^2\Bigl[\bigl((\al\!+\!\be y\!+\!\ga x\pm s^{-1}
\ga^2(20\!-\!4x^2\!-\!4y^2))^2\!+\!y^2\!+\!a^2\bigr)^{-1/2}\!+\!2\Bigr],
\label{us4eq17}\\
B&=\frac{\bigl(\ga\mp 8s^{-1}\ga^2x\bigr)^2
\bigl(\al+\be y+\ga x\pm s^{-1}\ga^2(20-4x^2-4y^2)\bigr)}{
\bigl((\al+\be y+\ga x\pm s^{-1}\ga^2(20-4x^2-4y^2))^2
+y^2+a^2\bigr)^{3/2}}.
\label{us4eq18}
\end{align}
We shall show that $A\ge\md{B}$, which gives $Q(v_+)\le 0$
and~$Q(v_-)\ge 0$.

Define $w=\ga x\pm s^{-1}\ga^2(20-4x^2-4y^2)$. Then as $\md{x}\le 1$, 
$x^2+y^2\le 1$ and $\ga\le\frac{1}{20}s$ we see that $\md{w}\le 
\frac{1}{10}s$. Thus Lemma \ref{us4lem} implies that
\begin{equation*}
(\al+\be y+\ga x\pm s^{-1}\ga^2(20-4x^2-4y^2))^2+y^2+a^2\ge
\ts\frac{1}{4}(y^2+s^2)\ge\frac{1}{4}s^2.
\end{equation*}
Raising this to the power $-\ha$ and using \eq{us4eq18} and
the inequalities
\begin{align*}
\bmd{\al+\be y+&\ga x\pm s^{-1}\ga^2(20-4x^2-4y^2)}\le\\
&\bigl((\al+\be y+\ga x\pm s^{-1}\ga^2(20-4x^2-4y^2))^2+y^2+a^2\bigr)^{1/2}
\end{align*}
and $\md{\ga\mp 8s^{-1}\ga^2x}\le 2\md{\ga}$ implies that
\begin{equation*}
\md{B}\le 4\ga^2\cdot 2s^{-1}\cdot\bigl((\al+\be y+\ga x\pm 
s^{-1}\ga^2(20-4x^2-4y^2))^2+y^2+a^2\bigr)^{-1/2}.
\end{equation*}
Comparing this with \eq{us4eq17} gives $A\ge\md{B}$, as we have to prove.
\end{proof}

Now $v=\al+\be x+\ga y+\frac{\pd\phi}{\pd x}$ and 
$\md{\frac{\pd\phi}{\pd x}}\le 16s^{-1}\ga^2$ on $\pd D$ by Corollary
\ref{us4cor1}. Thus \eq{us4eq16} gives $v_-\le v\le v_+$ on $\pd D$. 
But as $v$ satisfies \eq{us3eq8}, and $Q(v_+)\le 0$ and $Q(v_-)\ge 0$ 
from above, \cite[Prop.~8.5]{Joyc3} implies that $v_-\le v\le v_+$.
Subtracting $\al+\be y+\ga x$ from each side then yields

\begin{cor} In Definition \ref{us4def1}, $\md{\frac{\pd\phi}{\pd x}}
=\md{\psi}\le s^{-1}\ga^2(20-4x^2-4y^2)$ on $D$. 
Hence~$\cnm{\psi}{0}\le 20s^{-1}\ga^2$.
\label{us4cor3}
\end{cor}

As $\md{\ga}\le\frac{1}{20}s$ this gives $\md{\psi}\le\md{\ga}$, so 
that $\md{\ga x+\psi}\le 2\md{\ga}\le\frac{1}{10}s$ as $\md{x}\le 1$. 
Thus, applying Lemma \ref{us4lem} with $w=\ga x+\psi$ shows that

\begin{cor} $\frac{1}{4}(y^2+s^2)\le(\al+\be y+\ga x+\psi)^2
+y^2+a^2\le 4(y^2+s^2)$ on~$D$.
\label{us4cor4}
\end{cor}

This gives an a priori bound on the factor in front of 
$\frac{\pd^2\psi}{\pd x^2}$ in \eq{us4eq9}, showing that
\eq{us4eq9} is {\it uniformly elliptic}. The analogues of
Proposition \ref{us4prop3} and Corollary \ref{us4cor3} for
Definition \ref{us4def2} are:

\begin{prop} Define $v'_\pm\in C^\iy(D)$ by $v'_\pm=\al+\be y+\ga xy
\pm\ga^2(\frac{9}{4}-\frac{1}{4}y^2)$. Then $Q(v'_+)\le 0$ and\/ 
$Q(v'_-)\ge 0$, where $Q$ is defined in~\eq{us3eq8}.
\label{us4prop4}
\end{prop}

\begin{proof} From \eq{us3eq8} we find that $Q(v'_\pm)=\mp A-B$,
where $A=\ga^2$ and
\begin{equation*}
B=\frac{\ga^2y^2
\bigl(\al+\be y+\ga xy\pm\ga^2(\frac{9}{4}-\frac{1}{4}y^2)
\bigr)}{\bigl((\al+\be y+\ga xy\pm\ga^2(\frac{9}{4}-\frac{1}{4}y^2))^2
+y^2+a^2\bigr)^{3/2}}.
\end{equation*}
As $\md{y},\md{\al+\be y+\ga xy\pm\ga^2(\frac{9}{4}-\frac{1}{4}y^2)}
\le((\al+\be y+\ga xy\pm\ga^2(\frac{9}{4}-\frac{1}{4}y^2))^2+y^2
+a^2)^{1/2}$ we see that $\md{B}\le\ga^2=A$, and so $Q(v'_+)\le 0$ 
and~$Q(v'_-)\ge 0$.
\end{proof}

\begin{cor} In Definition \ref{us4def2},~$\cnm{\frac{\pd\phi'}{\pd
x}}{0}=\cnm{\psi'}{0}\le\frac{9}{4}\ga^2$.
\label{us4cor5}
\end{cor}

\subsection{Estimates for higher derivatives of $\psi$ and $\psi'$}
\label{us45}

For the next part of the proof we make use of an {\it overall
symmetry} of equations \eq{us3eq2}--\eq{us3eq9}, arising
geometrically from {\it dilations} on $\C^3$. In the situation
of Proposition \ref{us3prop1}, suppose $a\ne 0$ and $u,v$ satisfy
\eq{us3eq4}, so that $N$ defined in \eq{us3eq2} is an SL 3-fold.
Let $r>0$, and define $\ti N=r^{-1}N$. Since $\ti N$ is a dilation
of $N$ it is also an SL 3-fold, and also $\U(1)$-invariant.

Now $\ti N$ is of the form \eq{us3eq2} with $u,v$ and $a$
replaced by
\begin{equation*}
\ti u(x,y)=r^{-1}u(rx,r^2y),\quad\ti v(x,y)=r^{-2}v(rx,r^2y)
\quad\text{and}\quad \ti a=r^{-2}a.
\end{equation*}
Thus $\ti u,\ti v$ and $\ti a$ also satisfy \eq{us3eq4},
and $\ti v,\ti a$ satisfy \eq{us3eq8}. If $f$ is a
potential for $u,v$ as in Proposition \ref{us3prop2},
so that $f,a$ satisfy \eq{us3eq5}, then $\ti f(x,y)=
r^{-3}f(rx,r^2y)$ is a potential for $\ti u,\ti v$,
and $\ti f,\ti a$ satisfy~\eq{us3eq5}.

The fact that $x$ scales by $r$ and $y$ by $r^2$ is one
reason why $x$ {\it and\/ $y$ derivatives behave differently
in this problem}, so that for instance in \S\ref{us65} we
obtain $L^p$ estimates for $\frac{\pd u}{\pd x}$ and $L^q$
estimates for $\frac{\pd u}{\pd y}$ with $p,q$ in different
ranges.

Here is how we use the rescaling idea in this section. The
function $\psi$ of \S\ref{us42} transforms like $v$, so its
rescaled version is $\ti\psi(x,y)=r^{-2}\psi(rx,r^2y)$,
which satisfies a rescaled version of \eq{us4eq9}. It
turns out that on a ball away from the $x$ axis, this
rescaled equation is {\it uniformly elliptic uniformly in}
$r\in(0,1]$. Thus we can use elliptic regularity results
to estimate $\ti\psi$ and its derivatives away from the
$x$-axis independently of $r\in(0,1]$. Transforming back
gives estimates of $\psi$ and its derivatives near the
$x$-axis, in terms of powers of~$y$.

We shall work in the situation of Definition \ref{us4def1} until
Theorem~\ref{us4thm3}.

\begin{dfn} In Definition \ref{us4def1}, let $\sqrt{s}\le r\le 1$. Write
\begin{equation*}
\ti D=\bigl\{(x,y)\in\R^2:r^2x^2+r^4y^2\le 1\bigr\},
\end{equation*}
and define $\ti\psi\in C^\iy(\ti D)$ by~$\ti\psi(x,y)=r^{-2}\psi(rx,r^2y)$.
\label{us4def3}
\end{dfn}

From Corollaries \ref{us4cor3} and \ref{us4cor4} and equation
\eq{us4eq9} we deduce:

\begin{prop} In the situation above, $\ti\psi$ satisfies
$\cnm{\ti\psi}{0}\le 20r^{-2}s^{-1}\ga^2$,
\begin{align}
\begin{split}
&\Bigl(\bigl(r^{-2}\al+\be y+r^{-1}\ga x+\ti\psi\,\bigr)^2+y^2+r^{-4}a^2
\Bigr)^{-1/2}\frac{\pd^2\ti\psi}{\pd x^2}+2\,\frac{\pd^2\ti\psi}{\pd y^2}=\\
&\frac{\bigl(r^{-2}\al+\be y+r^{-1}\ga x+\ti\psi
\bigr)\bigl(\frac{\pd\ti\psi}{\pd x}+r^{-1}\ga\bigr)^2}{
\bigl((r^{-2}\al+\be y+r^{-1}\ga x+\ti\psi)^2+y^2+
r^{-4}a^2\bigr)^{3/2}},\quad\text{and}
\end{split}
\label{us4eq19}\\
&\ts\frac{1}{4}\bigl(y^2\!+\!r^{-4}s^2\bigr)\le
\bigl(r^{-2}\al\!+\!\be y\!+\!r^{-1}\ga x\!+\!\ti\psi\,\bigr)^2
\!+\!y^2\!+\!r^{-4}a^2\le 4\bigl(y^2\!+\!r^{-4}s^2\bigr).
\label{us4eq20}
\end{align}
\label{us4prop5}
\end{prop}

Now \eq{us4eq19} involves constants $r^{-2}\al$, $\be$, $r^{-1}\ga$ 
and $r^{-2}a$ which are all uniformly bounded independently of $r$, as
\begin{equation*}
\md{r^{-2}\al}\le 1, \quad \md{\be}\le 1,\quad
\md{r^{-1}\ga}\le\frac{1}{20} \quad\text{and}\quad \md{r^{-2}a}\le 1.
\end{equation*}
Because of this we will be able to treat \eq{us4eq19} as a quasilinear 
elliptic equation satisfying elliptic regularity bounds that are 
{\it independent\/} of $a,\al,\be,\ga,r$ and $s$, no matter how small 
$r$ and $s$ are.

In the next few results we shall construct a priori bounds for the 
derivatives of $\ti\psi$ in the interior of $\ti D$. We do this by 
using interior elliptic regularity results to bound H\"older norms 
of $\ti\phi$ and $\ti\psi$ on a series of small balls $E_n$ in~$\ti D$.

\begin{dfn} Let $\ovB_R(x,y)$ be the closed ball of radius $R$ 
about $(x,y)$ in $\R^2$. Suppose $x'\in\R$ and $y'=\pm 1$ are such 
that $\ovB_{1/2}(x',y')\subset\ti D$. Define a decreasing series of 
balls $E_2\supset E_3\supset\cdots$ in $\ti D$ by $E_n=\ovB_{1/n}
(x',y')$. Fix~$\ep\in(0,1)$.
\label{us4def4}
\end{dfn}

\begin{prop} There exists $A_n>0$ independent of\/ $a,\al,\be,\ga,r,s$
such that\/ $\bcnm{\ti\psi\vert_{E_{n+3}}}{n,\ep}\le A_nr^{-2}s^{-1}\ga^2$
for each\/~$n=0,1,2,3,\ldots$.
\label{us4prop6}
\end{prop}

\begin{proof} Let us rewrite equation \eq{us4eq19} in the form
\begin{align*}
\frac{\pd}{\pd x}\Bigl[\bigl((r^{-2}\al+\be y+r^{-1}\ga x+\ti\psi\,)^2
+y^2+r^{-4}a^2\bigr)^{-1/2}\frac{\pd\ti\psi}{\pd x}\Bigr]
+\frac{\pd}{\pd y}\Bigl[2\,\frac{\pd\ti\psi}{\pd y}\Bigr]&\\
-\frac{r^{-1}\ga\bigl(r^{-2}\al+\be y+r^{-1}\ga x+\ti\psi
\bigr)}{\bigl((r^{-2}\al+\be y+r^{-1}\ga x+\ti\psi)^2+y^2+
r^{-4}a^2\bigr)^{3/2}}\,\frac{\pd\ti\psi}{\pd x}=&\\
\frac{r^{-2}\ga^2}{\bigl((r^{-2}\al+\be y+r^{-1}\ga x+\ti\psi)^2
+y^2+r^{-4}a^2\bigr)^{3/2}}.&
\end{align*}
Regard this as a linear elliptic equation $Pu=f$ on $\ti\psi$, where
\begin{equation*}
\bigl(Pu\bigr)(x)=\sum_{i,j=1}^2\frac{\pd}{\pd x_i}\Bigl[a^{ij}(x)
\frac{\pd u}{\pd x_j}\Bigr]+\sum_{i=1}^2b^i(x)\frac{\pd u}{\pd x_i}+c(x)u(x),
\end{equation*}
so that $u=\ti\psi$, $f=r^{-2}\ga^2\bigl((r^{-2}\al+\be y+r^{-1}\ga x
+\ti\psi)^2+y^2+r^{-4}a^2\bigr)^{-3/2}$,
\begin{gather}
\begin{gathered}
a^{11}=\bigl(\bigl(r^{-2}\al+\be y+r^{-1}\ga x+
\ti\psi\bigr)^2+y^2+r^{-4}a^2\bigr)^{-1/2},\\
a^{12}=a^{21}=0,\quad a^{22}=2,
\end{gathered}
\label{us4eq23}\\
b^1=-\frac{r^{-1}\ga\bigl(r^{-2}\al+\be y+r^{-1}\ga x+\ti\psi
\bigr)}{\bigl((r^{-2}\al\!+\!\be y\!+\!r^{-1}\ga x\!+\!\ti\psi)^2
\!+\!y^2\!+\!r^{-4}a^2\bigr)^{3/2}},
\;\>\text{and}\;\> b^2=c=0.
\label{us4eq24}
\end{gather}
To prove the proposition we shall apply three interior elliptic
regularity results of Gilbarg and Trudinger \cite{GiTr} for second
order linear elliptic operators with principal part in divergence
form, to prove the cases $n=0$, $n=1$ and $n\ge 2$ respectively.

From \cite[Th.~8.24, p.~202]{GiTr}, if $\la,\La_0,\nu_0>0$ and
$\ep\in(0,1)$ with
\begin{align}
\ts\sum_{i,j=1}^2a^{ij}\xi_i\xi_j\ge\la\ms{\xi} 
\quad\text{in $E_2$ for all $\xi=(\xi_1,\xi_2)\in\R^2$,}
\label{us4eq25}\\
\bcnm{a^{ij}\vert_{E_2}}{0}\le\La_0,\quad
\bcnm{b^i\vert_{E_2}}{0}\le\nu_0\quad\text{for all $i,j$, and}\quad
\bcnm{c\vert_{E_2}}{0}\le\nu_0,
\label{us4eq26}
\end{align}
then there exists $C>0$ depending only on $E_2,E_3,\la,\La_0,\nu_0$
and $\ep$ such that if $u\in C^2(E_2)$ and $Pu=f$, then~$\cnm{u
\vert_{E_3}}{0,\ep}\le C\bigl(\cnm{u}{0}+\cnm{f}{0}\bigr)$.

As $\ha\le\md{y}\le\frac{3}{2}$ for $(x,y)\in E_2$ and $0\le r^{-4}s^2
\le 1$, by \eq{us4eq20} we have
\e
\ts\frac{1}{16}\le
\bigl(r^{-2}\al\!+\!\be y\!+\!r^{-1}\ga x\!+\!\ti\psi\,\bigr)^2
\!+\!y^2\!+\!r^{-4}a^2\le 13 \quad\text{on $E_2$.}
\label{us4eq27}
\e
Using \eq{us4eq23}--\eq{us4eq24} and \eq{us4eq27} we can find
$\la,\La_0,\nu_0$ independent of $a,\al,\be,\ga,r,s$ such that
\eq{us4eq25}--\eq{us4eq26} hold.

Thus there exists $C_0>0$ independent of $a,\al,\be,\ga,r,s$ such that
\begin{align*}
\bcnm{\ti\psi\vert_{E_3}}{0,\ep}&\le C_0\bigl(\bcnm{\ti\psi\vert_{E_2}}{0}
+\bcnm{f\vert_{E_2}}{0}\bigr)\\
&\le C_0(20r^{-2}s^{-1}\ga^2+64r^{-2}\ga^2)\le A_0r^{-2}s^{-2}\ga^2,
\end{align*}
setting $A_0=84C_0$, using Corollary \ref{us4cor3}, $0<s\le 1$
and $\cnm{f\vert_{E_2}}{0}\le 64r^{-2}\ga^2$, which follows from
\eq{us4eq27} and the definition of $f$. This proves the case~$n=0$.

Combining the case $n=0$ with \eq{us4eq23} and $\cnm{a^{ij}
\vert_{E_2}}{0}\le 4$ we can find $\La_1>0$ depending on $A_0$
such that $\bcnm{a^{ij}\vert_{E_3}}{0,\ep}\le\La_1$. Applying
\cite[Th.~8.32, p.~210]{GiTr} shows that there exists $C_1>0$
depending only on $E_3,E_4,\la,\La_1,\nu_0$ and $\ep$ such that
if $u\in C^2(E_3)$ and $Pu=f$, then~$\cnm{u\vert_{E_4}}{1,\ep}\le
C_1\bigl(\cnm{u}{0}\!+\!\cnm{f}{0}\bigr)$. The argument above then
shows that the case $n=1$ holds with~$A_1=84C_1$.

The case $n\ge 2$ is proved using the {\it Schauder interior
estimates} \cite[Th.s~6.2 \& 6.17]{GiTr}, by a technique known
as {\it bootstrapping}. Suppose by induction that the proposition
holds for $n=k-1$, for $k\ge 2$. Then
\begin{equation*}
\bcnm{a^{ij}\vert_{E_{k+2}}}{k-1,\ep}\le\La_k,\quad
\bcnm{b^i\vert_{E_{k+2}}}{k-2,\ep}\le\nu_k
\quad\text{and}\quad
\bcnm{c\vert_{E_{k+2}}}{k-2,\ep}\le\nu_k,
\end{equation*}
for some $\La_k,\nu_k$ depending on $k,A_{k-1}$, and all $i,j$.
Therefore \cite[Th.s~6.2 \& 6.17]{GiTr} gives $C_k>0$ depending
on $k,E_{k+2},E_{k+3},\la,\La_k,\nu_k$ such that if
$u\in C^2(E_{k+2})$, $f\in C^{k-2,\ep}(E_{k+2})$ and $Pu=f$,
then~$\cnm{u\vert_{E_{k+3}}}{k,\ep}\le C_k\bigl(\cnm{u}{0}\!+
\!\cnm{f}{k-2,\ep}\bigr)$. 

Using \eq{us4eq27}, the definition of $f$ and the case $n=k-1$
we can find $C_k'>0$ depending only on $k,\ep$ and $A_{k-1}$
such that $\bcnm{f\vert_{E_{k+2}}}{k-2,\ep}\le C_k'r^{-2}\ga^2$.
The case $n=k$ of the proposition then follows with
$A_k=C_k(20+C_k')$. Hence by induction, the proof is complete.
\end{proof}

As $\bmd{\pd^n\ti\psi(x',y')}\le\cnm{\ti\psi\vert_{E_{n+3}}}{n,\ep}$,
we deduce:

\begin{cor} There exist constants $A_n$ for $n\ge 0$ independent 
of\/ $a,\al,\be,\ga,r$ and\/ $s$, such that if\/ $x'\in\R$ and\/ 
$y'=\pm 1$ with\/ $\ovB_{1/2}(x',y')\subset\ti D$, then
\e
\Bigl\vert\frac{\pd^{j+k}\ti\psi}{\pd x^j\pd y^k}(x',y')\Bigr\vert
\le A_{j+k}r^{-2}s^{-1}\ga^2 \quad\text{for all\/ $j,k\ge 0$.}
\label{us4eq30}
\e
\label{us4cor6}
\end{cor}

Here the $A_n$ are given in Proposition \ref{us4prop6}. Now the only
reason for setting $y'=\pm 1$ above was to be able to prove \eq{us4eq27},
using \eq{us4eq20} and the inequality $\ha\le\md{y}\le\frac{3}{2}$ on
$E_2$. In the special case $r=\sqrt{s}$, the terms $r^{-4}s^2$ in
\eq{us4eq20} are 1, and we then only need $\md{y}\le\frac{3}{2}$ on
$E_2$ to prove \eq{us4eq27}. Thus, when $r=\sqrt{s}$ the proofs above
are valid for $\md{y'}\le 1$ rather than just $\md{y'}=1$, and we get:

\begin{cor} Let\/ $r=\sqrt{s}$ in Definition \ref{us4def3}. If\/ $x'\in\R$ 
and\/ $\md{y'}\le 1$ with\/ $\ovB_{1/2}(x',y')\subset\ti D$, then
\e
\Bigl\vert\frac{\pd^{j+k}\ti\psi}{\pd x^j\pd y^k}(x',y')\Bigr\vert
\le A_{j+k}s^{-2}\ga^2 \quad\text{for all\/ $j,k\ge 0$,}
\label{us4eq31}
\e
where the $A_n$ are as in Corollary~\ref{us4cor6}.
\label{us4cor7}
\end{cor}

We shall use the last two corollaries to estimate the derivatives 
of $\psi$. As by definition $\ti\psi(x,y)=r^{-2}\psi(rx,r^2y)$, we 
see that for~$(x,y)\in D$,
\begin{equation*}
\frac{\pd^{j+k}\psi}{\pd x^j\pd y^k}(x,y)=
r^{2-j-2k}\frac{\pd^{j+k}\ti\psi}{\pd x^j\pd y^k}(x',y'),
\quad\text{where $(x',y')=(r^{-1}x,r^{-2}y)$.}
\end{equation*}
Put $r=\md{y}^{1/2}$ if $\md{y}\ge s$ and $r=\sqrt{s}$ if $\md{y}<s$.
It turns out that $x^2+y^2\le\frac{1}{4}$ is sufficient to ensure that
$\ovB_{1/2}(x',y')\subset\ti D$. Thus, Corollaries \ref{us4cor6} and 
\ref{us4cor7} give

\begin{thm} Suppose $(x,y)\in D$ with\/ $x^2+y^2\le\frac{1}{4}$. Then 
for all\/~$j,k\ge 0$,
\e
\Bigl\vert\frac{\pd^{j+k}\psi}{\pd x^j\pd y^k}(x,y)\Bigr\vert\le 
\begin{cases} A_{j+k}\md{y}^{-j/2-k}s^{-1}\ga^2, & \md{y}\ge s, \\
A_{j+k}s^{-1-j/2-k}\ga^2, & \md{y}<s, \end{cases}
\label{us4eq32}
\e
where $A_0,A_1,\ldots$ are positive constants independent of\/
$x,y,a,\al,\be,\ga$ and~$s$.
\label{us4thm3}
\end{thm}

This is a measure of how closely the solution $v$ of \eq{us3eq8}
defined in Definition \ref{us4def1} approximates $\al+\be y+\ga x$
in the interior of $D$. With some more work, the interior estimates
\eq{us4eq32} can be extended to the whole of $D$, using elliptic
regularity results for regions with boundary, but we will not need this.

Extending the material above to the functions of Definition
\ref{us4def2} is straightforward, so we just indicate the main
differences. In Definition \ref{us4def3} we replace
the assumption $\sqrt{s}\le r\le 1$ by $\max\bigl(
\md{\ga},\sqrt{s}\,\bigr)\le r\le 1$. We then follow through
to Corollary \ref{us4cor6} with little essential change.

To prove an analogue of Corollary \ref{us4cor7} we need an additional
assumption, that $\md{\ga}\le\frac{1}{10}s$. This is because for
Definition \ref{us4def1}, when $r=\sqrt{s}$ the first inequality of
\eq{us4eq20} gives $\frac{1}{4}\le(r^{-2}\al+\be y+r^{-1}\ga x+
\ti\psi)^2+y^2+r^{-4}a^2$, and we use this to bound the coefficient
of $\frac{\pd^2\ti\psi}{\pd x^2}$ in \eq{us4eq19}, and so to show
that \eq{us4eq19} is uniformly elliptic. However, for Definition
\ref{us4def2} we need an extra condition like $\md{\ga}\le\frac{1}{10}s$
to show that the analogue of \eq{us4eq19} is uniformly elliptic when
$y$ is small.

The analogue of Theorem \ref{us4thm3} we get is:

\begin{thm} Let\/ $(x,y)\in D$ with\/ $x^2\!+\!y^2\le\frac{1}{4}$,
$\md{y}\ge\ga^2$ and\/ $\md{y}\ge s$. Then 
\begin{equation*}
\Bigl\vert\frac{\pd^{j+k}\psi'}{\pd x^j\pd y^k}(x,y)\Bigr\vert\le 
A_{j+k}\md{y}^{-j/2-k}\ga^2
\quad\text{for all\/~$j,k\ge 0$.}
\end{equation*}
If\/ $\md{\ga}\le\frac{1}{10}s$, then whenever $(x,y)\in D$
with\/ $x^2\!+\!y^2\le\frac{1}{4}$ and\/ $\md{y}\le s$, we have
\begin{equation*}
\Bigl\vert\frac{\pd^{j+k}\psi'}{\pd x^j\pd y^k}(x,y)\Bigr\vert\le 
A_{j+k}s^{-j/2-k}\ga^2
\quad\text{for all\/~$j,k\ge 0$.}
\end{equation*}
Here $A_0,A_1,\ldots$ are positive constants independent
of\/ $x,y,a,\al,\be$ and\/~$\ga$.
\label{us4thm4}
\end{thm}

\subsection{Solutions of \eq{us3eq8} with $v,\frac{\pd v}{\pd x},
\frac{\pd v}{\pd y}$ prescribed at $(x_0,y_0)$}
\label{us46}

To prove a priori bounds for derivatives of solutions of 
\eq{us3eq4} in \S\ref{us5}, we will need to find examples 
of solutions $u,v$ of \eq{us3eq4} in $D$ such that $u,v,
\frac{\pd v}{\pd x},\frac{\pd v}{\pd y}$ take prescribed 
values $u_0,v_0,p_0,q_0$ at a given point $(x_0,y_0)$ in $D$.
As we are free to add a constant to $u$, it is enough to
consider only $v$, regarded as a solution of \eq{us3eq8}, and 
ensure that $v,\frac{\pd v}{\pd x},\frac{\pd v}{\pd y}$
take prescribed values.

\begin{thm} Let\/ $a,x_0,y_0,v_0,p_0$ and\/ $q_0$ be real numbers. Define
\e
\begin{gathered}
\al_0=v_0-x_0p_0-y_0q_0, \quad \be_0=q_0,\quad
\ga_0=p_0, \quad s_0=\sqrt{\al_0^2+a^2}\\
\text{and}\quad
\ts J=\min\bigl(\frac{1}{40},\ha A_0^{-1},\frac{1}{16}A_1^{-1}\bigr),
\end{gathered}
\label{us4eq33}
\e
where $A_0,A_1>0$ are as in Theorem \ref{us4thm3}. Suppose
\e
\ts a\ne 0,\;\>
x_0^2+y_0^2\le\frac{1}{4},\;\>
s_0^2\le\ha,\;\>
\md{\be_0}\le\frac{4}{5}\;\>
\text{and\/}\;\>
\md{\ga_0}\le J\max(s_0\md{y_0}^{1/2},s_0^{3/2}).
\label{us4eq34}
\e
Then there exist\/ $\al,\be,\ga\in\R$ satisfying \eq{us4eq5} and
\e
\ts\md{\al-\al_0}\le\frac{1}{4}s_0, \quad
\md{\be-\be_0}\le\frac{1}{4}s_0 \quad\text{and}\quad
\md{\ga-\ga_0}\le \md{\ga_0},
\label{us4eq35}
\e
such that the solution $v$ of\/ \eq{us3eq8} constructed
in Definition \ref{us4def1} using $a,\al,\be,\ga$ satisfies
\e
v(x_0,y_0)=v_0,\quad
\frac{\pd v}{\pd x}(x_0,y_0)=p_0 \quad\text{and}\quad
\frac{\pd v}{\pd y}(x_0,y_0)=q_0.
\label{us4eq36}
\e
\label{us4thm5}
\end{thm}

\begin{proof} Let $\al,\be,\ga$ satisfy \eq{us4eq35}, and define
$s=\sqrt{\al^2+a^2}$. Using $\md{\al-\al_0}\le\frac{1}{4}s_0$ one
can easily show that
\e
\ha s_0^2\le s^2\le 2s_0^2.
\label{us4eq37}
\e
As $s_0^2\le\ha$ by \eq{us4eq34} this gives $0<s\le 1$. The other
inequalities $\md{\be}\le 1$ and $\md{\ga}\le\ts\frac{1}{20}s$ in
\eq{us4eq5} also follow from \eq{us4eq34}. Thus \eq{us4eq5} holds,
and Definition \ref{us4def1} gives $v,\psi\in C^\iy(D)$ with
$v=\al+\be y+\ga x+\psi$, such that $v$ satisfies \eq{us3eq8} and
$\psi=0$ on~$\pd D$.

We shall show that there exist $\al,\be,\ga$ satisfying \eq{us4eq35}
for which $v$ satisfies \eq{us4eq36}. Using $v=\al+\be y+\ga x+\psi$
and \eq{us4eq33} we find that \eq{us4eq36} is equivalent to
\begin{align}
F_1(\al,\be,\ga)&=(\al-\al_0)+y_0(\be-\be_0)+x_0(\ga-\ga_0)+\psi(x_0,y_0)=0, 
\label{us4eq38}\\
F_2(\al,\be,\ga)&=(\be-\be_0)+\frac{\pd\psi}{\pd y}(x_0,y_0)=0,
\label{us4eq39}\\
F_3(\al,\be,\ga)&=(\ga-\ga_0)+\frac{\pd\psi}{\pd x}(x_0,y_0)=0.
\label{us4eq40}
\end{align}
Define $\al_\pm=\al_0\pm\frac{1}{4}s_0$, $\be_\pm=\be_0\pm\frac{1}{4}
s_0$ and $\ga_\pm=\ga_0\pm \md{\ga_0}$. Then \eq{us4eq35} is equivalent 
to $\al_-\le\al\le\al_+$, $\be_-\le\be\le\be_+$ and $\ga_-\le\ga\le
\ga_+$. Thus, \eq{us4eq38}--\eq{us4eq40} define functions $F_1,F_2,F_3:
[\al_-,\al_+]\t[\be_-,\be_+]\t[\ga_-,\ga_+]\ra\R$. Using 
\cite[Th.~7.7]{Joyc3} one can show that the $F_j$ are {\it
continuous} functions.

\begin{prop} Suppose $\ga_0\ne 0$. For all\/ $\al,\be,\ga$ satisfying 
\eq{us4eq35} we have
\e
\begin{gathered}
F_1(\al_-,\be,\ga)<0<F_1(\al_+,\be,\ga),\quad
F_2(\al,\be_-,\ga)<0<F_2(\al,\be_+,\ga),\\
\text{and}\quad F_3(\al,\be,\ga_-)<0<F_3(\al,\be,\ga_+).
\end{gathered}
\label{us4eq41}
\e
\label{us4prop7}
\end{prop}

\begin{proof} To prove the first pair of inequalities, we shall show 
that for $\al=\al_\pm,\be,\ga$ satisfying \eq{us4eq35}, we have
\e
\md{\al_\pm-\al_0}>\md{y_0}\md{\be-\be_0}
+\md{x_0}\md{\ga-\ga_0}+\bmd{\psi(x_0,y_0)}.
\label{us4eq42}
\e
Thus by \eq{us4eq38}, $F_1(\al_\pm,\be,\ga)$ has the same sign as
$\al_\pm-\al_0$, and the first part of \eq{us4eq41} follows. First
suppose $\md{y_0}\ge s_0$. Then from \eq{us4eq34}--\eq{us4eq35}
we have
\begin{equation*}
\ts\md{\al_\pm-\al_0}=\frac{1}{4}s_0,\quad
\md{\be-\be_0}\le\frac{1}{4}s_0\quad\text{and}\quad
\md{\ga-\ga_0}\le \md{\ga_0}\le Js_0\md{y_0}^{1/2},
\end{equation*}
and Theorem \ref{us4thm3} gives
\begin{equation*}
\bmd{\psi(x_0,y_0)}\le A_0s^{-1}\ga^2\le A_0(\ha s_0)^{-1}
\bigl(2Js_0\md{y_0}^{1/2}\bigr)=8A_0J^2s_0\md{y_0},
\end{equation*}
using $\ha s_0\le s$ and $\md{\ga}\le 2Js_0\md{y_0}^{1/2}$. Thus
\eq{us4eq42} holds if
\begin{equation*}
\ts\frac{1}{4}s_0> \frac{1}{4}s_0\md{y_0}+
Js_0\md{x_0}\md{y_0}^{1/2}+8A_0J^2s_0\md{y_0},
\end{equation*}
which follows from $s_0>0$, $\md{x_0}\le\ha$, $\md{y_0}\le\ha$,
$J\le\frac{1}{40}$ and $J\le\ha A_0^{-1}$. The other four
inequalities are proved in a similar way.
\end{proof}

The reason for supposing $\ga_0\ne 0$ is to get strict inequalities
in the third part of \eq{us4eq41}. We can now finish the proof of
Theorem \ref{us4thm5}. If $\ga_0=0$ then $\al=\al_0$, $\be=\be_0$ and
$\ga=\ga_0=0$ satisfy the conditions of the theorem, as then
$v=\al+\be y+\ga x=\al_0+\be_0 y$ is an {\it exact} solution of
\eq{us3eq8}, and $\psi\equiv 0$, so \eq{us4eq38}--\eq{us4eq40}
hold. So suppose $\ga_0\ne 0$. Write $B=[\al_-,\al_+]\t[\be_-,\be_+]\t
[\ga_-,\ga_+]$, and consider the map ${\bf F}=(F_1,F_2,F_3):B\ra\R^3$.
By Proposition \ref{us4prop7}, $\bf F$ maps $\pd B$ to $\R^3\sm\{0\}$.
Furthermore, both $\pd B$ and $\R^3\sm\{0\}$ are homotopic to
${\cal S}^2$, and one can show from the proposition that
${\bf F}_*:H_2(\pd B,\Z)\ra H_2(\R^3\sm\{0\},\Z)$ is the identity
map~$\Z\ra\Z$.

Suppose $\bf F$ maps $B\ra\R^3\sm\{0\}$. Then $F(\pd B)$ is homologous 
to zero in $\R^3\sm\{0\}$, as it bounds $F(B)$. But $F_*([\pd B])$
generates $H_2(\R^3\sm\{0\},\Z)$ and so is nonzero, a contradiction.
Thus $\bf F$ cannot map $B\ra\R^3\sm\{0\}$, and there exists 
$(\al,\be,\ga)\in B^\circ$ with ${\bf F}(\al,\be,\ga)=0$. From above
this is equivalent to \eq{us4eq36}, and so $\al,\be,\ga$ satisfy the
conditions in Theorem \ref{us4thm5}, completing the proof.
\end{proof}

Here is the analogue of Theorem \ref{us4thm5} for
Definition \ref{us4def2}, proved similarly.

\begin{thm} Let\/ $a,x_0,y_0,v_0,p_0$ and\/ $q_0$ be real numbers
with\/ $y_0\ne 0$. Define
\begin{gather*}
\al_0=v_0-2x_0p_0-y_0q_0, \quad \be_0=q_0-\frac{x_0p_0}{y_0},
\quad \ga_0=\frac{p_0}{y_0}\\
\text{and}\quad
\ts J=\min(\frac{1}{80},\frac{1}{4}A_0^{-1/2},\frac{1}{4}A_1^{-1}),
\end{gather*}
where $A_0,A_1>0$ are as in Theorem \ref{us4thm4}. Suppose
\begin{equation*}
\ts a\ne 0,\;\>
x_0^2+y_0^2\le\frac{1}{4},\;\>
\md{\al_0}\le\ha,\;\>
\md{\be_0}\le\frac{1}{80}\;\>
\text{and\/}\;\>
\md{\ga_0}\le J\,\md{y_0}^{3/2}.
\end{equation*}
Then there exist\/ $\al,\be,\ga\in\R$ satisfying \eq{us4eq10} and
\begin{equation*}
\ts\md{\al-\al_0}\le\ha, \quad \md{\be-\be_0}\le\frac{1}{80} 
\quad\text{and}\quad \md{\ga-\ga_0}\le\md{\ga_0},
\end{equation*}
such that the solution $v'$ of\/ \eq{us3eq8} constructed in Definition
\ref{us4def2} using $a,\al,\be,\ga$ satisfies
\begin{equation*}
v'(x_0,y_0)=v_0,\quad
\frac{\pd v'}{\pd x}(x_0,y_0)=p_0 \quad\text{and}\quad
\frac{\pd v'}{\pd y}(x_0,y_0)=q_0.
\end{equation*}
\label{us4thm6}
\end{thm}

\subsection{Proof of Theorems \ref{us4thm1} and \ref{us4thm2}}
\label{us47}

We now prove Theorem \ref{us4thm1}, by the rescaling method
of \S\ref{us45}. Define
\begin{gather*}
r=\min\bigl(1,(3L)^{-1},(3K)^{-1/2}\bigr),\quad
\hat a=r^2a,\quad \hat v_0=r^2v_0,\\
\hat x_0=rx_0,\quad \hat y_0=r^2y_0,\quad
\hat p_0=rp_0\quad\text{and}\quad \hat q_0=q_0.
\end{gather*}
Let $J$ be as in \eq{us4eq33}, and let $A,B,C>0$ be chosen such that
\begin{equation*}
A\le \ha Jr^2,\quad
A\le \frac{1}{8L^{3/2}},\quad
A\le \frac{1}{8L(M^2+K^2)^{1/4}},\quad
B\le{\ts\frac{1}{8}}
\quad\text{and}\quad
C\le{\ts\frac{4}{5}}.
\end{equation*}
Define $\hat\al_0,\hat\be_0,\hat\ga_0$ and $\hat s_0$ as in \eq{us4eq33},
using $\hat a,\hat v_0,\hat p_0,\hat q_0,\hat x_0$ and $\hat y_0$. We will
show that \eq{us4eq1} implies that \eq{us4eq34} holds for $\hat a,\hat x_0,
\hat y_0,\hat s_0,\hat\be_0$ and~$\hat\ga_0$.

As $r,a\ne 0$ we have $\hat a\ne 0$, and $r\le(3L)^{-1}$, $r\le 1$ and 
$x_0^2+y_0^2\le L^2$ imply that $\hat x_0^2+\hat y_0^2\le\frac{1}{4}$.
So the first two inequalities of \eq{us4eq34} hold. Now $\hat\al_0=
r^2(v_0-x_0p_0-y_0q_0)$. Using $\md{a}\le K$, $\md{x_0}\le L$,
$\md{y_0}\le L$, $\md{v_0}\le M$ and the first inequality of
\eq{us4eq1}, we find that
\begin{equation*}
\md{x_0p_0}\le A\max\bigl(L^{3/2},L(K^2+M^2)^{1/4}\bigr)
(v_0^2+a^2)^{1/2}\le\ts\frac{1}{8}(v_0^2+a^2)^{1/2},
\end{equation*}
as $A\le1/8L^{3/2}$ and~$A\le 1/8L(M^2+K^2)^{1/4}$.

Also $\md{y_0q_0}\le\frac{1}{8}(v_0^2+a^2)^{1/2}$ from the second
inequality of \eq{us4eq1} and $B\le\frac{1}{8}$. Thus
$\md{x_0p_0+y_0q_0}\le\frac{1}{4}(v_0^2+a^2)^{1/2}$, and following
the proof of \eq{us4eq37} we find that
\e
\ha r^4(v_0^2+a^2)\le \hat s_0^2=r^4(v_0-x_0p_0-y_0q_0)^2+r^4a^2
\le 2r^4(v_0^2+a^2).
\label{us4eq43}
\e
Now $r\le 1$, $r\le(3L)^{-1}$ and $\md{v_0}\le L$ imply that
$2r^4v_0^2\le\frac{2}{9}$, and $r\le(3K)^{-1/2}$, $\md{a}\le K$ yield
$2r^4a^2\le\frac{2}{9}$. Thus $\hat s_0^2\le\frac{4}{9}<\ha$, the third
inequality of \eq{us4eq34}. The fourth inequality $\md{\hat\be_0}\le
\frac{4}{5}$ follows from $\hat\be_0=q_0$, $\md{q_0}\le C$
and~$C\le\frac{4}{5}$.

The final inequality $\md{\hat\ga_0}\le J\max(\hat s_0
\md{\hat y_0}^{1/2},\hat s_0^{3/2})$ follows using \eq{us4eq43},
the first inequality of \eq{us4eq1}, and $A\le \ha Jr^2$. Thus
\eq{us4eq34} holds for $\hat a,\hat x_0,\hat y_0,\hat v_0,\hat p_0$
and $\hat q_0$. Therefore, by Theorem \ref{us4thm5}, there exist
$\hat\al,\hat\be,\hat\ga\in\R$ satisfying \eq{us4eq35} such that if
$\hat f,\hat\phi,\hat\psi,\hat u,\hat v\in C^\iy(D)$ are constructed
in Definition \ref{us4def1} using $\hat a,\hat\al,\hat\be,\hat\ga$ then
\e
\hat v(\hat x_0,\hat y_0)=\hat v_0,\quad
\frac{\pd\hat v}{\pd x}(\hat x_0,\hat y_0)=\hat p_0 
\quad\text{and}\quad \frac{\pd\hat v}{\pd y}(\hat x_0,\hat y_0)=\hat q_0.
\label{us4eq44}
\e

Define $u,v\in C^\iy(D_L)$ by
\begin{equation*}
u(x,y)=r^{-1}\hat u(rx,r^2y)-r^{-1}\hat u(\hat x_0,\hat y_0)+u_0
\quad\text{and}\quad
v(x,y)=r^{-2}\hat v(rx,r^2y).
\end{equation*}
As $r\le 1$ and $r\le(3L)^{-1}$ it
follows that if $(x,y)\in D_L$ then $(rx,r^2y)\in D$, so $u,v$ are
well-defined. Also, $u,v$ and $a$ satisfy \eq{us3eq4} as $\hat u,\hat v$
and $\hat a$ do, and $u(x_0,y_0)=u_0$ follows from the definition of
$u$, so \eq{us4eq44} implies~\eq{us4eq2}.

It remains only to show that $(u-u_0)^2+(v-v_0)^2<N^2$ on $D_L$.
This is implied by $\md{u-u_0}\le\ha N$ and $\md{v-v_0}\le\ha N$
on $D_L$, which in turn follows from
\e
\bmd{\hat u-\hat u(\hat x_0,\hat y_0)}\le\ha rN
\quad\text{and}\quad \bmd{\hat v-\hat v_0}\le\ha r^2N
\label{us4eq45}
\e
on $D$. But by \cite[Cor.~4.4]{Joyc3} the maxima and minima of
$\hat u$ and $\hat v$ are achieved on $\pd D$. Thus it is enough
for \eq{us4eq45} to hold on $\pd D$, and further the maxima of
$\hat u-\hat u(x_0,y_0),\hat v-\hat v_0$ on $\pd D$ are nonnegative,
and the minima nonpositive.

To prove this we use \eq{us4eq8} to write $\hat u,\hat v$ in terms of
$x,y,\frac{\pd\hat\phi}{\pd y}$ and $\frac{\pd\hat\phi}{\pd x}$, and
then apply \eq{us4eq13} to show that $\smash{\bmd{\frac{\pd\hat\phi}{\pd
x}},\bmd{\frac{\pd\hat\phi}{\pd y}}}\le 16\hat s^{-1}\hat\ga^2$ on $\pd D$.
Using \eq{us4eq33}, \eq{us4eq35} and \eq{us4eq1} we can derive upper
bounds for $\hat u-\hat u(x_0,y_0)$ and $\hat v-\hat v_0$ on $\pd D$.
If they are less than $\ha rN$ and $\ha r^2N$ respectively then
\eq{us4eq45} holds on $\pd D$, and we are finished. This will be
true provided $A,B,C>0$ are chosen small enough to satisfy certain
inequalities involving $N$. We leave the details to the reader.

This completes the proof of Theorem \ref{us4thm1}. The proof of
Theorem \ref{us4thm2} is similar, using Theorem \ref{us4thm6}
rather than Theorem~\ref{us4thm5}.

\section{A priori estimates for $\frac{\pd u}{\pd x},
\frac{\pd u}{\pd y},\frac{\pd v}{\pd x}$ and $\frac{\pd v}{\pd y}$}
\label{us5}

We can now use the results of \S\ref{us4} and
\cite[Th.~6.9]{Joyc3} to derive a priori interior and global
estimates for derivatives of solutions $u,v$ of \eq{us3eq4}
satisfying a $C^0$ bound.

\subsection{Interior estimates for $\frac{\pd u}{\pd x},
\frac{\pd u}{\pd y},\frac{\pd v}{\pd x}$ and $\frac{\pd v}{\pd y}$}
\label{us51}

Here are {\it a priori interior estimates} for $\pd u,\pd v$
when $u,v$ satisfy~\eq{us3eq4}.

\begin{thm} Let\/ $K,L>0$ be given, and\/ $S,T$ be domains in $\R^2$
with\/ $T\subset S^\circ$. Then there exists\/ $R>0$ depending only
on $K,L,S$ and\/ $T$ such that the following is true.

Suppose that\/ $a\in\R$ with\/ $a\ne 0$ and\/ $\md{a}\le K$, and
that\/ $u,v\in C^1(S)$ satisfy \eq{us3eq4} and\/ $u^2+v^2<L^2$.
Then whenever $(x_0,y_0)\in T$, we have
\begin{align}
\Bigl\vert\frac{\pd u}{\pd x}(x_0,y_0)\Bigr\vert
&\le\sqrt{2}\,R\bigl(v(x_0,y_0)^2+y_0^2+a^2\bigr)^{-1},
\label{us5eq1}\\
\Bigl\vert\frac{\pd u}{\pd y}(x_0,y_0)\Bigr\vert &\le R
\bigl(v(x_0,y_0)^2+y_0^2+a^2\bigr)^{-5/4},
\label{us5eq2}\\
\Bigl\vert\frac{\pd v}{\pd x}(x_0,y_0)\Bigr\vert &\le 
2R\bigl(v(x_0,y_0)^2+y_0^2+a^2\bigr)^{-3/4},
\label{us5eq3}\\
\text{and}\quad \Bigl\vert\frac{\pd v}{\pd y}(x_0,y_0)\Bigr\vert &\le
\sqrt{2}\,R\bigl(v(x_0,y_0)^2+y_0^2+a^2\bigr)^{-1}.
\label{us5eq4}
\end{align}
\label{us5thm1}
\end{thm}

This gives good estimates of $\frac{\pd u}{\pd x},\frac{\pd u}{\pd y},
\frac{\pd v}{\pd x},\frac{\pd v}{\pd y}$ except when $v,y,a\approx 0$.
But the equations \eq{us3eq4} are singular exactly when $v=y=a=0$. So,
we have good estimates of the derivatives of $u,v$ except when we are
close to a singular point.

The proof of Theorem \ref{us5thm1} uses \cite[Th.~6.9]{Joyc3}, which
was explained in \S\ref{us1}, and which readers are advised to consult
at this point. Applying this result with $(\hat u,\hat v)$ equal to one
of the families $(u,v)$, $(u',v')$ constructed in \S\ref{us4} gives
nonexistence results for $(u,v)$ with $(u,v)(x_0,y_0)=(u_0,v_0)$
and prescribed $\pd v(x_0,y_0)$. In this way we exclude all values
of $\pd v(x_0,y_0)$ except those allowed by \eq{us5eq3}--\eq{us5eq4},
which are equivalent to \eq{us5eq1}--\eq{us5eq2} by~\eq{us3eq4}.

Theorem \ref{us5thm1} follows from the next two theorems, which
combine \cite[Th.~6.9]{Joyc3} with Theorems \ref{us4thm1} and
\ref{us4thm2} respectively. Note that we need {\it both\/} Theorems
\ref{us4thm1} and \ref{us4thm2} to prove Theorem~\ref{us5thm1}.

\begin{thm} Let\/ $K,L>0$ be given, and\/ $S,T$ be domains in $\R^2$
with\/ $T\subset S^\circ$. Then there exists\/ $R'>0$ depending only
on $K,L,S$ and\/ $T$ such that the following is true.

Suppose that\/ $a\in\R$ with\/ $a\ne 0$ and\/ $\md{a}\le K$, and
that\/ $u,v\in C^1(S)$ satisfy \eq{us3eq4} and\/ $u^2+v^2<L^2$.
Then whenever $(x_0,y_0)\in T$, we have
\begin{align}
\Bigl\vert\frac{\pd u}{\pd x}(x_0,y_0)\Bigr\vert
&\le\sqrt{2}\,R'\bigl(y_0^2+a^2\bigr)^{-1/2},
\label{us5eq5}\\
\Bigl\vert\frac{\pd u}{\pd y}(x_0,y_0)\Bigr\vert &\le R'\bigl(
v(x_0,y_0)^2+y_0^2+a^2\bigr)^{-1/4}\bigl(y_0^2+a^2\bigr)^{-1/2},
\label{us5eq6}\\
\Bigl\vert\frac{\pd v}{\pd x}(x_0,y_0)\Bigr\vert &\le 2R'\bigl(
v(x_0,y_0)^2+y_0^2+a^2\bigr)^{1/4}\bigl(y_0^2+a^2\bigr)^{-1/2},
\label{us5eq7}\\
\text{and}\quad \Bigl\vert\frac{\pd v}{\pd y}(x_0,y_0)\Bigr\vert 
&\le\sqrt{2}\,R'\bigl(y_0^2+a^2\bigr)^{-1/2}.
\label{us5eq8}
\end{align}
\label{us5thm2}
\end{thm}

\begin{proof} As $u,v$ satisfy \eq{us3eq4}, equations \eq{us5eq5}
and \eq{us5eq6} follow from \eq{us5eq7} and \eq{us5eq8}. So it is
enough to prove \eq{us5eq7} and \eq{us5eq8}. Define
\begin{equation*}
M=\sup_{(x,y)\in T}(x^2+y^2)^{1/2}\quad\text{and}\quad
N=\sup\bigl\{\ep>0:B_\ep(x,y)\subset S\;\> \forall (x,y)\in T\bigr\}.
\end{equation*}
Then $M,N>0$ are well-defined, as $T$ is compact and~$T\subset S^\circ$.

Let $A,B,C>0$ be as in Theorem \ref{us4thm1}, using these $K,L,M,N$,
and define
\e
R'=\max\bigl(2^{1/4}A^{-1},2^{-1/2}LB^{-1},
2^{-1/2}(M^2+K^2)^{1/2}C^{-1}\bigr).
\label{us5eq9}
\e
Then $R'$ depends only on $K,L,S$ and $T$, as $M,N$ and $A,B,C$ do. Define
\begin{gather}
\begin{aligned}
u_0&=u(x_0,y_0),&\;\> v_0&=v(x_0,y_0), &\;\>
p_0&=\ts\frac{\pd v}{\pd x}(x_0,y_0),&\;\>
q_0&=\ts\frac{\pd v}{\pd y}(x_0,y_0),\\
\hat x_0&=u_0,&\;\> \hat y_0&=v_0,\;\>&
\hat u_0&=x_0,&\;\> \hat v_0&=y_0,
\end{aligned}
\label{us5eq10}\\
\hat p_0=-\,\frac{p_0}{\ha(v_0^2\!+\!y_0^2\!+\!a^2)^{-1/2}p_0^2\!+\!q_0^2}
\;\>\text{and}\;\>
\hat q_0=\frac{q_0}{\ha(v_0^2\!+\!y_0^2\!+\!a^2)^{-1/2}p_0^2\!+\!q_0^2}.
\label{us5eq11}
\end{gather}
A straightforward calculation now shows:

\begin{lem} Suppose either \eq{us5eq7} or \eq{us5eq8} does not 
hold. Then $a,\hat x_0,\hat y_0,\hat u_0$, $\hat v_0,\hat p_0$
and\/ $\hat q_0$ satisfy equation \eq{us4eq1}, replacing $x_0$ 
by $\hat x_0$, and so on.
\label{us5lem1}
\end{lem}

We can now finish the proof of Theorem \ref{us5thm2}. Suppose
that either \eq{us5eq7} or \eq{us5eq8} does not hold. Then by 
Lemma \ref{us5lem1}, the hypotheses of Theorem \ref{us4thm1} 
hold, with $x_0,y_0,u_0,v_0,p_0,q_0$ replaced by $\hat x_0,
\hat y_0,\hat u_0,\hat v_0,\hat p_0,\hat q_0$. Hence Theorem 
\ref{us4thm1} gives $\hat u,\hat v\in C^\iy(D_L)$ satisfying 
\eq{us3eq4}, $(\hat u-\hat u_0)^2+(\hat v-\hat v_0)^2<N^2$
and~\eq{us4eq2}.

But $(\hat u_0,\hat v_0)=(x_0,y_0)$ which lies in $T$, and so
the open ball $B_N(x_0,y_0)$ of radius $N$ about $(x_0,y_0)$
lies in $S$, by definition of $N$. Therefore $(\hat u,\hat v)$
maps $D_L\ra S$. Applying \cite[Th.~6.9]{Joyc3} then shows that
there do not exist $u,v\in C^1(S)$ satisfying \eq{us3eq4},
$u^2+v^2<L^2$ and \eq{us4eq2}, contradicting the definitions
of $u_0,v_0,p_0$ and $q_0$. Therefore both \eq{us5eq7} and
\eq{us5eq8} hold, and the theorem is complete.
\end{proof}

In the same way, combining Theorem \ref{us4thm2}
and \cite[Th.~6.9]{Joyc3} we prove:

\begin{thm} Let\/ $K,L>0$ be given, and\/ $S,T$ be domains in $\R^2$
with\/ $T\subset S^\circ$. Then there exists\/ $R''>0$ depending only
on $K,L,S$ and\/ $T$ such that the following is true.

Suppose that\/ $a\in\R$ with\/ $a\ne 0$ and\/ $\md{a}\le K$, and
that\/ $u,v\in C^1(S)$ satisfy \eq{us3eq4} and\/ $u^2+v^2<L^2$.
Then whenever $(x_0,y_0)\in T$, we have
\begin{align}
\Bigl\vert\frac{\pd u}{\pd x}(x_0,y_0)\Bigr\vert
&\le\sqrt{2}\,R''\bigl(v(x_0,y_0)^2+y_0^2+a^2\bigr)^{1/4}
\bmd{v(x_0,y_0)}^{-5/2},
\label{us5eq12}\\
\Bigl\vert\frac{\pd u}{\pd y}(x_0,y_0)\Bigr\vert &\le R''
\bmd{v(x_0,y_0)}^{-5/2},
\label{us5eq13}\\
\Bigl\vert\frac{\pd v}{\pd x}(x_0,y_0)\Bigr\vert &\le 
2R''\bigl(v(x_0,y_0)^2+y_0^2+a^2\bigr)^{1/2}\bmd{v(x_0,y_0)}^{-5/2},
\label{us5eq14}\\
\text{and}\quad \Bigl\vert\frac{\pd v}{\pd y}(x_0,y_0)\Bigr\vert &\le
\sqrt{2}\,R''\bigl(v(x_0,y_0)^2+y_0^2+a^2\bigr)^{1/4}\bmd{v(x_0,y_0)}^{-5/2}.
\label{us5eq15}
\end{align}
\label{us5thm3}
\end{thm}

We now prove Theorem \ref{us5thm1}. Let $R',R''$ be as in
Theorems \ref{us5thm2}, \ref{us5thm3}, and put
\begin{equation*}
\ts R=\max\bigl(2^{1/2}(L^2+\sup_T\ms{y}+K^2)^{1/2}R',2^{5/4}R''\bigr).
\end{equation*}
Let $(x_0,y_0)\in T$, and divide into the two cases (a) $v(x_0,y_0)^2
\le y_0^2+a^2$ and (b) $v(x_0,y_0)^2>y_0^2+a^2$. In case (a) we have
\begin{align*}
R'(y_0^2+a^2)^{-1/2}&\le 2^{1/2}R'\bigl(v(x_0,y_0)^2+y_0^2+a^2\bigr)^{-1/2}\\
&\ts\le 2^{1/2}(L^2+\sup_T\ms{y}+K^2)^{1/2}R'
\bigl(v(x_0,y_0)^2+y_0^2+a^2\bigr)^{-1}\\
&\le R\bigl(v(x_0,y_0)^2+y_0^2+a^2\bigr)^{-1},
\end{align*}
and in case (b) we have
\begin{align*}
R''\bigl(v(x_0,y_0)^2+y_0^2+a^2\bigr)^{1/4}\bmd{v(x_0,y_0)}^{-5/2}
&\le 2^{5/4}R''\bigl(v(x_0,y_0)^2+y_0^2+a^2\bigr)^{-1}\\
&\le R\bigl(v(x_0,y_0)^2+y_0^2+a^2\bigr)^{-1}.
\end{align*}
Equations \eq{us5eq1}--\eq{us5eq4} then follow from
\eq{us5eq5}--\eq{us5eq8} in case (a), and
\eq{us5eq12}--\eq{us5eq15} in~(b).

\subsection{Global estimates for $\frac{\pd u}{\pd x},
\frac{\pd u}{\pd y},\frac{\pd v}{\pd x}$ and $\frac{\pd v}{\pd y}$}
\label{us52}

In analysis, interior estimates can usually be be extended to global
estimates on the whole domain $S$ by imposing suitable boundary
conditions on $\pd S$. We shall now extend the results of \S\ref{us51}
to all of $S$, provided $\pd u,\pd v$ satisfy certain inequalities
on $\pd S$. Here is our main result.

\begin{thm} Let\/ $J,K,L>0$ be given, and\/ $S$ be a domain in $\R^2$
such that\/ $T_{(x,0)}\pd S$ is parallel to the $y$-axis for each\/
$(x,0)$ in $\pd S$. Then there exists\/ $H>0$ depending only
on $J,K,L$ and\/ $S$ such that the following is true.

Suppose that\/ $a\in\R$ with\/ $a\ne 0$ and\/ $\md{a}\le K$, and
that\/ $u,v\in C^1(S)$ satisfy \eq{us3eq4} and\/ $u^2+v^2<L^2$
on $S$, and\/ $\bmd{\frac{\pd v}{\pd x}}\le J$, 
$\bmd{\frac{\pd v}{\pd y}}\le J(y^2\!+\!a^2)^{-1/2}$ and\/
$\bmd{\frac{\pd v}{\pd y}}\le J(v^2\!+\!y^2\!+\!a^2)^{1/4}
\md{v}^{-5/2}$ on $\pd S$. Then for all\/ $(x_0,y_0)\in S$, we have
\begin{align}
\Bigl\vert\frac{\pd u}{\pd x}(x_0,y_0)\Bigr\vert=
\Bigl\vert\frac{\pd v}{\pd y}(x_0,y_0)\Bigr\vert
&\le H\bigl(v(x_0,y_0)^2+y_0^2+a^2\bigr)^{-1},
\label{us5eq16}\\
\Bigl\vert\frac{\pd v}{\pd x}(x_0,y_0)\Bigr\vert\le J
\quad\text{and\/}\quad
\Bigl\vert\frac{\pd u}{\pd y}(x_0,y_0)\Bigr\vert &\le\ha J
\bigl(v(x_0,y_0)^2+y_0^2+a^2\bigr)^{-1/2}.
\label{us5eq17}
\end{align}
\label{us5thm4}
\end{thm}

This is a global estimate for derivatives of solutions $u,v$ of
\eq{us3eq4} on a domain $S$ satisfying certain bounds on $\pd S$,
similar to the interior estimates in Theorem \ref{us5thm1}. Here
\eq{us5eq16} is essentially the same as \eq{us5eq1} and \eq{us5eq4},
but \eq{us5eq17} is stronger than \eq{us5eq2} and \eq{us5eq3}. This
is because $\frac{\pd v}{\pd x}$ satisfies a {\it maximum principle}
on $S$ \cite[Prop.~8.12]{Joyc3}, so $\bmd{\frac{\pd v}{\pd x}}\le J$
on $\pd S$ implies $\bmd{\frac{\pd v}{\pd x}}\le J$ on~$S$.

The bounds needed in Theorem \ref{us5thm4} for $\pd v$ on $\pd S$
are quite strong, in that we assume a bound on {\it all} of $\pd v$,
but in applications such as the Dirichlet problem for $v$, we initially
only know a bound on {\it half\/} of $\pd v$. In \S\ref{us6} we will
implicitly show how to extend this to bound all of $\pd v$, so that
Theorem \ref{us5thm4} applies.

The rest of the section proves Theorem \ref{us5thm4}. Consider
the following situation.

\begin{dfn} Let $S$ be a domain in $\R^2$, such that for every point
of the form $(x,0)$ in $\pd S$, the tangent line $T_{(x,0)}\pd S$ is
parallel to the $y$-axis. Let $J,K,L>0$ be given, and suppose that
$a\in\R$ with $a\ne 0$ and $\md{a}\le K$, and that $u,v\in C^1(S)$
satisfy \eq{us3eq4} and $u^2+v^2<L^2$ on $S$, and
$\bmd{\frac{\pd v}{\pd x}}\le J$, $\bmd{\frac{\pd v}{\pd y}}\le
J(y^2\!+\!a^2)^{-1/2}$ and $\bmd{\frac{\pd v}{\pd y}}\le
J(v^2\!+\!y^2\!+\!a^2)^{1/4}\md{v}^{-5/2}$ on~$\pd S$.
\label{us5def}
\end{dfn}

We shall show that under these assumptions, $\pd u,\pd v$ satisfy
estimates similar to \eq{us5eq1}--\eq{us5eq4} on the whole of $S$.
From \cite[Prop.~8.12]{Joyc3} and \eq{us3eq4} we have:

\begin{cor} We have $\bmd{\frac{\pd v}{\pd x}}\le J$ and\/
$\bmd{\frac{\pd u}{\pd y}}\le\ha J(v^2+y^2+a^2)^{-1/2}$ on~$S$.
\label{us5cor}
\end{cor}

Thus the problem is to estimate $\frac{\pd u}{\pd x}=\frac{\pd v}{\pd y}$.
We begin by bounding $\frac{\pd v}{\pd y}$ away from the $x$-axis.

\begin{prop} Let\/ $\ep>0$ be small, and set\/ $S_\ep=\bigl\{(x,y)\in
S:\md{y}>\ep\bigr\}$. Then there exists $G>0$ depending only on
$S,J,K,L,\ep$ such that\/~$\bcnm{\frac{\pd v}{\pd y}\vert_{S_\ep}}{0}\le G$.
\label{us5prop1}
\end{prop}

\begin{proof} The proposition will follow from an {\it interior
regularity result} for {\it quasilinear elliptic equations} on
{\it domains in $\R^2$ with boundary portions}. For {\it linear\/}
operators of a certain form this is done in Gilbarg and Trudinger
\cite[p.~302-4]{GiTr}, and will be discussed in the proof of
Proposition \ref{us7prop1}. For {\it quasilinear} operators we
can deduce what we need from the proof of~\cite[Th.~15.2]{GiTr}.

This says that if $Q$ is a quasilinear elliptic operator of the
form \eq{us2eq2}, then if $a^{ij}$ and $b$ and their first
derivatives satisfy certain complicated estimates on a domain
$S$, and $v\in C^2(S)$ with $Qv=0$, then $\cnm{\pd v}{0}\le G$
for some $G>0$ depending on $S$, $\cnm{v}{0}$, $\cnm{\pd v
\vert_{\pd S}}{0}$, and quantities in the estimates on $a^{ij}$
and~$b$.

We may extend this result to an interior regularity result for
{\it domains with boundary portions}, using the ideas of
\cite[p.~302-4]{GiTr}. That is, we suppose that the estimates
hold in $S_{\ep/2}$, but deduce the a priori bound on $S_\ep$.
The important point here is that $S_{\ep/2}$ is {\it noncompact},
and its boundary $\pd S_{\ep/2}$ is part of $\pd S$. The closure
$\overline S_{\ep/2}$ has an {\it extra boundary portion}, two
line segments with $y=\pm\ep/2$. But we do {\it not\/} need a
bound for $\pd v$ on these line segments.

The price of this is that we can only bound $\pd v$ away from
$y=\pm\ep/2$, which is why we end up with an a priori bound for
$\cnm{\pd v\vert_{S_\ep}}{0}$ rather than $\cnm{\pd v\vert_{
S_{\ep/2}}}{0}$. Since $\md{\frac{\pd v}{\pd x}}\le J$, $\md{a}\le K$,
$\md{v}<L$, and $\ha\ep\le y\le\sup_S\md{y}$ on $S_{\ep/2}$, it is not
difficult to show that the necessary estimates on $a^{ij}$ and $b$ hold
at $v$ with constants depending only on $S,J,K,L$ and $\ep$. It then
follows that $\bcnm{\frac{\pd v}{\pd y}\vert_{S_\ep}}{0}\le G$ for
some $G>0$ depending only on $S,J,K,L$ and~$\ep$.
\end{proof}

It remains to bound $\frac{\pd v}{\pd y}$ near the $x$-axis. We
do this by extending Theorems \ref{us5thm2} and \ref{us5thm3}
from interior domains $T$ to all points $(x_0,y_0)$ in $S$
near the $x$-axis. Here is the extension of Theorem~\ref{us5thm2}.

\begin{prop} There exist constants $\ep>0$ depending only on
$S$, and\/ $R'>0$ depending only on $S,J,K,L$ and\/ $\ep$, such
that if\/ $(x_0,y_0)\in S$ with\/ $\md{y_0}<\ep$ then equations
\eq{us5eq5}--\eq{us5eq8} hold.
\label{us5prop2}
\end{prop}

\begin{proof} Let small $\ep>0$ and large $R'>0$ be chosen, to
satisfy conditions we will give later. Since $u,v$ are $C^1$ and
\eq{us5eq5}--\eq{us5eq8} are closed conditions, it is enough to
prove the proposition for $(x_0,y_0)\in S^\circ$ rather than $S$.
So suppose for a contradiction that $(x_0,y_0)\in S^\circ$ with
$\md{y_0}<\ep$, and that \eq{us5eq5}--\eq{us5eq8} do not all hold.
As \eq{us5eq7}--\eq{us5eq8} imply \eq{us5eq5}--\eq{us5eq6} and
\eq{us5eq7} follows from Corollary \ref{us5cor} when $R'\gg J$,
this means that \eq{us5eq8} does not hold.

We follow the proof of Theorem \ref{us5thm2}. Set $M=\sup_{(x,y)\in
S}(x^2+y^2)^{1/2}$ and $N=\ep$. Let $A,B,C>0$ be as in Theorem
\ref{us4thm1}, using these $K,L,M,N$, and choose $R'$ greater than
or equal to the r.h.s.\ of \eq{us5eq9}. Define $u_0,v_0,p_0,q_0,
\hat x_0,\hat y_0,\hat u_0,\hat v_0,\hat p_0$ and $\hat q_0$ as in
\eq{us5eq10}--\eq{us5eq11}. Then Lemma \ref{us5lem1} shows that
$a,\hat x_0,\hat y_0,\hat u_0$, $\hat v_0,\hat p_0$ and $\hat q_0$
satisfy \eq{us4eq1}. Hence Theorem \ref{us4thm1} gives $\hat u,
\hat v\in C^\iy(D_L)$ satisfying \eq{us3eq4}, $(\hat u-\hat u_0)^2
+(\hat v-\hat v_0)^2<\ep^2$ and~\eq{us4eq2}.

Now $\md{p_0}\le J$ by Corollary \eq{us5cor} and $\md{q_0}>
\sqrt{2}R'(y_0^2+a^2)^{-1/2}$ as \eq{us5eq8} does not hold.
Hence $p_0=O(1)$ and $q_0$ is large. Equation \eq{us5eq11} then
gives $\hat q_0\approx q_0^{-1}$ and $\hat p_0=O(\hat q_0^2)$.
The material of \S\ref{us4} then implies that $\hat u,\hat v$
and their first derivatives approximate the affine maps
\e
\hat u(x,y)\approx \hat u_0+\hat q_0(x-\hat x_0),\qquad
\hat v(x,y)\approx \hat v_0+\hat q_0(y-\hat y_0).
\label{us5eq18}
\e

Define $U=(\hat u,\hat v)(D_L)$. Then $U\subset B_\ep(x_0,y_0)$,
the open ball of radius $\ep$ about $(x_0,y_0)$, as $\hat u_0=x_0$,
$\hat v_0=y_0$ and $(\hat u-\hat u_0)^2+(\hat v-\hat v_0)^2<\ep^2$.
Furthermore, \eq{us5eq18} implies that $U$ is approximately a
closed disc of radius $\md{\hat q_0}L$, and that the map
$(\hat u,\hat v):D_L\ra U$ is invertible with differentiable inverse.

Let $(u',v'):U\ra D_L$ be this inverse map. Then
\cite[Prop.~6.8]{Joyc3} implies that $u',v'$ satisfy \eq{us3eq4} in
$U$. Moreover, equations \eq{us4eq2}, \eq{us5eq10} and \eq{us5eq11}
imply that $u=u'$, $v=v'$, $\pd u=\pd u'$ and $\pd v=\pd v'$ at
$(x_0,y_0)$. Thus $(u',v')-(u,v)$ has a zero of {\it multiplicity}
at least 2 at $(x_0,y_0)$, in the sense of~\cite[Def.~6.3]{Joyc3}.

Suppose now that $U\subset S$. Then $(u',v')$ and $(u,v)$ satisfy
\eq{us3eq4} in $U$. Now $(u',v')$ takes $\pd U$ to $\pd D_L$, the
circle of radius $L$, and winds round $\pd D_L$ once in the
positive sense. Since $u^2+v^2<L^2$ it follows that the winding
number of $(u',v')-(u,v)$ about 0 along $\pd U$ is the same as
that of $(u',v')$, which is 1. So by \cite[Th.~6.7]{Joyc3} there
is 1 zero of $(u',v')-(u,v)$ in $U^\circ$, counted with
multiplicity. But this is a contradiction, as $(u',v')-(u,v)$
has a zero of multiplicity at least 2 at $(x_0,y_0)\in U^\circ$.
This proves the proposition when~$U\subset S$.

It remains to deal with the case that $U\not\subset S$. Then $U$
must intersect $\pd S$. Since $\md{y_0}<\ep$ and $U\subset B_\ep
(x_0,y_0)$ it follows that $\md{y}<2\ep$ on $U$. So if $\ep$ is
small enough, $(x_0,y_0)$ and $U$ must be close to a point of
the form $(x,0)$ in $\pd S$. By Definition \ref{us5def} the
tangent line $T_{(x,0)}\pd S$ is parallel to the $y$-axis.
Thus, by making $\ep$ small we can assume that $U$ intersects
a portion of $\pd S$ close to $(x,0)$, and with tangent spaces
nearly parallel to the $y$-axis.

For $s\in\R$, define $U_s=\bigl\{(x+s,y):(x,y)\in U\bigr\}$
and $u_s',v_s':U_s\ra\R$ by $(u_s',v_s')(x,y)=(u',v')(x-s,y)$.
Then $(u_s',v_s'):U_s\ra D_L$ satisfy \eq{us3eq4}. Suppose for
simplicity that $(1,0)$ points inward to $S$ at $(x,0)$. Then
$U_0=U$, and as $s$ increases from zero, $U_s$ moves inwards
into $S$, until $U_t$ lies wholly in $S$ for some~$t=O(\ep)$.

The argument above shows that the number of zeroes of
$(u'_t,v'_t)-(u,v)$ in $U_t$, counted with multiplicity, is 1.
We shall show that as $s$ increases from 0 to $t$, the number
of zeroes of $(u'_s,v'_s)-(u,v)$ in $U_s\cap S$, counted with
multiplicity, can only increase. Hence, the number of zeroes
of $(u',v')-(u,v)$ in $U\cap S$, counted with multiplicity,
is no more than 1, and as above we have a contradiction with
\cite[Th.~6.7]{Joyc3} because $(x_0,y_0)$ is a zero of
multiplicity at least 2, so the proof will be finished.

The only way in which the number of zeroes of $(u'_s,v'_s)-(u,v)$
in $U_s\cap S$ with multiplicity can change as $s$ changes, is if
a zero crosses the boundary $\pd(U_s\cap S)$. This consists of two
portions, $\pd U_s\cap S$ and $U_s\cap\pd S$. But on $\pd U_s\cap S$
we have $(u'_s)^2+(v'_s)^2=L^2$ and $u^2+v^2<L^2$, so no zeroes can
cross this part of the boundary. Thus, zeroes can only enter or
leave $U_s\cap S$ across~$U_s^\circ\cap\pd S$.

The next lemma computes $\frac{\d}{\d s}$ of a zero $\bigl(x(s),y(s)
\bigr)$ of $(u'_s,v'_s)-(u,v)$, in terms of $\pd v_s'$ and $\pd v$.
The proof is an elementary calculation using \eq{us3eq4}, and is left
to the reader.

\begin{lem} Suppose $\bigl(x(s),y(s)\bigr)$ is a zero of\/
$(u'_s,v'_s)-(u,v)$ in $U_s\cap S$. Then
\e
\begin{split}
\Bigl(\bigl({\ts\frac{\pd}{\pd y}(v'_s-v)}\bigr)^2+\ha(v^2+y^2+
a^2)^{-1/2}
\bigl({\ts\frac{\pd}{\pd x}(v'_s-v)}\bigr)^2\Bigr)\cdot
\frac{\d}{\d s}\begin{pmatrix} x(s) \\ y(s)\end{pmatrix}&=\\
\begin{pmatrix}
\ts\frac{\pd}{\pd y}v'_s\cdot
\frac{\pd}{\pd y}(v'_s-v)
+\ha(v^2+y^2+a^2)^{-1/2}\frac{\pd}{\pd x}v'_s
\frac{\pd}{\pd x}(v'_s-v)\\
\ts\frac{\pd}{\pd x}v'_s\cdot
\frac{\pd}{\pd y}(v'_s-v)
-\frac{\pd}{\pd y}v'_s\cdot
\frac{\pd}{\pd x}(v'_s-v)
\end{pmatrix}&.
\end{split}
\label{us5eq19}
\e
\label{us5lem2}
\end{lem}

When $\bigl(x(s),y(s)\bigr)$ actually crosses $\pd S$ we have
$\bmd{\frac{\pd v}{\pd y}}\le J(y^2\!+\!a^2)^{-1/2}$, by assumption.
Also, from \eq{us5eq18} we can show that $\frac{\pd}{\pd y}v'_s
\approx q_0$ and $\frac{\pd}{\pd x}v'_s=O(1)$. Therefore, provided
$R'\gg J$, careful calculation shows that dominant terms on both
sides of \eq{us5eq19} are $(\frac{\pd}{\pd y}v'_s)^2$, and so
$\frac{\d}{\d s}\bigl(x(s),y(s)\bigr)\approx(1,0)$. But
from above, $T_{(x,y)}\pd S$ is nearly parallel to the
$y$-axis, and $(1,0)$ points inwards to $S$. Hence, as $s$
increases from 0 to $t$, zeroes of $(u'_s,v'_s)-(u,v)$ can
only move into $U'_s\cap S$, not out. This completes the
proof of Proposition~\ref{us5prop2}.
\end{proof}

By using the solutions of Theorem \ref{us4thm2} instead of
Theorem \ref{us4thm1}, using a similar proof we obtain a
global analogue of Theorem \ref{us5thm3}. Combining these two
results using the method of Theorem \ref{us5thm1} then yields
Theorem~\ref{us5thm4}.

\section{The Dirichlet problem for $v$ when $a=0$}
\label{us6}

Theorem \ref{us3thm2} shows that the Dirichlet problem for equation
\eq{us3eq8} is uniquely solvable in arbitrary domains in $\R^2$ for
$a\ne 0$. In this section we will show that the Dirichlet problem for
\eq{us3eq8} also has a unique {\it weak} solution when $a=0$, for
strictly convex domains $S$ invariant under $(x,y)\mapsto(x,-y)$.

\subsection{The main results}
\label{us61}

The following theorem is an analogue of Theorem \ref{us3thm2}
for the case $a=0$, and the first main result of the paper.

\begin{thm} Let\/ $S$ be a strictly convex domain in $\R^2$ invariant
under the involution $(x,y)\mapsto(x,-y)$, let\/ $k\ge 0$ and\/ $\al\in
(0,1)$. Suppose $\phi\in C^{k+2,\al}(\pd S)$ with\/ $\phi(x,0)\ne 0$
for points $(x,0)$ in $\pd S$. Then there exists a unique weak solution
$v$ of\/ \eq{us3eq8} in $C^0(S)$ with\/ $a=0$ and\/~$v\vert_{\pd S}=\phi$.

Fix a basepoint\/ $(x_0,y_0)\in S$. Then there exists a unique $u\in
C^0(S)$ with\/ $u(x_0,y_0)=0$ such that\/ $u,v$ are weakly differentiable
in $S$ and satisfy \eq{us3eq3} with weak derivatives. The weak derivatives
$\frac{\pd u}{\pd x},\frac{\pd u}{\pd y},\frac{\pd v}{\pd x},
\frac{\pd v}{\pd y}$ satisfy $\frac{\pd u}{\pd x}=\frac{\pd v}{\pd y}
\in L^p(S)$ for $p\in[1,\frac{5}{2})$, and\/ $\frac{\pd u}{\pd y}\in
L^q(S)$ for $q\in[1,2)$, and\/ $\frac{\pd v}{\pd x}$ is bounded on $S$.
Also $u,v$ are $C^{k+2,\al}$ in $S$ and real analytic in $S^\circ$
except at singular points $(x,0)$ with\/~$v(x,0)=0$.
\label{us6thm1}
\end{thm}

Combined with Proposition \ref{us3prop1} the theorem can be used to
construct large numbers of {\it $\U(1)$-invariant singular special
Lagrangian $3$-folds} in $\C^3$. This is the principal motivation
for the paper. The singularities of these special Lagrangian 3-folds
will be studied in \cite{Joyc4}. The restriction to boundary data
$\phi$ with $\phi(x,0)\ne 0$ for points $(x,0)$ in $\pd S$ is to
avoid singular points on the boundary~$\pd S$. 

Our second main result extends \cite[Th.~8.9]{Joyc3} to include
the case~$a=0$.

\begin{thm} Let\/ $S$ be a strictly convex domain in $\R^2$ invariant
under the involution $(x,y)\mapsto(x,-y)$, let\/ $k\ge 0$, $\al\in(0,1)$,
and\/ $(x_0,y_0)\in S$. Define $X$ to be the set of\/ $\phi\in C^{k+2,\al}
(\pd S)$ with\/ $\phi(x,0)=0$ for some $(x,0)\in\pd S$. Then the map
$C^{k+2,\al}(\pd S)\t\R\sm X\t\{0\}\ra C^0(S)^2$ taking $(\phi,a)
\mapsto(u,v)$ is continuous, where $(u,v)$ is the unique solution of\/
\eq{us3eq4} (with weak derivatives when $a=0$) with\/ $v\vert_{\pd S}=\phi$
and\/ $u(x_0,y_0)=0$, constructed in Theorem \ref{us3thm2} when $a\ne 0$,
and in Theorem \ref{us6thm1} when $a=0$. This map is also continuous in
stronger topologies on $(u,v)$ than the $C^0$ topology.
\label{us6thm2}
\end{thm}

The proofs of Theorems \ref{us6thm1} and \ref{us6thm2} will
take up the rest of the section. Here is how they are laid out.
Let $S,\phi$ be as in Theorem \ref{us6thm1}. In \S\ref{us62},
for each $a\in(0,1]$ we define $v_a\in C^{k+2,\al}(S)$ to be
the unique solution of \eq{us3eq8} in $S$ with $v_a\vert_{\pd S}
=\phi$, and $u_a,f_a$ such that $u_a,v_a$ satisfy \eq{us3eq4}
and $\frac{\pd f_a}{\pd y}=u_a$, $\frac{\pd f_a}{\pd x}=v_a$.
The idea is to show that $u_a,v_a,f_a\ra u_0,v_0,f_0$ as $a\ra 0_+$,
for unique, suitably differentiable $u_0,v_0,f_0\in C^0(S)$. Then
$u_0,v_0$ are the weak solutions $u,v$ in Theorem~\ref{us6thm1}.

To show that these limits $u_0,v_0,f_0$ exist, the main issue
is to prove {\it a priori estimates} of $u_a,v_a,f_a$ that are
{\it uniform in} $a$. That is, we need bounds such as $\cnm{u_a}0
\le C$ for all $a\in(0,1]$, with $C$ independent of $a$. Given
strong enough uniform a priori estimates, the existence of {\it
some} weak limits $u_0,v_0,f_0$ becomes essentially trivial,
using compact embeddings of Banach spaces.

Getting such uniform a priori estimates is difficult, since
equations \eq{us3eq4} and \eq{us3eq8} {\it really are} singular
when $a=0$, so many norms of $u_a,v_a$ such as $\cnm{\pd u_a}0$,
$\cnm{\pd v_a}0$ can diverge to infinity as $a\ra 0_+$, and
uniform a priori estimates of these norms {\it do not exist}.
The part of Theorem \ref{us6thm1} that gave the author most
trouble was finding estimates strong enough to prove that $u$
is {\it continuous}.

This is important geometrically, as if $u,v$ are not continuous
then the SL 3-fold $N$ in \eq{us3eq2} is not locally closed,
and one singular point of $u,v$ will correspond to many singular
points of $\overline N$ rather than one. To show $u$ is continuous
we use the {\it nonstandard Sobolev embedding result\/} Theorem
\ref{us2thm}, which allows us to trade off a stronger $L^p$
estimate of $\frac{\pd u}{\pd x}$ against a weaker $L^q$
estimate of~$\frac{\pd u}{\pd y}$.

The a priori estimates that we need are built up step by step
in \S\ref{us63}--\S\ref{us65}. In \S\ref{us63} we construct
super- and subsolutions for $v_a$ near points $(x,0)$ in
$\pd S$. This gives a positive lower bound for $\md{v_a}$
near $(x,0)$, which proves that \eq{us3eq4} and \eq{us3eq8}
are {\it uniformly elliptic} in $a$ close to $\pd S$. Then
\S\ref{us64} proves uniform $C^0$ estimates for
$u_a,v_a,f_a$ and some derivatives on $S$ or~$\pd S$.

In \S\ref{us65} we use the results of \S\ref{us5} to
prove uniform $L^p$ estimates for $\smash{\frac{\pd u_a}{\pd x}}$,
$\smash{\frac{\pd u_a}{\pd y}}$, $\smash{\frac{\pd v_a}{\pd x}}$
and $\smash{\frac{\pd v_a}{\pd y}}$, and deduce uniform
continuity of the $u_a,v_a$ from Theorem \ref{us2thm}. Section
\ref{us66} proves the existence of limits $u_0,v_0,f_0$
in $C^0$, and that they satisfy \eq{us3eq3}, \eq{us3eq5}
and \eq{us3eq8} in the appropriate weak senses.
Section \ref{us67} proves uniqueness of the solutions,
and completes the proofs.

The requirement that $S$ be a strictly convex domain in $\R^2$
invariant under $(x,y)\mapsto(x,-y)$ is unnecessarily strong. All
the proofs above actually use is that $S$ should be a domain in
$\R^2$, and that for every point $(x,0)$ in $\pd S$, the tangent
to $\pd S$ at $(x,0)$ should be parallel to the $y$-axis, and $S$
should be strictly convex near~$(x,0)$.

The author believes that Theorem \ref{us6thm1} actually holds for
arbitrary domains $S$ in $\R^2$. To extend the proof to such $S$
would need suitable super- and subsolutions for $v_a$ near $(x,0)$
in $\pd S$, generalizing those of \S\ref{us63}. Perhaps one can
use the above Theorem \ref{us6thm1} to construct such super- and
subsolutions.

\subsection{A family of solutions of \eq{us3eq8}}
\label{us62}

We shall consider the following situation.

\begin{dfn} Let $S$ be a strictly convex domain in $\R^2$ which
is invariant under the involution $(x,y)\mapsto(x,-y)$. Then there
exist unique $x_1,x_2\in\R$ with $x_1<x_2$ and $(x_i,0)\in\pd S$
for $i=1,2$. Let $k\ge 0$ and $\al\in(0,1)$, and suppose $\phi\in
C^{k+2,\al}(\pd S)$ with $\phi(x_i,0)\ne 0$ for $i=1,2$. For each
$a\in(0,1]$, let $v_a\in C^{k+2,\al}(S)$ be the unique solution of
\eq{us3eq8} in $S$ with this value of $a$ and $v_a\vert_{\pd S}=\phi$,
which exists and is unique by Theorem~\ref{us3thm2}.

For each $a\in(0,1]$, let $u_a\in C^{k+2,\al}(S)$ be the unique
function such that $u_a,v_a$ and $a$ satisfy \eq{us3eq4}, and
$u_a(x_1,0)=0$. This exists by Proposition \ref{us3prop3}. Let
$f_a\in C^{k+3,\al}(S)$ be the unique solution of \eq{us3eq5}
satisfying $\frac{\pd f_a}{\pd y}=u_a$, $\frac{\pd f_a}{\pd x}=v_a$
and $f_a(x_1,0)=0$. This exists by Proposition~\ref{us3prop2}.
\label{us6def1}
\end{dfn}

We will show that $v_a$ converges in $C^0(S)$ to $v_0\in C^0(S)$ as
$a\ra 0_+$, and that $v_0$ is the unique weak solution of the Dirichlet
problem for \eq{us3eq8} on $S$ when $a=0$. The reason for supposing
that $\phi(x_i,0)\ne 0$ is to avoid having singular points on the
boundary~$\pd S$.

We begin by defining a family of solutions $v_{a,\ga}$ of
\eq{us3eq8} which we will use in \S\ref{us63} as super- and
subsolutions to bound the $v_a$ near the $(x_i,0)$ for~$i=1,2$.

\begin{dfn} Let $R>0$, and define $D_R$ to be the closed disc
of radius $R$ about $(0,0)$ in $\R^2$. For each $a\in(0,1]$ and
$\ga>0$, define $v_{a,\ga}\in C^\iy(D_R)$ to be the unique
solution of \eq{us3eq8} in $D_R$ with this value of $a$ and
$v_a\vert_{\pd S}=\ga x$, which exists and is unique by Theorem
\ref{us3thm2}. Considering how $v_{a,\ga}$ transforms under the
involutions $(x,y)\mapsto(-x,y)$ and $(x,y)\mapsto(x,-y)$, by
uniqueness we see that $v_{a,\ga}$ satisfies the identities
$v_{a,\ga}(-x,y)=-v_{a,\ga}(x,y)$, $v_{a,\ga}(x,-y)=v_{a,\ga}(x,y)$
and~$v_{a,\ga}(0,y)=0$.
\label{us6def2}
\end{dfn}

Provided $\ga$ is large enough, these $v_{a,\ga}$ satisfy
certain inequalities on~$D_R$.

\begin{prop} There exists $C>0$ depending only on $R$ such
that whenever $a\in(0,1]$ and\/ $\ga\ge C$, the function
$v_{a,\ga}$ of Definition \ref{us6def2} satisfies
$v_{a,\ga}(x,y)>0$ when $x>0$, $v_{a,\ga}(x,y)<0$ when $x<0$
and\/ $v_{a,\ga}(0,y)=0$, and
\e
\bmd{v_{a,\ga}}\le\ga\md{x}
\quad\text{and}\quad
v_{a,\ga}^2+y^2\ge x^4
\quad\text{on $D_R$.}
\label{us6eq1}
\e
\label{us6prop1}
\end{prop}

\begin{proof} In \cite[Th.~5.1]{Joyc3} the author defined an
explicit family of solutions $\hat u_a,\hat v_a$ of \eq{us3eq4} on $\R^2$
for all $a\ge 0$, with the properties that $\hat v_a(x,y)>0$ when
$x>0$, $\hat v_a(x,y)<0$ when $x<0$ and $\hat v_a(0,y)=0$ for all $y$, and
\e
\hat v_a^2+y^2\equiv(x^2+\hat u_a^2)(x^2+\hat u_a^2+2a).
\label{us6eq2}
\e
Choose $C>0$ such that $\md{\hat v_a(x,y)}\le C\md{x}$ whenever
$(x,y)\in\pd D_R$, for all $a\in[0,1]$. This is possible
because the $\hat v_a$ are smooth on $\pd D_R$ and depend
smoothly on $a$, and~$\hat v_a(0,y)\equiv 0$.

Define $S_R$ to be the semicircle $\bigl\{(x,y)\in D_R:x\ge 0
\bigr\}$. Let $\ga\ge C$. Then for each $a\in(0,1]$ we have
$\hat v_a\le Cx\le \ga x=v_{a,\ga}$ on $\pd S_R$. Since $\hat v_a$
and $v_{a,\ga}$ satisfy \eq{us3eq8} in $S_R$, we see from
\cite[Prop.~8.5]{Joyc3} on $S_R$ with $v=\hat v_a$ and $v'=
v_{a,\ga}$ that $\hat v_a\le v_{a,\ga}$ in $S_R$. But $\hat
v_a(x,y)>0$ for $x>0$, and thus $v_{a,\ga}(x,y)>0$ when~$x>0$.

Now $v_{a,\ga}^2\ge\hat v_a^2$ in $S_R$, as $v_{a,\ga}\ge\hat v_a\ge 0$
there. But it follows from \eq{us6eq2} that $\hat v_a^2+y^2\ge x^4$ in
$\R^2$. Hence $v_{a,\ga}^2+y^2\ge x^4$ in $S_R$. Also
$v_{a,\ga}\le\ga x$ on $\pd S_R$, and
\begin{equation*}
Q(\ga x)=-(\ga^2x^2+y^2+a^2)^{-3/2}\ga^3x\le 0 \quad\text{in $S_R$.}
\end{equation*}
Thus, applying \cite[Prop.~8.5]{Joyc3} on $S_R$ with $v=v_{a,\ga}$
and $v'=\ga x$ gives $v_{a,\ga}\le\ga x$ on $S_R$, so that
$\md{v_{a,\ga}}\le\ga\md{x}$ on $S_R$. This proves all the assertions of
the proposition on $S_R$, that is, when $x\ge 0$. The case $x<0$ follows
immediately using the identity~$v_{a,\ga}(-x,y)=-v_{a,\ga}(x,y)$.
\end{proof}

\subsection{Super- and subsolutions for $v_a$ near $(x_j,0)$}
\label{us63}

Next we find a uniform positive lower bound for $v_a^2+y^2+a^2$
near $\pd S$ for all $a\in(0,1]$. This will ensure that \eq{us3eq8}
is uniformly elliptic at the $v_a$ near $\pd S$. The difficulty is
to estimate the $v_a$ near points $(x,0)$ in $\pd S$, that is, near
$(x_1,0)$ and $(x_2,0)$. We will do this by using the solutions
$v_{a,\ga}$ of Proposition \ref{us6prop1} as super- and subsolutions.

\begin{prop} In the situation above, there exists $\de>0$
such that whenever $(x,y)\in S$ with\/ $x\le x_1+\de$ then
$v_a(x,y)^2+y^2\ge(x-x_1-\de)^4$ for all\/ $a\in(0,1]$, and
whenever $(x,y)\in S$ with\/ $x\ge x_2-\de$ then
$v_a(x,y)^2+y^2\ge(x-x_2+\de)^4$ for all\/~$a\in(0,1]$.
\label{us6prop2}
\end{prop}

\begin{proof} We first prove the estimate near $(x_1,0)$.
Suppose $\phi(x_1,0)>0$. Choose large $\ga>0$ and small
$\de,R>0$ such that the following conditions hold:
\begin{itemize}
\setlength{\itemsep}{0pt}
\setlength{\parsep}{0pt}
\item[(a)] $\ga\ge C$, where $C$ is given in Proposition~\ref{us6prop1};
\item[(b)] for all $(x,y)\in S$ with $(x-x_1-\de)^2+y^2=R^2$, we
have~$\ga(x_1-\de-x)\le\inf_{\pd S}\phi$;
\item[(c)] for all $(x,y)\in\pd S$ with $(x-x_1-\de)^2+y^2\le R^2$
and $x\ge x_1+\de$, we have $\phi(x,y)\ge 0$; and
\item[(d)] for all $(x,y)\in\pd S$ with $(x-x_1-\de)^2+y^2\le R^2$
and $x\le x_1+\de$, we have~$\phi(x,y)\ge\ga(x_1+\de-x)$.
\end{itemize}
For small enough $\de$ and $R$, part (c) holds automatically and
parts (b) and (d) are approximately equivalent to $\ga\de-\ga R^2
/2\ka\le\inf_{\pd S}\phi$ and $\ga\de\le\phi(x_1,0)$, where $\ka>0$
is the radius of curvature of $\pd S$ at $(x_1,0)$. It is then easy
to see that if $\ga,\de,R$ satisfy $1\ll R^{-2}\ll\ga\ll\de^{-1}$
then all the conditions hold.

Define $T=\bigl\{(x,y)\in S:(x-x_1-\de)^2+y^2\le R^2\bigr\}$, and
for each $a\in(0,1]$ define $v'_a\in C^\iy(T)$ by $v'_a(x,y)=
v_{a,\ga}(x_1+\de-x,y)$, where $v_{a,\ga}$ is given in Proposition
\ref{us6prop1}. Then $v'_a$ and $a$ satisfy \eq{us3eq8} in $T$, as
$v_{a,\ga}$ and $a$ do. Now $T$ is a domain with piecewise-smooth
boundary, which consists of two portions, the first an arc of the
circle of radius $R$ about $(x_1+\de,y)$, and the second a part
of~$\pd S$.

We claim that $v_a\ge v'_a$ on $\pd T$. On the first portion of
$\pd T$, the circle arc, we have $v'_a(x,y)=\ga(x_1-\de-x)$ by
definition of $v_{a,\ga}$, and thus the claim follows from part (b)
above. On the second portion of $\pd T$, the part of $\pd S$, the
claim follows from parts (c) and (d) above and the facts that
$v'_a(x,y)\le 0$ for $x\ge x_1+\de$ and $v'_a(x,y)\le\ga
(x_1+\de-x)$ for $x\ge x_1+\de$, which in turn follow from the
statements in Proposition \ref{us6prop1} that $v_{a,\ga}(x,y)\le 0$
when $x\le 0$ and $v_{a,\ga}(x,y)\le\ga x$ when~$x\ge 0$. 

Thus $v_a$ and $v'_a$ both satisfy \eq{us3eq8} in $T$, and
$v_a\ge v'_a$ on $\pd T$. So by \cite[Prop.~8.5]{Joyc3} we
have $v_a\ge v'_a$ on $T$. If $(x,y)\in S$ with $x\le x_1+\de$
then $(x,y)\in T$, so that $v_a(x,y)\ge v'_a(x,y)$, and also
$v'_a(x,y)\ge 0$ and $v'_a(x,y)^2+y^2\ge(x-x_1-\de)^4$ by
Proposition \ref{us6prop1}. Combining these
gives~$v_a(x,y)^2+y^2\ge(x-x_1-\de)^4$. 

This proves the estimate near $(x_1,0)$ in the case $\phi(x_1,0)>0$.
For the case $\phi(x_1,0)<0$ we instead define $v'_a(x,y)=-v_{a,\ga}
(x_1+\de-x,y)$, and use it as a supersolution rather than a subsolution.
For the estimate near $(x_2,0)$ we use the a similar argument, with
$v'_a(x,y)=\pm v_{a,\ga}(x-x_2+\de,y)$. This completes the proof.
\end{proof}

The proposition implies that if $(x,y)\in S$ is close to $(x_1,0)$
or $(x_2,0)$ then $v_a(x,y)^2+y^2$ is uniformly bounded below by a
positive constant. But if $(x,y)\in S$ is close to $\pd S$ then
either $(x,y)$ is close to $(x_1,0)$ or $(x_2,0)$, or else $\md{y}$
is bounded below by a positive constant, and hence $v_a(x,y)^2+y^2$
is uniformly bounded below by a positive constant. Thus we may prove:

\begin{cor} There exist\/ $\ep,J>0$ such that whenever $(x,y)\in S$ lies
within distance $2\ep$ of\/ $\pd S$, then $v_a(x,y)^2+y^2+a^2\ge J^2$
for all\/~$a\in(0,1]$.
\label{us6cor1}
\end{cor}

\subsection{Estimates on $u_a,v_a$ and $f_a$ on $\pd S$}
\label{us64}

Corollary \ref{us6cor1} implies that \eq{us3eq8} is {\it uniformly
elliptic} at $v_a$ near $\pd S$ for all $a\in(0,1]$. We shall use
this to prove estimates on $u_a,v_a$ and $f_a$ near $\pd S$. We begin
by bounding the $\pd v_a$ on~$\pd S$.

\begin{prop} There exists $K>0$ with\/ $\cnm{\pd v_a\vert_{\pd S}}{0}
\le K$ for all\/~$a\in(0,1]$.
\label{us6prop3}
\end{prop}

\begin{proof} Gilbarg and Trudinger \cite[Th.~14.1, p.~337]{GiTr}
show that if $v\in C^2(S)$ satisfies a quasilinear equation $Qv=0$
of the form \eq{us2eq2} on a domain $S$ and $v\vert_{\pd S}=\phi\in
C^2(\pd S)$, then $\cnm{\pd v\vert_{\pd S}}{0}\le K$ for some $K>0$
depending only on $S$, upper bounds for $\cnm{v}{0}$ and $\cnm{\phi}{2}$,
and certain constants to do with $Q$, which ensure that $Q$ is uniformly
elliptic and $b$ not too large.

Examining the proof shows that it is enough for the conditions to
hold within distance $2\ep$ of $\pd S$. We apply this to the solutions
$v_a$ of \eq{us3eq8} for all $a\in(0,1]$. Corollary \ref{us6cor1}
implies that \eq{us3eq8} is uniformly elliptic at the $v_a$ within
distance $2\ep$ of $\pd S$ for all $a\in(0,1]$, and the other conditions
of \cite[Th.~14.1]{GiTr} easily follow. Thus the theorem gives $K>0$
such that $\cnm{\pd v_a\vert_{\pd S}}{0}\le K$ for all~$a\in(0,1]$.
\end{proof}

This implies uniform $C^0$ bounds for $u_a,v_a$ and $C^1$ bounds
for~$f_a$.

\begin{cor} There exist constants $K_1,\ldots,K_5>0$ such that
\e
\begin{gathered}
\cnm{u_a}{0}\le K_1,\;\> \cnm{\pd u_a\vert_{\pd S}}{0}\le K_2,\;\>
\cnm{v_a}{0}\le K_3,\;\> \cnm{\pd v_a\vert_{\pd S}}{0}\le K_4,\\
\ts\cnm{\frac{\pd v_a}{\pd x}}{0}\le K_4,
\quad\text{and\/}\quad \cnm{f_a}{1}\le K_5
\quad\text{for all\/ $a\in(0,1]$.}
\end{gathered}
\label{us6eq3}
\e
\label{us6cor2}
\end{cor}

\begin{proof} As $\md{v_a}$ is maximum on $\pd S$ by
\cite[Cor.~4.4]{Joyc3} and $v_a\vert_{\pd S}=\phi$ we have
$\cnm{v_a}{0}=K_3=\cnm{\phi}{0}$ for all $a$. Proposition
\ref{us6prop3} gives $\cnm{\pd v_a\vert_{\pd S}}{0}\le K_4$ for
all $a$, with $K_4=K$. Thus $\bmd{\frac{\pd v_a}{\pd x}\vert_{\pd S}}
\le K_4$. But the maximum of $\bmd{\frac{\pd v_a}{\pd x}}$ is achieved
on $\pd S$ by \cite[Prop.~8.12]{Joyc3}. Hence $\cnm{\frac{\pd v_a}{\pd
x}}{0}\le K_4$. As $(v_a^2+y^2+a^2)^{-1/2}\le J^{-1}$ on $\pd S$ by
Corollary \ref{us6cor1}, we see from \eq{us3eq4} and $\cnm{\pd v_a
\vert_{\pd S}}{0}\le K_4$ that $\cnm{\pd u_a\vert_{\pd S}}{0}\le K_2$,
with~$K_2=\max(\ha J^{-1},1)K_4$.

Now $u_a(x_1,0)=0$ by definition and $\cnm{\pd u_a\vert_{\pd S}}{0}
\le K_2$, so $\cnm{u_a\vert_{\pd S}}{0}\le\ha K_2\,l(\pd S)$, where
$l(\pd S)$ is the length of $\pd S$. But $u_a$ is maximum on $\pd S$
by \cite[Cor.~4.4]{Joyc3}, so $\cnm{u_a}{0}\le K_1$, where $K_1=\ha
K_2\,l(\pd S)$. Similarly, $\cnm{\pd f_a}{0}\le\cnm{u_a}{0}+\cnm{v_a}{0}
\le K_1+K_3$, and $f_a(x_1,0)=0$, so $\cnm{f_a}{0}\le\cnm{\pd f_a}{0}
\cdot\diam(S)$. Hence $\cnm{f_a}{1}\le K_5$, with~$K_5=(K_1+K_3)(1+
\diam(S))$.
\end{proof}

Corollary \ref{us6cor2} bounded $\frac{\pd v_a}{\pd x}$ uniformly
in $S$ for all $a\in(0,1]$. We now use Theorem \ref{us5thm4} to
prove a slightly weaker result for~$\frac{\pd v_a}{\pd y}$.

\begin{prop} There exists $K_6>0$ such that
\e
\ts\bmd{\frac{\pd u_a}{\pd x}}=\bmd{\frac{\pd v_a}{\pd y}}\le
K_6\bigl(v_a^2+y^2+a^2\bigr)^{-1}
\quad\text{in $S$ for all\/ $a\in(0,1]$.}
\label{us6eq4}
\e
\label{us6prop4}
\end{prop}

\begin{proof} The bound $\cnm{\pd v_a\vert_{\pd S}}{0}\le K$ for all
$a\in(0,1]$ implies that there exists $J>0$ such that the boundary
conditions on $\pd v$ in Theorem \ref{us5thm4} hold for $a$ and $v_a$
for all $a\in(0,1]$. Set $K=1$ and $L=K_1+K_3+1$, so that $u_a^2+v_a^2
\le K_1^2+K_3^2<L^2$ in $S$ by \eq{us6eq3}. Applying Theorem
\ref{us5thm4} with these $J,K,L$ and $S$ gives $H>0$ such that
\eq{us5eq16} holds in $S$ with $v=v_a$ for all $a\in(0,1]$. Setting
$K_6=H$ completes the proof.
\end{proof}

\subsection{Estimates for $\frac{\pd u_a}{\pd x},\frac{\pd u_a}{\pd y},
\frac{\pd v_a}{\pd x},\frac{\pd v_a}{\pd y}$ in $L^p$ and $u_a,v_a$ in
$C^0$}
\label{us65}

Here is an exact expression for a weighted $L^2$ norm of~$\pd v_a$.

\begin{prop} Let\/ $\al\in[0,\ha)$, and define $J(a,v)=-\int_0^v
(w^2+a^2)^{-\al}\d w$. Then for each\/ $a\in(0,1]$ we have
\e
\int_S(v_a^2\!+\!a^2)^{-\al}\Bigl[\ha\bigl(v_a^2\!+\!y^2\!+\!a^2
\bigr)^{-1/2}\Bigl(\frac{\pd v_a}{\pd x}\Bigr)^2\!+\!\Bigl(
\frac{\pd v_a}{\pd y}\Bigr)^2\Bigr]\d x\,\d y
\!=\!\int_{\pd S}J(a,v_a)\d u_a.
\label{us6eq5}
\e
\label{us6prop5}
\end{prop}

\begin{proof} Stokes' Theorem gives
\begin{align*}
\int_{\pd S}J(a,v_a)\d u_a&=-\int_S\d u_a\w\d\bigl(J(a,v_a)\bigr)=
\int_S(v_a^2+a^2)^{-\al}\d u_a\w\d v_a\\
&=\int_S(v_a^2+a^2)^{-\al}\Bigl[\frac{\pd u_a}{\pd x}\frac{\pd v_a}{\pd y}
-\frac{\pd u_a}{\pd y}\frac{\pd v_a}{\pd x}\Bigr]\d x\w\d y\\
&=\int_S(v_a^2\!+\!a^2)^{-\al}\Bigl[\ha\bigl(v_a^2\!+\!y^2\!+\!a^2
\bigr)^{-1/2}\Bigl(\frac{\pd v_a}{\pd x}\Bigr)^2\!+\!\Bigl(
\frac{\pd v_a}{\pd y}\Bigr)^2\Bigr]\d x\,\d y,
\end{align*}
using \eq{us3eq4} to rewrite $\frac{\pd u_a}{\pd x},\frac{\pd u_a}{\pd
y}$ in terms of~$\frac{\pd v_a}{\pd y},\frac{\pd v_a}{\pd x}$.
\end{proof}

We use this to derive uniform $L^p$ estimates on $\frac{\pd u_a}{\pd x},
\frac{\pd u_a}{\pd y},\frac{\pd v_a}{\pd x}$ and~$\frac{\pd v_a}{\pd y}$.

\begin{prop} Let\/ $p\in[2,\frac{5}{2})$, $q\in[1,2)$ and\/
$r\in[1,\iy)$. Then there exist\/ $C_p,C'_q,C''_r>0$ such that
$\lnm{\frac{\pd u_a}{\pd x}}{p}=\lnm{\frac{\pd v_a}{\pd y}}{p}
\le C_p$, $\lnm{\frac{\pd u_a}{\pd y}}{q}\le C'_q$ and\/
$\lnm{\frac{\pd v_a}{\pd x}}{r}\le C''_r$ for all\/~$a\in(0,1]$.
\label{us6prop6}
\end{prop}

\begin{proof} Let $p\in[2,\frac{5}{2})$ and define $\al=p-2$, so
that $\al\in[0,\ha)$. Then \eq{us6eq4} gives
\begin{equation*}
\ts\bmd{\frac{\pd v_a}{\pd y}}^p=\bmd{\frac{\pd v_a}{\pd y}}^\al\cdot
\bms{\frac{\pd v_a}{\pd y}}\le K_6^\al(v_a^2\!+\!y^2\!+\!a^2)^{-\al}
\bms{\frac{\pd v_a}{\pd y}}\le K_6^\al(v_a^2\!+\!a^2)^{-\al}
\bms{\frac{\pd v_a}{\pd y}}.
\end{equation*}
Integrating this over $S$ gives
\begin{equation*}
\ts\blnm{\frac{\pd v_a}{\pd y}}{p}^p\le
K_6^\al\int_S(v_a^2+a^2)^{-\al}\bms{\frac{\pd v_a}{\pd y}}\d x\,\d y
\le K_6^\al\int_{\pd S}J(a,\phi)\d u_a,
\end{equation*}
using equation \eq{us6eq5} and~$v_a\vert_{\pd S}=\phi$.

Now from $J(a,v)=-\int_0^v(w^2+a^2)^{-\al}\d w$ we see that
$J(a,0)=0$ and $\bmd{\frac{\pd J}{\pd v}(a,v)}\le\md{v}^{-2\al}$, so
integrating yields $\bmd{J(a,v)}\le(1-2\al)^{-1}\md{v}^{1-2\al}$.
Putting $v=\phi$ in this and estimating $\d u_a$ using \eq{us6eq3}
then gives
\e
\int_{\pd S}J(a,\phi)\d u_a\le(1-2\al)^{-1}
\cnm{\phi}{0}^{1-2\al}K_2\,l(\pd S),
\label{us6eq6}
\e
where $l(\pd S)$ is the length of $\pd S$. Combining the last two
equations and taking $p^{\rm th}$ roots gives $\lnm{\frac{\pd
v_a}{\pd y}}{p}\le C_p$ with~$C_p=\bigl(K_6^\al(1\!-\!2\al)^{-1}
\cnm{\phi}{0}^{1-2\al}K_2\,l(\pd S)\bigr)^{1/p}$.

To prove the second inequality, let $q\in[1,2)$ and put
$\al=\frac{3}{4}-\frac{1}{2q}$, a different value of $\al$. Then
$\al\in[\frac{1}{4},\ha)$, and using \eq{us2eq2} we see that
\begin{equation*}
\ts(v_a^2+y^2+a^2)^{1/2-\al}\bms{\frac{\pd u_a}{\pd y}}\le\frac{1}{4}
(v_a^2+a^2)^{-\al}(v_a^2+y^2+a^2)^{-1/2}\bms{\frac{\pd v_a}{\pd x}}.
\end{equation*}
Integrating this over $S$ and using \eq{us6eq5} and \eq{us6eq6} shows that
\e
\begin{split}
\int_S(v_a^2+y^2+a^2)^{1/2-\al}\bms{\ts\frac{\pd u_a}{\pd y}}\d x\,\d y
&\le{\ts\frac{1}{4}}\int_{\pd S}J(a,\phi)\d u_a\\
&\le{\ts\frac{1}{4}}(1-2\al)^{-1}\cnm{\phi}{0}^{1-2\al}K_2\,l(\pd S).
\end{split}
\label{us6eq7}
\e
Also $\md{y}^{-1/2}\in L^1(S)$, and therefore
\e
\int_S(v_a^2+y^2+a^2)^{-1/4}\d x\,\d y\le\int_S\md{y}^{-1/2}\d x\,\d y
=\blnm{\md{y}^{-1/2}}{1}.
\label{us6eq8}
\e

If $f,g\in L^1(S)$ and $\be\in(0,1)$ then H\"older's inequality gives
\begin{equation*}
\int_S\md{f}^\be\md{g}^{1-\be}\d x\,\d y\le
\Bigl(\int_S\md{f}\d x\,\d y\Bigr)^\be\cdot
\Bigl(\int_S\md{g}\d x\,\d y\Bigr)^{1-\be}.
\end{equation*}
Applying this to equations \eq{us6eq7} and \eq{us6eq8} with
$\be=q/2$ yields
\begin{equation*}
\blnm{\ts\frac{\pd u_a}{\pd y}}{q}^q
\le\bigl({\ts\frac{1}{4}}(1-2\al)^{-1}\cnm{\phi}{0}^{1-2\al}
K_2\,l(\pd S)\bigr)^{q/2}\blnm{\md{y}^{-1/2}}{1}^{1-q/2},
\end{equation*}
as the powers of $(v_a^2+y^2+a^2)$ cancel. Defining $C_q'$ to be the
$q^{\rm th}$ root of the r.h.s.\ then gives the second inequality.
Finally, \eq{us6eq3} implies that $\lnm{\frac{\pd v_a}{\pd x}}{r}\le
K_4\vol(S)^{1/r}=C''_r$ for all $r\in[1,\iy)$ and~$a\in(0,1]$.
\end{proof}

We use the $L^p$ bounds for $\frac{\pd u_a}{\pd x},\frac{\pd u_a}{\pd y},
\frac{\pd v_a}{\pd x}$ and $\frac{\pd v_a}{\pd y}$ to show that the $u_a,
v_a$ are {\it uniformly continuous} uniformly in $a\in(0,1]$. This will
imply that the limits $u_0,v_0$ of $u_a,v_a$ as $a\ra 0_+$ are continuous.

\begin{prop} In the situation of Definition \ref{us6def1}, there
exist continuous functions $M,N:S\t S\ra[0,\iy)$ satisfying
$M(x',y',x,y)=M(x,y,x',y')$, $N(x',y',x,y)=N(x,y,x',y')$ and\/
$M(x,y,x,y)=N(x,y,x,y)=0$ for all\/ $(x,y),(x',y')\in S$, such
that for all\/ $(x,y),(x',y')\in S$ and\/ $a\in(0,1]$ we have
\e
\begin{split}
\bmd{u_a(x,y)-u_a(x',y')}&\le M(x,y,x',y')\\
\text{and\/}\quad
\bmd{v_a(x,y)-v_a(x',y')}&\le N(x,y,x',y').
\end{split}
\label{us6eq14}
\e
\label{us6prop7}
\end{prop}

\begin{proof} Choose $p\in[2,\frac{5}{2})$ and $q\in[1,2)$ with
$p^{-1}+q^{-1}<1$. Then Theorem \ref{us2thm} gives continuous
functions $G,H:S\t S\ra[0,\iy)$. Combining Proposition \ref{us6prop6}
and Theorem \ref{us2thm}, we see that for all $(x,y),(x',y')\in S$
and $a\in(0,1]$, we have
\begin{equation*}
\bmd{u_a(x,y)-u_a(x',y')}\le C_pG(x,y,x',y')+C_q'H(x,y,x',y')
=M(x,y,x',y').
\end{equation*}

Similarly, choosing $\hat p\in[1,\iy)$ and $\hat q\in[2,\frac{5}{2})$
with $\hat p^{-1}+\hat q^{-1}<1$, for all $(x,y),(x',y')\in S$ and
$a\in(0,1]$ Proposition \ref{us6prop6} and Theorem \ref{us2thm} give
\begin{equation*}
\bmd{v_a(x,y)-v_a(x',y')}\le C_{\hat p}''\hat G(x,y,x',y')
+C_{\hat q}\hat H(x,y,x',y')=N(x,y,x',y'),
\end{equation*}
where $\hat G,\hat H$ are the continuous functions $S\t S\ra[0,\iy)$
given by Theorem \ref{us2thm}, starting from $\hat p,\hat q$ rather
than $p,q$. The conditions $M(x',y',x,y)=M(x,y,x',y')$, $N(x',y',x,y)=
N(x,y,x',y')$ and $M(x,y,x,y)=N(x,y,x,y)=0$ follow from the corresponding
conditions on $G,H$ in Theorem~\ref{us2thm}.
\end{proof}

\subsection{Existence of weak solutions when $a=0$}
\label{us66}

We can now construct a weak solution $v_0$ of the Dirichlet problem
for \eq{us3eq8} when~$a=0$.

\begin{prop} Let\/ $(a_n)_{n=1}^\iy$ be a sequence in $(0,1]$
with\/ $a_n\ra 0$ as $n\ra\iy$. Then there exists a subsequence
$(a_{n_i})_{i=1}^\iy$ and\/ $f_0\in C^1(S)$ such that\/
$f_{a_{n_i}}\ra f_0$ in $C^1(S)$ as $i\ra\iy$. Let\/
$u_0=\frac{\pd f_0}{\pd y}$ and\/ $v_0=\frac{\pd f_0}{\pd x}$.
Then $u_{a_{n_i}}\ra u_0$ and\/ $v_{a_{n_i}}\ra v_0$ in
$C^0(S)$ as $i\ra\iy$, and\/~$v_0\vert_{\pd S}=\phi$.
\label{us6prop8}
\end{prop}

\begin{proof} By Ascoli's Theorem \cite[Th.~3.15]{Aubi}, the
inclusion $C^1(S)\hookra C^0(S)$ is {\it compact}. But Corollary
\ref{us6cor2} gives $\cnm{f_{a_n}}{1}\le K_5$ for all $n$. Thus
the $f_{a_n}$ all lie in a compact subset of $C^0(S)$, and there
exists a subsequence $(a_{n_i})_{i=1}^\iy$ such that $f_{a_{n_i}}
\ra f_0$ in $C^0(S)$ for some $f_0\in C^0(S)$ as~$i\ra\iy$.

Define $A=\bigl\{u_a:a\in(0,1]\bigr\}$, and regard $A$ as a subset
of $C^0(S)$. Then $A$ is {\it bounded}, as $\cnm{u_a}{0}\le K_1$ by
\eq{us6eq3}, and {\it equicontinuous}, by Proposition \ref{us6prop7}.
Hence by Ascoli's Theorem \cite[Th.~3.15]{Aubi} $A$ is precompact in
$C^0(S)$. So the sequence $(u_{a_{n_i}})_{i=1}^\iy$ in $A$ must have
a convergent subsequence $(u_{a_{\hat n_i}})_{i=1}^\iy$ in $C^0(S)$,
converging to some $u_0\in C^0(S)$. We shall show that
$u_0=\frac{\pd f_0}{\pd y}$ in~$S$.

Suppose $y_1<y_2$ and $x\in\R$ with $(x,y)\in S$ whenever
$y\in[y_1,y_2]$. Then
\begin{equation*}
f_{a_{\hat n_i}}(x,y_2)-f_{a_{\hat n_i}}(x,y_1)=
\int_{y_1}^{y_2}u_{a_{\hat n_i}}(x,y)\,\d y.
\end{equation*}
Let $i\ra\iy$ in this equation. As $(a_{\hat n_i})_{i=1}^\iy$ is a
subsequence of $(a_{n_i})_{i=1}^\iy$ and $f_{a_{n_i}}\ra f_0$ in
$C^0(S)$ as $i\ra\iy$, the left hand side converges to $f_0(x,y_2)
-f_0(x,y_1)$. Thus, as $u_{a_{\hat n_i}}$ converges uniformly to
$u_0$ in $S$, we have
\begin{equation*}
f_0(x,y_2)-f_0(x,y_2)=\int_{y_1}^{y_2}u_0(x,y)\,\d y.
\end{equation*}
Differentiating this equation with respect to $y_2$ and using continuity
of $u_0$ shows that $\frac{\pd f_0}{\pd y}$ exists in $S$ and equals~$u_0$.

Thus the limit $u_0$ is unique, and so the whole sequence
$(u_{a_{n_i}})_{i=1}^\iy$ converges in $C^0(S)$ to $u_0$, rather than
just some subsequence $(u_{a_{\hat n_i}})_{i=1}^\iy$, as otherwise we
could choose a subsequence converging to a different limit. Using the
same argument we prove that $v_0=\frac{\pd f_0}{\pd x}$ exists in
$S$ and lies in $C^0(S)$, and $v_{a_{n_i}}\ra v_0$ in $C^0(S)$. As
$v_{a_{n_i}}\vert_{\pd S}=\phi$ it follows that~$v_0\vert_{\pd S}=\phi$.

Since $f_0\in C^0(S)$ and $\frac{\pd f_0}{\pd x},\frac{\pd f_0}{\pd y}$
exist and lie in $C^0(S)$ we have $f_0\in C^1(S)$, and as $f_{a_{n_i}}
\ra f_0$, $\frac{\pd}{\pd x}f_{a_{n_i}}\ra\frac{\pd}{\pd x}f_0$ and
$\frac{\pd}{\pd y}f_{a_{n_i}}\ra\frac{\pd}{\pd y}f_0$ in $C^0(S)$
as $i\ra\iy$, we have $f_{a_{n_i}}\ra f_0$ in $C^1(S)$ as $i\ra\iy$.
This completes the proof.
\end{proof}

In the next result, when we say that $u_0,v_0$ satisfy \eq{us3eq3}
and $f_0$ satisfies \eq{us3eq6} {\it with weak derivatives}, we mean
that the corresponding derivatives exist weakly and satisfy the equations.

\begin{prop} The function $v_0$ satisfies \eq{us3eq8} weakly with\/
$a=0$. The derivatives $\frac{\pd u_0}{\pd x}$, $\frac{\pd u_0}{\pd y}$,
$\frac{\pd v}{\pd x}$, $\frac{\pd v}{\pd y}$ exist weakly and
satisfy \eq{us3eq3}, with\/ $\frac{\pd u_0}{\pd x}=\frac{\pd
v_0}{\pd y}\in L^p(S)$ for $p\in[1,\frac{5}{2})$, and\/
$\frac{\pd u_0}{\pd y}\in L^q(S)$ for $q\in[1,2)$, and\/
$\frac{\pd v_0}{\pd x}$ bounded on $S$. The function $f_0$
satisfies \eq{us3eq6} with weak derivatives and\/~$a=0$.
\label{us6prop9}
\end{prop}

\begin{proof} We shall show that $v_0$ satisfies \eq{us3eq9}
weakly, which is equivalent to \eq{us3eq8}. Let $A(a,y,v)$ be as
in \eq{us3eq7}. Then by \eq{us3eq4}, for all $a\in(0,1]$ we have
\begin{equation*}
\frac{\pd}{\pd x}\bigl(A(a,y,v_a)\bigr)=
(v_a^2+y^2+a^2)^{-1/2}\frac{\pd v_a}{\pd x}=2\frac{\pd u_a}{\pd y}.
\end{equation*}
Thus Proposition \ref{us6prop6} gives $\lnm{\frac{\pd}{\pd x}
A(a,y,v_a)}{q}\le 2C_q'$ and $\lnm{\frac{\pd v_a}{\pd y}}{p}\le C_p$
for $p\in[2,\frac{5}{2})$, $q\in(1,2)$ and all $a\in(0,1]$.

The inclusions $L^p(S)\hookra L^1(S)$ and $L^q(S)\hookra
L^1(S)$ are {\it compact\/} by Aubin \cite[Th.~2.33]{Aubi}.
Thus the functions $\frac{\pd}{\pd x}A(a,y,v_a)$ for all
$a\in(0,1]$ and $\frac{\pd v_a}{\pd y}$ for all $a\in(0,1]$
lie in compact subsets of $L^1(S)$. So we may choose
subsequences of the sequences $(\frac{\pd}{\pd x}A(a_{n_i},
y,v_{a_{n_i}}))_{i=1}^\iy$ and $(\frac{\pd v_{a_{n_i}}}{\pd y}
)_{i=1}^\iy$ which converge in~$L^1(S)$.

By an argument similar to that in the proof of Proposition
\ref{us6prop8}, we can show that these limits are the weak derivatives
$\frac{\pd}{\pd x}A(0,y,v_0)$ and $\frac{\pd}{\pd y}v_0$. So the limits
are unique, and the whole sequence converges to them, not just a
subsequence. That is, $\frac{\pd}{\pd x}A(0,y,v_0)$ and
$\frac{\pd}{\pd y}v_0$ both exist weakly in $L^1(S)$, and
$\frac{\pd}{\pd x}A(a_{n_i},y,v_{a_{n_i}})\ra\frac{\pd}{\pd x}
A(0,y,v_0)$ and $\frac{\pd}{\pd y} v_{a_{n_i}}\ra\frac{\pd}{\pd y}v_0$
in $L^1(S)$ as~$i\ra\iy$.

Let $\psi\in C^1_0(S)$. As $v_{a_{n_i}}$ and $a_{n_i}$ satisfy
\eq{us3eq9}, multiplying by $\psi$ and integrating by parts gives
\begin{equation*}
-\int_S\frac{\pd\psi}{\pd x}\cdot\frac{\pd}{\pd x}A(a_{n_i},y,
v_{a_{n_i}})\,\d x\,\d y-2\int_S\frac{\pd\psi}{\pd y}\cdot
\frac{\pd v_{a_{n_i}}}{\pd y}\,\d x\,\d y=0,
\end{equation*}
Letting $i\ra\iy$ in this equation shows that
\begin{equation*}
-\int_S\frac{\pd\psi}{\pd x}\cdot\frac{\pd}{\pd x}A(0,y,v_0)\,\d x\,\d y
-2\int_S\frac{\pd\psi}{\pd y}\cdot\frac{\pd v_0}{\pd y}\,\d x\,\d y=0,
\end{equation*}
since $\frac{\pd}{\pd x}A(a_{n_i},y,v_{a_{n_i}})\ra\frac{\pd}{\pd x}
A(0,y,v_0)$ and $\frac{\pd}{\pd y}v_{a_{n_i}}\ra\frac{\pd}{\pd y}v_0$
in $L^1(S)$ as $i\ra\iy$, and $\frac{\pd\psi}{\pd x},\frac{\pd\psi}{\pd
y}\in C^0(S)$. Thus $v_0$ satisfies \eq{us3eq9} weakly with~$a=0$.

The proof above shows that $\frac{\pd v_0}{\pd y}$ exists weakly and
is the limit in $L^1(S)$ of $\frac{\pd}{\pd y}v_{a_{n_i}}$ as $i\ra\iy$.
Since $\frac{\pd}{\pd y}v_{a_{n_i}}=\frac{\pd}{\pd x}u_{a_{n_i}}$, it
easily follows that $\frac{\pd u_0}{\pd x}$ exists weakly and is the
limit in $L^1(S)$ of $\frac{\pd}{\pd x}u_{a_{n_i}}$ as $i\ra\iy$. Also,
as $\frac{\pd}{\pd x}A(a_{n_i},y,v_{a_{n_i}})=-2\frac{\pd}{\pd y}
u_{a_{n_i}}$, a similar argument shows that $\frac{\pd u_0}{\pd y}$
exists weakly and is the limit in $L^1(S)$ of $\frac{\pd}{\pd y}
u_{a_{n_i}}$ as~$i\ra\iy$.

Since $v_{a_{n_i}}^2+y^2+a_{n_i}^2$ is bounded above we can use this
and \eq{us3eq4} to deduce that $\frac{\pd v_0}{\pd x}$ exists weakly
and is the limit in $L^1(S)$ of $\frac{\pd}{\pd x}v_{a_{n_i}}$ as
$i\ra\iy$. Thus, the first derivatives of $u_0,v_0$ and the second
derivatives of $f_0$ exist weakly, and are the limits in $L^1(S)$ of
the corresponding derivatives of $u_{a_{n_i}},v_{a_{n_i}},f_{a_{n_i}}$
as $i\ra\iy$. The estimates in Proposition \ref{us6prop6} and
Corollary \ref{us6cor2} then imply that $\frac{\pd u_0}{\pd x}=
\frac{\pd v_0}{\pd y}\in L^p(S)$ for $p\in[1,\frac{5}{2})$,
$\frac{\pd u_0}{\pd y}\in L^q(S)$ for $q\in[1,2)$, and
$\frac{\pd v_0}{\pd x}$ is bounded on $S$. Taking the limit in
$L^1(S)$ as $i\ra\iy$ of equations \eq{us3eq4} and \eq{us3eq6}
then completes the proof.
\end{proof}

Away from singular points $(x,0)$ with $v_0(x,0)=0$ we can prove
much stronger regularity of $u_0,v_0$. Near a nonsingular point
$(x,y)$, equation \eq{us3eq5} with $a=0$ is strictly elliptic.
We use results of \cite[\S 6]{GiTr} to show that $f_0$ is
$C^{2,\al}$ near $(x,y)$ and satisfies \eq{us3eq5} with $a=0$,
and then \cite[Th.~3.6]{Joyc3} gives:

\begin{prop} Except at singular points $(x,0)$ with\/ $v_0(x,0)=0$,
the functions $u_0,v_0$ are $C^{k+2,\al}$ in $S$ and real analytic
in~$S^\circ$.
\label{us6prop10}
\end{prop}

\subsection{Uniqueness of weak solutions when $a=0$}
\label{us67}

Weak solutions of the Dirichlet problem for \eq{us3eq8} are unique.

\begin{prop} Let\/ $v,v'\in C^0(S)\cap L^1_1(S)$ be weak solutions of\/
\eq{us3eq8} on $S$ with\/ $a=0$ and\/ $v\vert_{\pd S}=v'\vert_{\pd S}
=\phi$. Then~$v=v'$.
\label{us6prop11}
\end{prop}

\begin{proof} Following Proposition \ref{us3prop3} but using weak
solutions we find that there exist $u,u'\in L^1_1(S)$ such that $u,v$
and $u',v'$ satisfy \eq{us3eq3} with weak derivatives. Using the ideas
of \S\ref{us64} we can show that $u,u',v,v'$ are $C^{k+2,\al}$ near
$\pd S$, so $u,u',v,v'$ are bounded. From Proposition \ref{us3prop2}
we see that there exist $f,f'\in C^{0,1}(S)\cap L^2_1(S)$ with
$\frac{\pd f}{\pd y}=u$, $\frac{\pd f}{\pd x}=v$, $\frac{\pd f'}{\pd
y}=u'$ and $\frac{\pd f'}{\pd x}=v'$ weakly, that satisfy \eq{us3eq6}
with weak derivatives, and $f,f'$ are $C^{k+3,\al}$ near~$\pd S$.

Let $\ga\in\R$. Then as in \cite[Prop.~7.5]{Joyc3} we find that
$f-f'+\ga y$ weakly satisfies an equation $L(f-f'+\ga y)=0$, where
$L$ is a linear elliptic operator of the form \eq{us2eq1} with
$c\equiv 0$. Thus by the maximum principle for elliptic operators
\cite[Th.~3.1]{GiTr}, which holds for weak solutions by
\cite[p.~45-6]{GiTr}, the maximum and minimum of $f-f'+\ga y$ occur
on $\pd S$. Furthermore, one can use \cite[Lem.~3.4]{GiTr} to show
that either the normal derivatives of $f-f'+\ga y$ at the maximum
and minimum are nonzero, or else $f-f'+\ga y$ is constant in~$S$.

Let the maximum of $f-f'+\ga y$ occur at $(x,y)\in\pd S$. Then
\begin{equation*}
\ts\frac{\pd}{\pd x}\bigl(f-f'+\ga y\bigr)\big\vert_{(x,y)}=
v(x,y)-v'(x,y)=\phi(x,y)-\phi(x,y)=0.
\end{equation*}
But the derivative of $f-f'+\ga y$ at $(x,y)$ tangent to $\pd S$
is also zero. Thus, if $\frac{\pd}{\pd x}$ is not tangent to
$\pd S$ at $(x,y)$ then the normal derivative of $f-f'+\ga y$ at
$(x,y)$ is zero, and so $f-f'+\ga y$ is constant.

Therefore either $f-f'+\ga y$ is constant, or else the maximum
(and similarly the minimum) of $f-f'+\ga y$ occur at points $(x,y)$
in $\pd S$ where $\frac{\pd}{\pd x}$ is tangent to $\pd S$. But as
$S$ is strongly convex there are only two points ${\bf x}_1,{\bf x}_2$
in $\pd S$ with $\frac{\pd}{\pd x}$ tangent to $\pd S$. Choose $\ga\in\R$
uniquely so that $\frac{\pd}{\pd y}\bigl(f-f'+\ga y\bigr)=0$ at ${\bf x}_1$.
Then either $f-f'+\ga y$ is constant, or else the maximum and minimum of
$f-f'+\ga y$ both occur at ${\bf x}_2$, again implying that $f-f'+\ga y$
is constant. Taking $\frac{\pd}{\pd x}$ gives~$v=v'$.
\end{proof}

This implies that the limit $v_0$ chosen in Proposition \ref{us6prop8}
is {\it unique}. Thus, the entire sequence $(v_{a_n})_{n=1}^\iy$
converges to $v_0$ in $C^0(S)$ rather than just the subsequence
$(v_{a_{n_i}})_{i=1}^\iy$, since otherwise we could have chosen a
different limit $v_0$. As this is true for an arbitrary sequence
$(a_n)_{n=1}^\iy$ in $(0,1]$, this shows that $v_a\ra v_0$ in
$C^0(S)$ as $a\ra 0_+$. Similarly, the limits $u_0,f_0$ are also
unique, as the freedoms to add constants are fixed by $u_0(x_1,0)=
f_0(x_1,0)=0$. So we deduce:

\begin{cor} As $a\ra 0_+$ in $(0,1]$ we have $u_a\ra u_0$, $v_a\ra v_0$
in $C^0(S)$ and\/ $f_a\ra f_0$ in~$C^1(S)$.
\label{us6cor3}
\end{cor}

Theorem \ref{us6thm1} now follows from Propositions
\ref{us6prop8}--\ref{us6prop11}. For Theorem \ref{us6thm2},
continuity in the $a$ variable at $a=0$ for fixed $\phi$ follows from
Corollary \ref{us6cor3}, since $u_a\ra u_0$ and $v_a\ra v_0$ in $C^0(S)$
as $a\ra 0_+$. To prove continuity in the $\phi$ variable we need to show
that small $C^{k+2,\al}$ changes in $\phi$ result in small $C^0$ changes
in $u,v$, and this can be seen by examining the proofs above.

\section{The Dirichlet problem for $f$ when $a=0$}
\label{us7}

Theorem \ref{us3thm1} shows that the Dirichlet problem for equation
\eq{us3eq5} is uniquely solvable in strictly convex domains for
$a\ne 0$. In this section we will use the material of \S\ref{us5}
and \S\ref{us6} to show that the Dirichlet problem also has a unique
solution when $a=0$, but with weak second derivatives.

\subsection{The main results}
\label{us71}

Here is an analogue of Theorem \ref{us3thm1} in the case~$a=0$.

\begin{thm} Let\/ $S$ be a strictly convex domain in $\R^2$ invariant
under the involution $(x,y)\mapsto(x,-y)$, let\/ $k\ge 0$ and\/ $\al\in
(0,1)$. Then for each\/ $\phi\in C^{k+3,\al}(\pd S)$ there exists a
unique weak solution $f$ of\/ \eq{us3eq6} in $C^1(S)$ with\/ $a=0$ and\/
$f\vert_{\pd S}=\phi$. Furthermore $f$ is twice weakly differentiable
and satisfies \eq{us3eq5} with weak derivatives.

Let\/ $u=\frac{\pd f}{\pd y}$ and\/ $v=\frac{\pd f}{\pd x}$. Then
$u,v\in C^0(S)$ are weakly differentiable and satisfy \eq{us3eq3}
with weak derivatives, and\/ $v$ satisfies \eq{us3eq8} weakly with\/
$a=0$. The weak derivatives $\frac{\pd u}{\pd x},\frac{\pd u}{\pd y},
\frac{\pd v}{\pd x},\frac{\pd v}{\pd y}$ satisfy $\frac{\pd u}{\pd x}=
\frac{\pd v}{\pd y}\in L^p(S)$ for $p\in[1,2]$, and\/ $\frac{\pd u}{\pd
y}\in L^q(S)$ for $q\in[1,2)$, and\/ $\frac{\pd v}{\pd x}$ is bounded
on $S$. Also $u,v$ are $C^{k+2,\al}$ in $S$ and real analytic in
$S^\circ$ except at singular points $(x,0)$ with\/~$v(x,0)=0$.
\label{us7thm1}
\end{thm}

Combined with Proposition \ref{us3prop1} the theorem can be used to
construct large numbers of {\it $\U(1)$-invariant singular special
Lagrangian $3$-folds} in $\C^3$. This is the principal motivation
for the paper. The singularities of these special Lagrangian 3-folds
will be studied in~\cite{Joyc4}.

Our second theorem extends \cite[Th.~7.7]{Joyc3} to include the
case~$a=0$.

\begin{thm} Let\/ $S$ be a strictly convex domain in $\R^2$ invariant
under the involution $(x,y)\mapsto(x,-y)$, let\/ $k\ge 0$ and\/
$\al\in(0,1)$. Then the map $C^{k+3,\al}(\pd S)\t\R\ra C^1(S)$ taking
$(\phi,a)\mapsto f$ is continuous, where $f$ is the unique solution
of\/ \eq{us3eq5} (with weak derivatives) with\/ $f\vert_{\pd S}=\phi$
constructed in Theorem \ref{us3thm1} when $a\ne 0$, and in Theorem
\ref{us7thm1} when $a=0$. This map is also continuous in stronger
topologies on $f$ than the $C^1$ topology.
\label{us7thm2}
\end{thm} 

The proofs of these theorems take up the rest of the paper,
and are similar to those of Theorems \ref{us6thm1} and
\ref{us6thm2}. Here is how they are laid out. Let $S,\phi$
be as in Theorem \ref{us7thm1}. In \S\ref{us72}, for each
$a\in(0,1]$ we define $f_a\in C^{k+3,\al}(S)$ to be the
unique solution of \eq{us3eq5} in $S$ with $f_a\vert_{\pd
S}=\phi$, and $u_a=\frac{\pd f_a}{\pd y}$, $v_a=\frac{\pd
f_a}{\pd x}$. As in \S\ref{us6}, we aim to show that
$u_a,v_a,f_a\ra u_0,v_0,f_0$ as $a\ra 0_+$, where $u_0,v_0,f_0$
are $u,v,f$ in Theorem \ref{us7thm1}, and the main issue is
to prove {\it a priori estimates} of $u_a,v_a,f_a$ that are
{\it uniform in}~$a$.

However, there are some important differences with \S\ref{us6}.
On the positive side, Theorem \ref{us3thm1} immediately gives
uniform bounds for $\cnm{f_a}1,\cnm{u_a}0$ and $\cnm{v_a}0$.
But as we have no analogue of the $\phi(x,0)\ne 0$ assumption
in Theorem \ref{us6thm1}, we must allow $v_a$ to be zero at
or near points $(x,0)$ in $\pd S$. Because of this, equation
\eq{us3eq5} at $f_a$ need {\it not\/} be uniformly elliptic in
$a$ close to $\pd S$, so the methods of \S\ref{us63}--\S\ref{us64}
for estimating $u_a,v_a,f_a$ near $\pd S$ do not work.

Instead, we do something different. The second derivatives
of $\phi$ give us a bound on half of $\pd u_a,\pd v_a$ on
$\pd S$. Using this, and supposing $\pd S$ positively
curved at points $(x,0)$, in \S\ref{us73} we use a
boundary version of the method of \S\ref{us5} to
estimate the other half of $\pd u_a,\pd v_a$ on $\pd S$
near points~$(x,0)$.

Section \ref{us74} then gives partial analogues of
the a priori estimates of \S\ref{us64}--\S\ref{us65}.
As our bounds on $\frac{\pd v_a}{\pd y}$ on $\pd S$
are not strong enough to apply the global estimates
Theorem \ref{us5thm4} uniformly in $a$, we instead
have to work in interior domains $T\subset S^\circ$
for our $L^p$ estimates of $\frac{\pd u_a}{\pd x},
\frac{\pd v_a}{\pd y}$ when~$p\in(2,\frac{5}{2})$.

In \S\ref{us75} we establish uniform continuity of
the $u_a,v_a$, as in \S\ref{us65}. But because we
only have interior $L^p$ estimates of $\frac{\pd
u_a}{\pd x}$, we must do some extra work to show
$u_a$ is uniformly continuous near points $(x,0)$
in $\pd S$. Finally, \S\ref{us76} proves existence
of the limit solutions $f_0,u_0,v_0$, which is
just as in \S\ref{us66}, and uniqueness, completing
the proofs.

The hypotheses of Theorems \ref{us7thm1} and \ref{us7thm2} can be
relaxed slightly, without changing the proofs: rather than requiring
$S$ invariant under $(x,y)\mapsto(x,-y)$ we can ask only that each
point $(x,0)$ in $\pd S$ has tangent $T_{(x,0)}\pd S$ parallel
to the $y$-axis, and rather than requiring $\phi\in C^{k+3,\al}(\pd S)$
we can ask only that~$\phi\in C^3(\pd S)$.

Note that these two theorems do not have the awkward restriction
that $\phi(x,0)\ne 0$ for points $(x,0)\in\pd S$ in Theorems
\ref{us6thm1} and \ref{us6thm2}. For this reason we find them
more convenient for applications such as constructing special
Lagrangian fibrations on subsets of $\C^3$, and we will generally
use them in preference to Theorems \ref{us6thm1} and \ref{us6thm2}
in the sequel~\cite{Joyc4}.

\subsection{A family of solutions $f_a$ to \eq{us3eq5}}
\label{us72}

Consider the following situation:

\begin{dfn} Let $S$ be a strictly convex domain in $\R^2$ which
is invariant under the involution $(x,y)\mapsto(x,-y)$. Then there
exist unique $x_1,x_2\in\R$ with $x_1<x_2$ and $(x_i,0)\in\pd S$
for $i=1,2$. Let $k\ge 0$ and $\al\in(0,1)$, and suppose $\phi\in
C^{k+3,\al}(\pd S)$. For each $a\in(0,1]$, let $f_a\in C^{k+3,\al}(S)$
be the unique solution of \eq{us3eq5} in $S$ with this value of $a$ and
$f_a\vert_{\pd S}=\phi$, which by Theorem \ref{us3thm1} exists and
satisfies $\cnm{f_a}{1}\le C\cnm{\phi}{2}$ for all $a\in(0,1]$, where
$C>0$ depends only on $S$, and in particular is independent of $a$. Set
\e
X=C\cnm{\phi}{2}\quad\text{and}\quad Y=\sup_{(x,y)\in S}\md{y}.
\label{us7eq1}
\e
Define $u_a,v_a\in C^{k+2,\al}(S)$ by $u_a=\frac{\pd f_a}{\pd y}$ and
$v_a=\frac{\pd f_a}{\pd x}$. Then
\e
\cnm{f_a}{1},\cnm{u_a}{0},\cnm{v_a}{0}\le X \quad\text{for all $a\in(0,1]$,}
\label{us7eq2}
\e
and $u_a,v_a$ and $a$ satisfy \eq{us3eq4}, and $v_a$
and $a$ satisfy~\eq{us3eq5}.
\label{us7def1}
\end{dfn}

We will show that $f_a$ converges in $C^1(S)$ to $f_0\in C^1(S)$ as
$a\ra 0_+$, and that $f_0$ is the unique weak solution of the Dirichlet
problem for \eq{us3eq5} on $S$ when $a=0$. The main difficulty in doing
this is to establish uniform continuity of $u_a$ and $v_a$ for all
$a\in(0,1]$, as we did for the $v$ Dirichlet problem in \S\ref{us65}.
Once we have done this, we can follow the proofs of \S\ref{us6} with
few changes.

First we bound the $f_a$ in $C^2$ away from the $x$-axis.

\begin{prop} Let\/ $\ep>0$ be small, and set\/ $S_\ep=\bigl\{
(x,y)\in S:\md{y}>\ep\bigr\}$. Then there exists $G>0$ such
that\/~$\cnm{f_a\vert_{S_\ep}}{2}\le G$ for all\/~$a\in(0,1]$.
\label{us7prop1}
\end{prop}

\begin{proof} Let $S_{\ep/2}$ and $S_{\ep/4}$ be defined
in the obvious way. We prove the proposition in two steps.
Regard $P$ as in \eq{us3eq5} as a {\it linear} operator,
with coefficients $a^{ij}$ depending on $v_a$. Firstly
we use estimates on $S_{\ep/4}$ to bound $\cnm{f_a
\vert_{S_{\ep/2}}}{1,\ga}$ uniformly in $a\in(0,1]$ for
some $\ga\in(0,1)$. This gives a uniform bound $\La'$ for
$\cnm{a^{ij}\vert_{S_{\ep/2}}}{0,\ga}$. Secondly, we bound
$\cnm{f_a\vert_{S_{\ep}}}{2,\ga}$ uniformly in~$a$.

The results we need are {\it interior regularity results}
for linear elliptic operators on {\it noncompact regions}
in $\R^2$ with {\it boundary portions}. In the first step
we use Gilbarg and Trudinger \cite[Th.~12.4, p.~302]{GiTr},
extended to the boundary case as in \cite[p.~303-4]{GiTr}.
This deals with equations of the form $Pu=f$, where $P$
is a linear elliptic operator of the form \eq{us2eq1} with
$a^{ij}\in C^0(S)$ and $b^i=c=0$. Note that \eq{us3eq5}
is of this form.

Now $\frac{1}{16}\ep^2\le v_a^2+y^2+a^2\le X^2+Y^2+1$ on $S_{\ep/4}$.
Raising this to the power $-\ha$ we see from \eq{us3eq5} that
$\sum_{i,j=1}^2a^{ij}\xi_i\xi_j\ge\la(\xi_1^2+\xi_2^2)$ for all
$\xi=(\xi_1,\xi_2)\in\R^2$ and $\cnm{a^{ij}}{0}\le\La$ hold in
$S_{\ep/4}$, with $\la=\min\bigl((X^2+Y^2+1)^{-1/2},2\bigr)$
and $\La=\max(4\ep^{-1},2)$. Therefore \cite[p.~302-4]{GiTr}
gives $\ga\in(0,1)$ depending only on $\La/\la$, and $G'>0$
such that $\cnm{f_a\vert_{S_{\ep/2}}}{1,\ga}\le G'$ for
all~$a\in(0,1]$.

This gives a uniform bound $\La'$ for $\cnm{a^{ij}
\vert_{S_{\ep/2}}}{0,\ga}$. Following \cite[Lem.~6.18]{GiTr},
which gives a priori $C^{2,\ga}$ interior estimates for solutions
of linear elliptic equations with $C^{0,\ga}$ coefficients
in domains with boundary portions, we find that there exists
$G>0$ depending on $\la,\La,\ga,G',S_{\ep/2}$ and $S_\ep$
such that $\cnm{f_a\vert_{S_{\ep}}}{2,\ga}\le G$ for all
$a\in(0,1]$. Hence $\cnm{f_a\vert_{S_{\ep}}}{2}\le G$ for
all $a$, as we want.
\end{proof}

Thus, $\pd^2f_a$ is uniformly bounded on $\pd S$ except arbitrarily
close to the points $(x_i,0)$ for $i=1,2$. So we shall study $u_a,v_a$
and $f_a$ near these points. Fix $i=1$ or 2. As $S$ is invariant under
$(x,y)\mapsto(x,-y)$, the tangent to $\pd S$ at $(x_i,0)$ is parallel
to the $y$-axis. Thus we may use $y$ as a parameter on $\pd S$ near
$(x_i,0)$, and write $(x,y)\in\pd S$ near $(x_i,0)$ as $\bigl(x(y),y
\bigr)$. Regard $\phi$ as a function of $y$ near $(x_i,0)$. Then
differentiating the equation $f_a\bigl(x(y),y\bigr)=\phi(y)$ once
and twice w.r.t.\ $y$ gives
\begin{align}
\dot\phi&=u_a(x,y)+v_a(x,y)\dot x \qquad\text{and}
\label{us7eq3}\\
\begin{split}
\ddot\phi&=
\frac{\pd u_a}{\pd x}\,\dot x+\frac{\pd u_a}{\pd y}+
\frac{\pd v_a}{\pd x}\,\dot x^2+\frac{\pd v_a}{\pd y}\,\dot x+
v_a\,\ddot x\\
&=\frac{\pd u_a}{\pd y}\bigl(1-2(v_a^2+y^2+a^2)^{1/2}\dot x^2\bigr)
+2\frac{\pd u_a}{\pd x}\,\dot x+v_a\,\ddot x,
\end{split}
\label{us7eq4}
\end{align}
writing `$\,\dot{}\,$' for $\frac{\d}{\d y}$, and using \eq{us3eq4}
in the final line.

Now when $y=0$ we have $\dot x=0$ and $\ddot x=\ka_i$, where
$\ka_i$ is the curvature of $\pd S$ at $(x_i,0)$, measured in the
direction of increasing $x$. Thus $\dot x\approx\ka_iy$ to
leading order in $y$ near $(x_i,0)$. As $S$ is strictly convex,
it follows that $\ka_1>0$ and $\ka_2<0$. We can use \eq{us7eq4}
to prove:

\begin{prop} There exist\/ $\ep,H>0$ and smooth functions
$F_{i,a}:[-\ep,\ep]\ra\R$ for $i=1,2$ and $a\in(0,1]$ with\/
$\bmd{F_{i,a}(y)-2\ka_iy}\le\md{\ka_iy}$, such that for all\/
$(x,y)\in\pd S$ close to $(x_i,0)$ and with\/ $\md{y}\le\ep$, we have
\e
\Big\vert\frac{\pd u_a}{\pd y}(x,y)+F_{i,a}(y)\frac{\pd u_a}{\pd x}
(x,y)\Big\vert<H\quad\text{for all\/ $a\in(0,1]$.}
\label{us7eq5}
\e
\label{us7prop2}
\end{prop}

Here $F_{i,a}(y)=2\dot x\bigl(1-2(v_a^2+y^2+a^2)^{1/2}\dot x^2\bigr)^{-1}$,
so it does depend on $a$ and $v_a$. But the approximations above give
$F_{i,a}(y)\approx 2\ka_iy$ for small $y$, and one can show that for
$\ep>0$ depending only on $S$ and upper bounds $1,X$ for $\md{a},\md{v_a}$,
if $\md{y}\le\ep$ then~$\bmd{F_{i,a}(y)-2\ka_iy}\le\md{\ka_iy}$.

\subsection{An a priori bound for $\frac{\pd u_a}{\pd x},
\frac{\pd v_a}{\pd y}$ on $\pd S$}
\label{us73}

Proposition \ref{us7prop1} bounds $\pd u_a,\pd v_a$ on $\pd S$
away from the $x$-axis, and Proposition \ref{us7prop2} in effect
bounds half of $\pd u_a$, and hence $\pd v_a$, on $\pd S$ near
the $a$-axis. The following theorem in effect bounds the other
half of $\pd u_a$ and $\pd v_a$ near the $x$-axis. The proof,
which is unfortunately rather long and complicated, adapts the
method of~\S\ref{us52}.

\begin{thm} There exist\/ $\de,J>0$ such that for all\/ $(x_0,y_0)
\in\pd S$ close to $(x_i,0)$ for $i=1$ or $2$ and with\/
$\md{y_0}\le\de$, we have
\e
\Big\vert\frac{\pd u_a}{\pd x}(x_0,y_0)\Big\vert=
\Big\vert\frac{\pd v_a}{\pd y}(x_0,y_0)\Big\vert\le
J\bigl(y_0^2+a^2\bigr)^{-1/2}\quad\text{for all\/ $a\in(0,1]$.}
\label{us7eq6}
\e
\label{us7thm3}
\end{thm}

\begin{proof} Choose small $\de>0$ and large $J>0$, to satisfy
conditions we will give during the proof. Suppose, for a contradiction,
that $(x_0,y_0)\in\pd S$ is close to $(x_i,0)$ for $i=1$ or $2$ with
$\md{y_0}\le\de$, and that \eq{us7eq6} does not hold for some given
$a\in(0,1]$. Let $\ep$ be as in Proposition \ref{us7prop2}, and
suppose $\de\le\ep$. Then Proposition \ref{us7prop2} implies that
\eq{us7eq5} holds at $(x_0,y_0)$. Define
\begin{equation*}
u_0=u_a(x_0,y_0),\;\> v_0=v_a(x_0,y_0), \;\>
p_0=\ts\frac{\pd v_a}{\pd x}(x_0,y_0)\quad\text{and}\quad
q_0=\ts\frac{\pd v_a}{\pd y}(x_0,y_0).
\end{equation*}
Then as \eq{us7eq6} does not hold we have 
\e
\md{q_0}>J(y_0^2+a^2)^{-1/2},
\label{us7eq7}
\e
and from \eq{us3eq4} and Proposition \ref{us7prop2} we deduce that
\e
\begin{gathered}
\md{p_0}<2(v_0^2+y_0^2+a^2)^{1/2}\bmd{F_{i,a}(y_0)}\cdot\md{q_0}
+2(v_0^2+y_0^2+a^2)^{1/2}H,\\
\text{where}\quad \bmd{F_{i,a}(y_0)-2\ka_iy_0}\le\md{\ka_iy_0}.
\end{gathered}
\label{us7eq8}
\e
These imply that $q_0$ is large, and that $p_0$ is small compared
to~$q_0$.

Set $L=10X$, and define $\hat x_0,\hat y_0,\hat u_0,\hat v_0,\hat p_0,
\hat q_0$ as in \eq{us5eq11}. Then as in the proof of Proposition
\ref{us5prop2}, we can show that if $\de>0$ is small enough and $J>0$
large enough there exist $\hat u,\hat v\in C^\iy(D_L)$ satisfying
\eq{us3eq4} and \eq{us4eq2}, and with $(\hat u-\hat u_0,\hat v-\hat v_0)$
bounded by a small constant. Furthermore, as $q_0$ is large and $p_0$
small compared to $q_0$, we find that $\hat u,\hat v$ approximate the
affine maps \eq{us5eq18} up to their first derivatives.

Define $U=(\hat u,\hat v)(D_L)$. Then as in the proof of Proposition
\ref{us5prop2}, $U$ is approximately a closed disc of radius $\md{\hat
q_0}L$, and $(\hat u,\hat v):D_L\ra U$ is invertible with differentiable
inverse $(u',v'):U\ra D_L$. Moreover, $u',v'$ satisfy \eq{us3eq4} in $U$,
and by construction we have $u_a=u'$, $v_a=v'$, $\pd u_a=\pd u'$ and
$\pd v_a=\pd v'$ at $(x_0,y_0)$. Thus $(u',v')-(u_a,v_a)$ has a zero
of {\it multiplicity} at least 2 at $(x_0,y_0)$, in the sense
of~\cite[Def.~6.3]{Joyc3}.

To complete the proof, we follow a similar strategy to Proposition
\ref{us5prop2}. Roughly speaking, we shall show that the winding
number of $(u',v')-(u_a,v_a)$ about 0 along the boundary of $U\cap S$
is at most 1. But this contradicts \cite[Th.~6.7]{Joyc3}, as the
number of zeroes of $(u',v')-(u_a,v_a)$ in $U\cap S$ counted with
multiplicity should be at most 1, but we already know there is a zero
of multiplicity at least 2 at $(x_0,y_0)$. This then proves the theorem.

There are two problems with this strategy. The first is that
$(x_0,y_0)$ lies on the boundary of $U\cap S$ rather than the
interior, and so the winding number of $(u',v')-(u_a,v_a)$ about
0 along $\pd(U\cap S)$ is not defined, and \cite[Th.~6.7]{Joyc3}
does not apply. To deal with this we perturb $(x_0,y_0)$ a very
little way into the interior of $U\cap S$, and construct a
slightly different $(u',v')$ to intersect $(u_a,v_a)$ with
multiplicity 2 at the new point $(x_0,y_0)$ instead. Since
$u_a,v_a$ are $C^1$ and \eq{us7eq7} and \eq{us7eq8} are open
conditions, we can still assume that \eq{us7eq7} and \eq{us7eq8}
hold at the new~$(x_0,y_0)$.

The second problem is how to prove that the winding number of
$(u',v')-(u_a,v_a)$ about 0 along $\pd(U\cap S)$ is at most 1,
given that we do not know much about the behaviour of $(u_a,v_a)$.
We shall use the method of \cite[Th.~7.10]{Joyc3}, which bounds the
number of zeroes of $(u_1,v_1)-(u_2,v_2)$ in $S$ in terms of the
stationary points of the difference $f_1-f_2$ of their potentials.
Here is the crucial step in the proof.

\begin{prop} Let\/ $f'\in C^\iy(U)$ be a potential for $u',v'$,
as in Proposition \ref{us3prop2}. Then $f'-f_a$ has at most two
stationary points on the curve~$U\cap\pd S$.
\label{us7prop3}
\end{prop}

\begin{proof} From above, $\hat u,\hat v$ approximate the affine
maps in \eq{us5eq18} up to their first derivatives. Inverting this,
we find that $u',v'$ also approximate the affine maps
\e
u'(x,y)\approx u_0+q_0(x-x_0),\qquad
v'(x,y)\approx v_0+q_0(y-y_0)
\label{us7eq9}
\e
up to their first derivatives. Hence the potential $f'$ approximates
the quadratic
\e
f'(x,y)\approx u_0(y-y_0)+v_0(x-x_0)+q_0(x-x_0)(y-y_0)+c
\label{us7eq10}
\e
up to its second derivatives, for some $c\in\R$. We are being vague
about what we mean by `approximates' here. An exact statement can
be derived by using Theorem \ref{us4thm3} to estimate $(\hat u,\hat v)$
and its derivatives, and then inverting.

As in \S\ref{us72}, we may parametrize $\pd S$ near $(x_i,0)$ as
$\bigl(x(y),y\bigr)$ with $x(y)\approx x_i+\ha\ka_iy^2$ for small
$y$, where $\ka_i\ne 0$ is the curvature of $\pd S$ at $(x_i,0)$.
It follows that the restriction of $f'$ to $\pd S$ is approximately
\e
f'\bigl(x(y),y\bigr)\approx \ha\ka_iq_0y^3+
\ha\ka_i(v_0-q_0y_0)y^2+(u_0+q_0x_i-q_0x_0)y+c',
\label{us7eq11}
\e
where $c'=c-u_0y_0+(q_0+v_0)(x_i-x_0)$. By making $J$ large enough
we can suppose that $\md{q_0}\gg 2\md{\ka_i}^{-1}$ by \eq{us7eq7}.
Thus $f'\vert_{U\cap\pd S}$ approximates a cubic polynomial in $y$
with {\it large third derivative}.

Now $f_a\vert_{\pd S}=\phi\in C^{k+3,\al}(\pd S)$. So by
choosing $J$ large compared to $\cnm{\phi}{3}$ we can ensure that
$\bmd{\frac{\d^3}{\d y^3}f_a\vert_{U\cap\pd S}}\ll 3\bmd{\ka_iq_0}$.
Therefore, provided the approximations are valid,
$\frac{\d^3}{\d y^3}(f'-f_a)\vert_{U\cap\pd S}$ has the same
sign as $\ka_iq_0$ on $U\cap\pd S$, and it easily follows that
$f'-f_a$ has at most two stationary points on~$U\cap\pd S$.

To make this into a rigorous argument, we need to consider the
approximations above very carefully using the estimates of
\S\ref{us45}, and make sure that for small $\ep$ and large $J$,
the effect on the second and third derivatives of $f'$ of the
`error terms' we have neglected are always significantly smaller
than the `leading terms' given in \eq{us7eq10} and \eq{us7eq11}.
We will not do this in detail, as it is long and dull, but here
is a sketch of the important steps.

Firstly, $U$ is approximately a closed disc of radius
$\md{\hat q_0}L$ with centre $(x_0,y_0)-\hat q_0(u_0,v_0)$.
Since $u_0^2+v_0^2\le X^2$ and $L=10X$, it follows that for
$(x,y)\in U$ we have $\md{y-y_0}\le\frac{11}{10}\md{\hat q_0}L$,
approximately. As $\hat q_0\approx q_0^{-1}$, we see from
\eq{us7eq7} that if $(x,y)\in U$ then $\md{y-y_0}\le\frac{11}{10}
J^{-1}L(y_0^2+a^2)^{1/2}$, approximately. So by choosing $J\gg L$
we easily show that
\begin{equation*}
\ts\frac{1}{4}(y_0^2+a^2)\le y^2+a^2\le 4(y_0^2+a^2)
\quad\text{on $U$.}
\end{equation*}
As $(\hat u,\hat v)$ maps $D_L\ra U$, this implies that
\begin{equation*}
\ts\frac{1}{4}(y_0^2+a^2)\le \hat v^2+a^2\le 4(y_0^2+a^2)
\quad\text{on $D_L$.}
\end{equation*}

Now $(\hat u,\hat v)$ is constructed using the solutions of
Theorem \ref{us4thm1}, and these are estimated in terms of a
parameter $s$, which is roughly a lower bound for
$(v^2+y^2+a^2)^{1/2}$. The equation above implies that
$(\hat v^2+y^2+a^2)^{1/2}\ge\ha(y_0^2+a^2)^{1/2}$ on
$D_L$. This means that for the solutions $(\hat u,\hat v)$
the parameter $s$ in \S\ref{us4} is approximately
$(y_0^2+a^2)^{1/2}$, up to a bounded factor. Therefore
we can use Theorem \ref{us4thm3} to estimate the derivatives
of $\hat v$ in terms of powers of~$(y_0^2+a^2)^{1/2}$.

Secondly, some detail on the estimates we need. For the argument
above to work, we need the approximation \eq{us7eq10} to hold up
to third derivatives with errors small compared to $q_0$. That
is, we need
\begin{equation*}
\ts\frac{\pd u'}{\pd x}\!=\!\frac{\pd v'}{\pd y}\!=\!q_0\!+\!o(q_0),\;\>
\frac{\pd u'}{\pd y}=o(q_0),\;\>
\frac{\pd v'}{\pd x}=o(q_0),\;\>
\pd^2u'=o(q_0),\;\>
\pd^2v'=o(q_0),
\end{equation*}
using the $o(\dots)$ order notation in the obvious way. Since
$(\hat u,\hat v)$ is the inverse map of $(u',v')$, calculation
shows these are equivalent to:
\e
\begin{gathered}
\ts\frac{\pd\hat u}{\pd x}=\frac{\pd\hat v}{\pd y}=q_0^{-1}
+o(q_0^{-1}),\quad
\frac{\pd\hat u}{\pd y}=o(q_0^{-1}),\\
\ts\frac{\pd\hat v}{\pd x}=o(q_0^{-1}),\quad
\pd^2\hat u=o(q_0^{-2})\quad\text{and}\quad
\pd^2\hat v=o(q_0^{-2}).
\end{gathered}
\label{us7eq12}
\e

Now $(\hat u,\hat v)$ were constructed in Theorem \ref{us4thm1} by
rescaling the solutions $u,v$ studied in \S\ref{us42}--\S\ref{us46}.
For simplicity let us ignore this rescaling process, and identify
$\hat u,\hat v$ with the solutions $u,v$ of \S\ref{us42}--\S\ref{us46}.
Then from \eq{us4eq8} we have $\hat v=\al+\be y+\ga x+\psi$, where
$\al,\be,\ga\in\R$ with $\be\approx q_0^{-1}$ and $\ga\approx p_0
q_0^{-2}$. Hence we can use Theorem \ref{us4thm3} to estimate the
derivatives of $\hat v$, in terms of powers of $\ga\approx p_0
q_0^{-2}$ and $s\approx(y_0^2+a^2)^{1/2}$. Using \eq{us3eq4} we
can then deduce estimates for the derivatives of~$\hat u$.

Thirdly, we divide into the two cases (a) $\md{v_0}\le(y_0^2+a^2)^{3/8}$,
and (b) $\md{v_0}>(y_0^2+a^2)^{3/8}$. In case (a), we can use Theorem
\ref{us4thm3} to prove \eq{us7eq12}, because \eq{us7eq8} implies that
$p_0=O\bigl((y_0^2+a^2)^{3/8}\md{y_0}\md{q_0}\bigr)$, and this gives
an estimate for $\ga$, which turns out to be exactly what we need.

However, in case (b), the most direct approach is insufficient to
prove \eq{us7eq12} on all of $U\cap\pd S$, so we have to do something
different. We divide $U\cap\pd S$ into three connected regions,
(i) $\md{y-y_0}\le\ha\md{v_0q_0}$; (ii) $y>y_0+\ha\md{v_0q_0}$; and
(iii) $y<y_0-\ha\md{v_0q_0}$. In region (i) we see from \eq{us7eq9}
that $\md{v'}\ge\ha\md{v_0}$, approximately. So when we invert $(u',v')$
to get $(\hat u,\hat v)$, the curve in $D_L$ corresponding to region (i)
will satisfy $\md{y}\ge\ha\md{v_0}$, approximately. Thus when we apply
Theorem \ref{us4thm3} to estimate the `errors', as above, the top line
of \eq{us4eq32} is strong enough to prove that $\frac{\d^3}{\d y^3}
(f'-f_a)\vert_{U\cap\pd S}$ has the same sign as $\ka_iq_0$ on region~(i).

In regions (ii) and (iii), we instead use similar arguments to show
that $\frac{\d^2}{\d y^2}(f'-f_a)\vert_{U\cap\pd S}$ has the same sign
as $\ka_iq_0$ on region (ii), and the opposite sign on region (iii).
Thus we see that $\frac{\d^2}{\d y^2}(f'-f_a)\vert_{U\cap\pd S}$ has
exactly one zero on $U\cap\pd S$, and so $\frac{\d}{\d y}(f'-f_a)
\vert_{U\cap\pd S}$ has at most two zeroes on $U\cap\pd S$. This
completes the proof of Proposition~\ref{us7prop3}.
\end{proof}

Divide the boundary $\pd(U\cap S)$ into two curves $\ga_1=\pd U\cap S$
and $\ga_2=U\cap\pd S$. As $U$ is approximately a closed disc of small
radius $\md{\hat q_0}L$, and $\pd S$ is approximately parallel to the
$y$-axis near $U$, it easily follows that $\ga_1$ and $\ga_2$ are
homeomorphic to $[0,1]$. As we have moved $(x_0,y_0)$ a little way
into the interior of $U\cap S$ it follows that $(u',v')\ne(u_a,v_a)$
on $\pd(U\cap S)$. Define $\th_j$ to be the angle that $(u',v')-(u,v)$
winds around zero along the curve $\ga_j$ for $j=1,2$, where $\ga_j$
has the orientation induced from $\pd U$ or $\pd S$. Then $\th_1+\th_2
=2\pi k$, where $k$ is the winding number of $(u',v')-(u,v)$ about 0
along $\pd U$. Theorem \ref{us7thm3} will follow from the next two
lemmas.

\begin{lem} In the situation above, $\th_1<2\pi$.
\label{us7lem1}
\end{lem}

\begin{proof} From above, $U$ is approximately a closed disc of
radius $\md{\hat q_0}L$ with centre $(x_0,y_0)-\hat q_0(u_0,v_0)$,
and $\pd S$ is near $U$ approximately parallel to the $y$-axis,
and passes very near $(x_0,y_0)$. Since $u_0^2+v_0^2\le X^2$ and
$L=10X$, it follows that $\pd S$ is to first approximation a
straight line which passes within distance $\frac{1}{10}\md{\hat q_0}L$
of the centre of the approximate closed disc $U$. Thus $\ga_1$ occupies
between $\pi-2\sin^{-1}\frac{1}{10}$ and $\pi+2\sin^{-1}\frac{1}{10}$
radians, approximately, of the approximate circle~$\pd U$.

But $(u',v'):U\ra D_L$ takes $\ga_1$ to a portion of the circle
$\pd D_L$, approximately preserving angles. Therefore $(u',v')(\ga_1)$
sweeps out an arc of $\pd D_L$ with angle approximately in the interval
$[2\cdot 94,3\cdot 34]$ radians. Arguing in more detail, using the
estimates on $\hat u,\hat v$ in \S\ref{us4}, one can show that if $J$
is big enough then $(u',v')(\ga_1)$ is no more than 5 radians of the
circle~$\pd D_L$.

Now on $\ga_1$ we have $\bmd{(u',v')}=L$ and $\bmd{(u_a,v_a)}\le
X=\frac{1}{10}L$. As the angle swept out by $(u',v')$ about 0 along
$\ga_1$ is no more than 5 radians, it follows that the angle swept
out by $(u',v')-(u_a,v_a)$ about 0 along $\ga_1$ is no more than $5+2
\sin^{-1}\frac{1}{10}=5\cdot 20$ radians. Thus~$\th_1\le 5\cdot 20<2\pi$.
\end{proof}

\begin{lem} In the situation above, $\th_2<2\pi$.
\label{us7lem2}
\end{lem}

\begin{proof} As above, we may parametrize $\pd S$ near $(x_i,0)$ as
$\bigl(x(y),y\bigr)$ with $x(y)\approx x_i+\ha\ka_iy^2$ for small
$y$. So $\ga_2$ is parametrized by $\bigl(x(y),y\bigr)$ for $y$
in some small interval $[y_1,y_2]$. Therefore
\begin{equation*}
\frac{\d}{\d y}\bigl(f'-f_a\bigr)\vert_{U\cap\pd S}=
\bigl(u'-u_a,v'-v_a\bigr)\cdot\bigl(1,\dot x(y)\bigr)
\quad\text{for $y\in[y_1,y_2]$,}
\end{equation*}
where $\dot x(y)=\frac{\d}{\d y}x(y)$. Now Proposition \ref{us7prop3}
implies that $\frac{\d}{\d y}\bigl(f'-f_a\bigr)\vert_{U\cap\pd S}$
can change sign at most twice on $\ga_2$. This shows that the angle
between $(u'-u_a,v'-v_a)$ and $(1,\dot x)$ can cross over $\pm\frac{
\pi}{2}$ modulo $2\pi$ at most twice on $\ga_2$. So we shall prove
the lemma by comparing the angles that $(u'-u_a,v'-v_a)$ and
$(1,\dot x)$ rotate through about 0 along~$\ga_2$.

For simplicity, suppose that $i=2$, so that $\ga_2$ is oriented
in the direction of increasing $y$ and $\ka_i=\ka_2<0$, and suppose
that $q_0>0$. The other possibilities of $i=1$ or 2 and $q_0>0$ or
$q_0<0$ follow in the same way. Then from the proof of Proposition
\ref{us7prop3}, as $\ka_iq_0<0$ we see that $\frac{\d}{\d y}\bigl(
f'-f_a\bigr)\vert_{U\cap\pd S}$ is negative near $y=y_1$ and $y=y_2$,
and from the proof of Lemma \ref{us7lem1} we see that at $y=y_1$ we
have $v'-v_a<0$, and at $y=y_2$ we have $v'-v_a>0$, and $u'-u_a$
is small compared to $v'-v_a$ at both $y=y_1$ and~$y=y_2$.

If $\frac{\d}{\d y}\bigl(f'-f_a\bigr)\vert_{U\cap\pd S}<0$ on all
of $\ga_2$, then $(u'-u_a,v'-v_a)\cdot(1,\dot x)<0$ on $\ga_2$,
which approximately says that $u'-u_a<0$ on $\ga_2$, as $\dot x$
is small. A little thought shows that the angle $\th_2$ which
$(u'-u_a,v'-v_a)$ rotates through along $\ga_2$ is approximately~$-\pi$.

If $\frac{\d}{\d y}\bigl(f'-f_a\bigr)\vert_{U\cap\pd S}$ changes
sign twice on $\ga_2$, then $(u'-u_a,v'-v_a)$ rotates through an
extra angle of $2\pi,0$ or $-2\pi$ compared to $(1,\dot x)$. Hence
$\th_2$ is approximately $\pi,-\pi$ or $-3\pi$. By explaining what
we mean by $u'-u_a$ being small compared to $v'-v_a$ at $y=y_1$ and
$y=y_2$, we may easily show that~$\th_2<2\pi$.
\end{proof}

We can now finish the proof of Theorem \ref{us7thm3}. From
above $\th_1+\th_2=2\pi k$, where $k$ is the winding number of
$(u',v')-(u,v)$ about 0 along $\pd(U\cap S)$. But $\th_1,\th_2<
2\pi$ by the last two lemmas. Hence $k<2$, and so the winding
number of $(u',v')-(u_a,v_a)$ about 0 along $\pd(U\cap S)$ is
at most 1. But $(u',v')-(u_a,v_a)$ has a zero of multiplicity
at least 2 at $(x_0,y_0)$ in $(U\cap S)^\circ$, so
\cite[Th.~6.7]{Joyc3} gives a contradiction.
\end{proof}

\subsection{Estimates for $\frac{\pd u_a}{\pd x},\frac{\pd u_a}{\pd y},
\frac{\pd v_a}{\pd x},\frac{\pd v_a}{\pd y}$ on $\pd S$ and in $L^p(S)$}
\label{us74}

Next we prove analogues of the estimates of \S\ref{us64}--\S\ref{us65}
in the situation of this section. Combining Theorem \ref{us7thm3} and
the results of \S\ref{us72}, we deduce:

\begin{cor} There exist constants $K_1,K_2>0$ such that
\e
\begin{split}
&\Big\vert\frac{\pd u_a}{\pd x}\Big\vert=
\Big\vert\frac{\pd v_a}{\pd y}\Big\vert\le K_1(y^2+a^2)^{-1/2},\quad
\Big\vert\frac{\pd u_a}{\pd y}\Big\vert\le K_2\quad\text{and}\\
&
\Big\vert\frac{\pd v_a}{\pd x}\Big\vert\le 2K_2(v_a^2+y^2+a^2)^{1/2}
\quad\text{on $\pd S$, for all\/ $a\in(0,1]$.}
\end{split}
\label{us7eq13}
\e
\label{us7cor1}
\end{cor}

\begin{proof} Theorem \ref{us7thm3} implies that 
$\bmd{\frac{\pd v_a}{\pd y}}\le J(y^2+a^2)^{-1/2}$ for
all $a\in(0,1]$ and $(x,y)\in\pd S$ with $\md{y}\le\de$,
for some $\de,J>0$. Apply Proposition \ref{us7prop1} with
$\ep$ replaced by this $\de$. This gives $G>0$ such that
$\bmd{\frac{\pd v_a}{\pd y}}\le G$ for all $a\in(0,1]$ and
$(x,y)\in\pd S$ with $\md{y}>\de$. Combining these shows that
$\bmd{\frac{\pd v_a}{\pd y}}\le K_1(y^2+a^2)^{-1/2}$ on $\pd S$
for all $a\in(0,1]$, where $K_1=\max\bigl(J,G(Y^2+1)^{1/2}\bigr)$.
This proves the first inequality of~\eq{us7eq13}.

Proposition \ref{us7prop2} gives $\ep,H>0$ and functions
$F_{i,a}$ such that if $(x,y)\in\pd S$ is close to $(x_i,0)$
and $\md{y}\le\ep$ then by \eq{us7eq5} we have
\begin{align*}
\Big\vert\frac{\pd u_a}{\pd y}(x,y)\Big\vert
&\le\bmd{F_{i,a}(y)}\Big\vert\frac{\pd u_a}{\pd x}(x,y)\Big\vert+H\\
&\le 3\md{\ka_i}\md{y} K_1(y^2+a^2)^{-1/2}+H\le 3\md{\ka_i}K_1+H,
\end{align*}
using the first inequality of \eq{us7eq13} and $\bmd{F_{i,a}(y)-2\ka_iy}
\le\md{\ka_iy}$, so that $\bmd{F_{i,a}(y)}\le 3\md{\ka_i}\md{y}$.
Applying Proposition \ref{us7prop1} with this value of $\ep$ gives
$G>0$ such that $\bmd{\frac{\pd u_a}{\pd y}}\le G$ for all $a\in(0,1]$
and $(x,y)\in\pd S$ with $\md{y}>\ep$. Combining these shows that
$\bmd{\frac{\pd u_a}{\pd y}}\le K_2$ on $\pd S$ for all $a\in(0,1]$,
where $K_2=\max(3\md{\ka_1}K_1+H,3\md{\ka_2}K_1+H,G)$. This proves
the second inequality of \eq{us7eq13}, and the third follows
from~\eq{us3eq4}.
\end{proof}

As $\md{v_a}\le X$ by \eq{us7eq2} we see from \eq{us7eq13} that
$\bmd{\frac{\pd v_a}{\pd x}}\le K_3$ on $\pd S$ for all $a\in(0,1]$,
where $K_3=2K_2(X^2+Y^2+1)^{1/2}$. Thus from \cite[Prop.~8.12]{Joyc3}
we deduce:

\begin{cor} There exists $K_3>0$ such that\/ $\bmd{\frac{\pd v_a}{\pd
x}}\le K_3$ in $S$ for all\/~$a\in(0,1]$.
\label{us7cor2}
\end{cor}

These two corollaries will serve as an analogue of Corollary
\ref{us6cor2}. Unfortunately we do not have an analogue of 
Proposition \ref{us6prop4} on all of $S$, as Corollary \ref{us7cor1}
is not strong enough to apply Theorem \ref{us5thm4}. However, by
applying Theorem \ref{us5thm1} with $K=1$ and $L=X$ we deduce
the following analogue of Proposition \ref{us6prop4} in interior
domains~$T\subset S^\circ$.

\begin{prop} Let\/ $T\subset S^\circ$ be a subdomain. Then there
exists $K_4>0$ with
\e
\ts\bmd{\frac{\pd u_a}{\pd x}}=\bmd{\frac{\pd v_a}{\pd y}}
\le K_4\bigl(v_a^2+y^2+a^2\bigr)^{-1}
\quad\text{in $T$ for all\/ $a\in(0,1]$.}
\label{us7eq14}
\e
\label{us7prop4}
\end{prop}

The integrals of \eq{us6eq5} are bounded uniformly in~$a$.

\begin{prop} Let\/ $\al\in[0,\ha)$, and define $J(a,v)$ as
in Proposition \ref{us6prop5}. Then there exists $D_\al>0$
depending only on $\al$ such that\/ $\int_{\pd S}J(a,v_a)
\d u_a\le D_\al$ for all\/~$a\in(0,1]$.
\label{us7prop5}
\end{prop}

\begin{proof} We have $\bmd{J(a,v_a)}\le(1-2\al)^{-1}
\md{v_a}^{1-2\al}\le(1-2\al)^{-1}X^{1-2\al}$ as in the
proof of Proposition \ref{us6prop6}, and so
\begin{align*}
\int_{\pd S}J(a,v_a)\d u_a&\le(1-2\al)^{-1}
X^{1-2\al}\int_{\pd S}\md{\d u_a}\\
&=(1-2\al)^{-1}X^{1-2\al}\int_0^1\Big\vert\frac{\pd u_a}{\pd x}
\cdot\frac{\d x}{\d s}+\frac{\pd u_a}{\pd y}\cdot\frac{\d y}{\d s}
\Big\vert\d s,
\end{align*}
where $s\mapsto\bigl(x(s),y(s)\bigr)$ is a parametrization
$[0,1]\mapsto\pd S$ of $\pd S$.

Now $\bmd{\frac{\pd u_a}{\pd x}}\le K_1(y^2+a^2)^{-1/2}$ and
$\bmd{\frac{\pd u_a}{\pd y}}\le K_2$ on $\pd S$ by Corollary
\ref{us7cor1}. Also, when $y$ is small we see that $\bigl(x(s),
y(s)\bigr)$ is close to $(x_i,0)$ for $i=1$ or 2, and then
$\frac{\d x}{\d s}\approx\ka_i y\frac{\d y}{\d s}$. Thus
$\bmd{\frac{\d x}{\d s}}\le C\md{y}\bmd{\frac{\d y}{\d s}}$
for $s\in[0,1]$ and some $C>0$ depending only on $S$. Hence
\begin{equation*}
\int_0^1\Big\vert\frac{\pd u_a}{\pd x}
\cdot\frac{\d x}{\d s}+\frac{\pd u_a}{\pd y}\cdot\frac{\d y}{\d s}
\Big\vert\d s\le
\int_0^1\bigl(K_1(y^2+a^2)^{-1/2}C\md{y}+K_2\bigr)
\Big\vert\frac{\d y}{\d s}\Big\vert\d s.
\end{equation*}
As $(y^2+a^2)^{-1/2}\md{y}\le 1$, combining the last two equations
proves the proposition, with~$D_\al=(1-2\al)^{-1}X^{1-2\al}(K_1C+K_2)
\int_0^1\bmd{\frac{\d y}{\d s}}\d s$.
\end{proof}

Following the proof of Proposition \ref{us6prop6}, we show:

\begin{prop} Let\/ $q\!\in\![1,2)$ and\/ $r\!\in\![1,\iy)$. Then there
exist\/ $C_2,C'_q,C''_r>0$ such that $\lnm{\frac{\pd u_a}{\pd x}}{2}=
\lnm{\frac{\pd v_a}{\pd y}}{2}\le C_2$, $\lnm{\frac{\pd u_a}{\pd
y}}{q}\le C'_q$ and\/ $\lnm{\frac{\pd v_a}{\pd x}}{r}\le C''_r$ for
all\/ $a\in(0,1]$. Let\/ $T\subset S^\circ$ be a subdomain, and let\/
$p\in(2,\frac{5}{2})$. Then there exists $C_p^T>0$ such that
$\lnm{\frac{\pd u_a}{\pd x}\vert_T}{p}=\lnm{\frac{\pd v_a}{\pd
y}\vert_T}{p}\le C_p^T$ for all\/~$a\in(0,1]$.
\label{us7prop6}
\end{prop}

The difference between this and Proposition \ref{us6prop6} is
because the proof of Proposition \ref{us6prop6} for $p>2$
involved Proposition \ref{us6prop4}, but its analogue in this
section is Proposition \ref{us7prop4}, which only holds in
interior domains $T\subset S^\circ$. Hence we can only prove
$L^p$ estimates for $\frac{\pd u_a}{\pd x},\frac{\pd v_a}{\pd y}$
for $p>2$ on such subdomains~$T$.

\subsection{Uniform continuity of the $u_a,v_a$}
\label{us75}

Now we can prove the analogue of Proposition \ref{us6prop7},
which shows that the $u_a,v_a$ are {\it uniformly continuous}
for all $a\in(0,1]$. The proof for $v_a$ is as in \S\ref{us65},
but because we can only bound $\frac{\pd u_a}{\pd x}$ in $L^p$
in interior domains $T\subset S^\circ$, we have to introduce
some new ideas to show that $u_a$ are uniformly continuous
near $(x_i,0)$ for~$i=1,2$.

\begin{thm} In the situation of Definition \ref{us7def1},
there exist continuous $M,N:S\t S\ra[0,\iy)$ satisfying
$M(x',y',x,y)=M(x,y,x',y')$, $N(x',y',x,y)=N(x,y,x',y')$ and\/
$M(x,y,x,y)=N(x,y,x,y)=0$ for all\/ $(x,y),(x',y')\in S$, such
that for all\/ $(x,y),(x',y')\in S$ and\/ $a\in(0,1]$ we have
\e
\begin{split}
\bmd{u_a(x,y)-u_a(x',y')}&\le M(x,y,x',y')\\
\text{and\/}\quad
\bmd{v_a(x,y)-v_a(x',y')}&\le N(x,y,x',y').
\end{split}
\label{us7eq15}
\e
\label{us7thm4}
\end{thm}

\begin{proof} Choose $p>2$ and set $q=2$. Then Theorem \ref{us2thm}
gives continuous functions $G,H:S\t S\ra[0,\iy)$. Combining Proposition
\ref{us7prop6} and Theorem \ref{us2thm}, we see that for all $(x,y),
(x',y')\in S$ and $a\in(0,1]$, we have
\begin{equation*}
\bmd{v_a(x,y)-v_a(x',y')}\le C_p''G(x,y,x',y')+C_2H(x,y,x',y')
=N(x,y,x',y').
\end{equation*}
This defines $N:S\t S\ra[0,\iy)$ satisfying all the conditions of
the theorem.

The next three lemmas construct versions of $M$ on subsets of~$S$.

\begin{lem} Let\/ $T\subset S^\circ$ be a subdomain. Then
there exists a continuous $M^T:T\t T\ra[0,\iy)$ such that for
all\/ $(x,y),(x',y')\in T$ and\/ $a\in(0,1]$ we have
$M^T(x',y',x,y)=M^T(x,y,x',y')$, $M^T(x,y,x,y)=0$ 
and\/~$\bmd{u_a(x,y)-u_a(x',y')}\le M^T(x,y,x',y')$.
\label{us7lem3}
\end{lem}

\begin{proof} Choose $p\in(2,\frac{5}{2})$ and $q\in[1,2)$ with
$p^{-1}+q^{-1}<1$. Then Theorem \ref{us2thm} with $T$ in place of $S$
gives continuous functions $G,H:T\t T\ra[0,\iy)$. Combining Proposition
\ref{us7prop6} and Theorem \ref{us2thm}, we see that for all $(x,y),
(x',y')\in T$ and $a\in(0,1]$, we have
\begin{equation*}
\bmd{u_a(x,y)-u_a(x',y')}\le C_p^TG(x,y,x',y')+C_q'H(x,y,x',y')
=M^T(x,y,x',y').
\end{equation*}
This $M^T$ satisfies the conditions of the lemma.
\end{proof}

\begin{lem} Let\/ $\ep>0$ be small, and let\/ $G$ be as in Proposition
\ref{us7prop1}. Then whenever $(x,y),(x',y')\in S$ and either $y,y'>\ep$
or $y,y'<-\ep$ we have
\e
\bmd{u_a(x,y)-u_a(x',y')}\le G\bigl((x-x')^2+(y-y')^2\bigr)^{1/2}
\quad\text{for all\/ $a\in(0,1]$.}
\label{us7eq16}
\e
\label{us7lem4}
\end{lem}

\begin{proof} Let $S_\ep$ be as in Proposition \ref{us7prop1}. Then
$S_\ep$ has two connected components, with $y>\ep$ and $y<-\ep$,
each of which is convex. The condition that either $y,y'>\ep$ or
$y,y'<-\ep$ ensures that $(x,y)$ and $(x',y')$ lie in the same
connected component of $S_\ep$. Hence the straight line segment
joining $(x,y)$ and $(x',y')$ lies in $S_\ep$, and so $\md{\pd
u_a}\le G$ on this line segment by Proposition \ref{us7prop1}.
Thus, $\bmd{u_a(x,y)-u_a(x',y')}$ is bounded by $G$ times the
length of the segment.
\end{proof}

\begin{lem} There exists a continuous $M^{\pd S}:\pd S\t\pd S\ra
[0,\iy)$ such that for all\/ $(x,y),(x',y')\!\in\!\pd S$ and\/
$a\in(0,1]$ we have $M^{\pd S}(x',y',x,y)\!=\!M^{\pd S}(x,y,x',y')$,
$M^{\pd S}(x,y,x,y)=0$ and\/~$\bmd{u_a(x,y)-u_a(x',y')}\le M^{\pd S}
(x,y,x',y')$.
\label{us7lem5}
\end{lem}

\begin{proof} First suppose $(x,y),(x',y')\in\pd S$ are close
to $(x_i,0)$ for $i=1$ or 2. Then as in \S\ref{us72} we can
parametrize $\pd S$ near $(x_i,0)$ as $\bigl(x(y),y\bigr)$,
and \eq{us7eq3} gives
\begin{equation*}
u_a(x,y)=\dot\phi(y)-\dot x(y)v_a(x,y)\quad\text{and}\quad
u_a(x',y')=\dot\phi(y')-\dot x(y')v_a(x',y'),
\end{equation*}
where `$\,\dot{}\,$' is short for $\frac{\d}{\d y}$. Hence
\begin{align*}
u_a(x,y)-u_a(x',y')=\dot\phi(y)-\dot\phi(y')
-&\ha\bigl(\dot x(y)-\dot x(y')\bigr)
\bigl(v_a(x,y)+v_a(x',y')\bigr)\\
-&\ha\bigl(\dot x(y)+\dot x(y')\bigr)
\bigl(v_a(x,y)-v_a(x',y')\bigr).
\end{align*}
Since $\md{v_a}\le X$ by \eq{us7eq2}, using the second inequality
of \eq{us7eq15} we see that
\begin{align*}
\bmd{u_a(x,y)-u_a(x',y')}\le\bmd{\dot\phi(y)&-\dot\phi(y')}
+X\bmd{\dot x(y)-\dot x(y')}\\
&+\ha\bmd{\dot x(y)+\dot x(y')}N(x,y,x',y').
\end{align*}

As $\dot\phi$ and $\dot x$ are well-defined and continuous
near $(x_i,0)$ on $\pd S$, and $N$ is continuous, the r.h.s.\
of this equation defines a continuous function $M^{\pd S}$ near
$\bigl((x_i,0),(x_i,0)\bigr)$ in $\pd S\t\pd S$ which satisfies
the conditions of the lemma. Also, Lemma \ref{us7lem4} shows how
to define $M^{\pd S}$ on the subsets of $\pd S$ with $y>\ep$ and
$y<-\ep$ for small $\ep>0$. It is then not difficult to patch
together these functions on subsets of $\pd S\t\pd S$ using the
triangle inequality to construct a suitable function $M^{\pd S}$
on all of $\pd S\t\pd S$, and we leave this as an exercise.
\end{proof}

Let us review what we have proved so far. The last three lemmas
show that the $u_a$ are uniformly continuous for all $a\in(0,1]$
in interior domains $T\subset S^\circ$, and away from the $x$-axis
in $S$, and on $\pd S$, respectively. Taken together, these imply
that the $u_a$ are uniformly continuous except possibly at $(x_i,0)$
for $i=1,2$, and that the restrictions of the $u_a$ to $\pd S$ are
uniformly continuous at $(x_i,0)$. It remains only to show that
the $u_i$ are uniformly continuous at~$(x_i,0)$.

By \cite[Cor.~4.4]{Joyc3}, if $u,v$ satisfy \eq{us3eq4} on a
domain $S$ then the maximum of $u$ is achieved on $\pd S$.
It is proved by applying the maximum principle to a second-order
linear elliptic equation satisfied by $u$ in \cite[eq.~(27),
Prop.~4.3]{Joyc3}. Now this equation has no terms in $u$ and
$\frac{\pd u}{\pd x}$, and so if $\al\in\R$ then $u-\al x$
satisfies the same equation. Hence we may prove:

\begin{lem} Let\/ $T\subset S$ be a subdomain, $\al\in\R$
and\/ $a\in(0,1]$. Then the maximum of\/ $u_a-\al x$ on
$T$ is achieved on~$\pd T$.
\label{us7lem6}
\end{lem}

Using this we shall show that the $u_a$ are uniformly continuous
at~$(x_i,0)$.

\begin{lem} Let\/ $i=1,2$ and\/ $\ep>0$. Then there exists
$\de>0$ such that if\/ $(x,y)\in S$ and\/ $\md{x-x_i}<\de$
then $\bmd{u_a(x,y)-u_a(x_i,0)}<\ep$ for all\/~$a\in(0,1]$.
\label{us7lem7}
\end{lem}

\begin{proof} Let $i,\ep$ be as above. Then as $S$ is strictly
convex, and $T_{(x_i,0)}\pd S$ is parallel to the $y$-axis, and
$M^{\pd S}$ is continuous with $M^{\pd S}(x_i,0,x_i,0)=0$, there
exists $\ga>0$ such that $M^{\pd S}(x,y,x_i,0)\le\frac{\ep}{2}$
for all $(x,y)\in\pd S$ such that $\md{x-x_i}\le\ga$. Define
$T=\bigl\{(x,y)\in S:\md{x-x_i}\le\ga\bigr\}$. Then $T$ is a
subdomain of $S$, with piecewise-smooth boundary consisting of
a portion of $\pd S$, and a straight line segment on
which~$\md{x-x_i}=\ga$.

Observe that on $\pd T$ we have
\e
u_a(x_i,0)-\frac{\ep}{2}-\frac{2X}{\ga}\md{x-x_i}\le
u_a(x,y)\le u_a(x_i,0)+\frac{\ep}{2}+\frac{2X}{\ga}\md{x-x_i}.
\label{us7eq17}
\e
This is because on the part of $\pd T$ coming from $\pd S$ we have
$\bmd{u_a(x,y)-u_a(x_i,0)}\le M^{\pd S}(x,y,x_i,0)\le\frac{\ep}{2}$,
by Lemma \ref{us7lem5} and the definition of $\ga$, and on the part
of $\pd T$ with $\md{x-x_i}=\ga$ we have $\frac{2X}{\ga}\md{x-x_i}
=2X$, and the result follows from $\md{u_a}\le X$ by~\eq{us7eq2}.

Now the l.h.s.\ and r.h.s.\ of \eq{us7eq17} are both of the form
$\al x+\be$ on $T$ for $\al,\be\in\R$, since $\md{x-x_i}=x-x_i$
on $T$ if $i=1$, and $\md{x-x_i}=x_i-x$ on $T$ if $i=2$. Therefore
Lemma \ref{us7lem6} implies that as \eq{us7eq17} holds on $\pd T$,
it holds on $T$. Thus $\bmd{u_a(x,y)-u_a(x_i,0)}\le\frac{\ep}{2}+
\frac{2X}{\ga}\md{x-x_i}$ on $T$. Choosing $\de=\min(\ga,\ep\ga/4X)$,
the lemma easily follows.
\end{proof}

We have now shown that the $u_a$ for $a\in(0,1]$ are uniformly
continuous everywhere in $S$. With some effort it can be shown
that one can piece together the various functions constructed
above to construct a continuous function $M:S\t S\ra[0,\iy)$
satisfying the conditions of the theorem, and we leave this
as an exercise for the reader. In fact it is rather easier to
construct $M$ which is {\it lower semicontinuous} rather than
continuous, and lower semicontinuity is all we will need to
show that limits of the $u_a$ are continuous. This completes
the proof of Theorem~\ref{us7thm4}.
\end{proof}

\subsection{Existence and uniqueness of weak solutions
of \eq{us3eq6}}
\label{us76}

We can now follow the proofs of \S\ref{us65} more-or-less
unchanged in the situation of this section, to prove:

\begin{prop} Propositions \ref{us6prop8}--\ref{us6prop10}
hold for $f_a,u_a,v_a$ in \S\ref{us72}, except that in
Proposition \ref{us6prop9} we have $\frac{\pd u_0}{\pd x}
=\frac{\pd v_0}{\pd y}\in L^p(S)$ only for~$p\in[1,2]$.
\label{us7prop7}
\end{prop}

In particular, this gives {\it existence} of a solution $f_0$
of \eq{us3eq6} with weak derivatives, and $f_0\vert_{\pd S}=\phi$.
The reason for the difference with Proposition \ref{us6prop9} is
that in Proposition \ref{us7prop6} we have a priori $L^p$ estimates
for $\frac{\pd u_a}{\pd x},\frac{\pd v_a}{\pd y}$ when $p\in(2,
\frac{5}{2})$ only in {\it interior subdomains} $T\subset S^\circ$.
Here is an analogue of Proposition \ref{us6prop11}, which shows that
weak solutions of the Dirichlet problem for \eq{us3eq6} are {\it unique}.

\begin{prop} Suppose $f,f'\in C^{0,1}(S)$ are weak solutions of\/
\eq{us3eq6} with\/ $a=0$ and\/ $f\vert_{\pd S}=f'\vert_{\pd S}$.
Then~$f=f'$.
\label{us7prop8}
\end{prop}

\begin{proof} Define $u,v,u',v'\in L^1(S)$ to be the weak derivatives
$\frac{\pd f}{\pd y},\frac{\pd f}{\pd x},\frac{\pd f'}{\pd y},
\frac{\pd f'}{\pd x}$. Now as $f,f'\in C^{0,1}(S)$ and $f\vert_{\pd S}=
f'\vert_{\pd S}$, it follows that $f-f'$ lies in the closure of
$C_0^1(S)$ in $C^{0,1}(S)$. Hence there exists a sequence $(\psi_n)$
in $C_0^1(S)$ converging to $f-f'$ in $C^{0,1}(S)$, where
$C_0^1(S)$ is the subspace of $\psi\in C^1(S)$ supported in
$S^\circ$. As $f'$ satisfies \eq{us3eq6} weakly and $\psi_n\in
C^1_0(S)$ we have
\e
-\int_S\frac{\pd\psi_n}{\pd x}\cdot A(0,y,v')\,\d x\,\d y
-2\int_S\frac{\pd\psi_n}{\pd y}\cdot u'\,\d x\,\d y=0.
\label{us7eq18}
\e
But as $\psi_n\ra f-f'$ in $C^{0,1}(S)$ it follows that
$\frac{\pd\psi_n}{\pd x}\ra v-v'$ and $\frac{\pd\psi_n}{\pd x}
\ra u-u'$ in $L^\iy(S)$ as $n\ra\iy$. So letting $n\ra\iy$ in
\eq{us7eq18} and using the Dominated Convergence Theorem shows that
\begin{equation*}
-\int_S(v-v')\cdot A(0,y,v')\,\d x\,\d y-2\int_S(u-u')\cdot
u'\,\d x\,\d y=0.
\end{equation*}
Applying the same argument with $f$ instead of $f'$ yields
\begin{equation*}
-\int_S(v-v')\cdot A(0,y,v)\,\d x\,\d y-2\int_S(u-u')\cdot u\,\d x\,\d y=0.
\end{equation*}

Subtracting the last two equations gives
\begin{equation*}
\int_S(v-v')\cdot\bigl(A(0,y,v)-A(0,y,v')\bigr)
\,\d x\,\d y+2\int_S(u-u')^2\,\d x\,\d y=0.
\end{equation*}
Now $\frac{\pd A}{\pd v}=\bigl(v^2+y^2+a^2\bigr)^{-1/2}$ by
\eq{us3eq7}. Thus if $y\ne 0$ the Mean Value Theorem shows that
$A(0,y,v)-A(0,y,v')=(w^2+y^2)^{-1/2}(v-v')$, for some $w$ between
$v$ and $v'$. Therefore
\begin{equation*}
\int_S(w^2+y^2)^{-1/2}(v-v')^2\d x\,\d y+2\int_S(u-u')^2\,\d x\,\d y=0.
\end{equation*}
As the integrands are nonnegative they are zero, so $u=u'$
and $v=v'$ in $L^1(S)$, and hence $f=f'$ in $C^{0,1}(S)$,
as~$f\vert_{\pd S}=f'\vert_{\pd S}$.
\end{proof}

Note that we have shown in Proposition \ref{us7prop7} that a
solution $f_0$ of the Dirichlet problem for \eq{us3eq6} on
$S$ exists that is twice weakly differentiable and satisfies
\eq{us3eq6} with weak derivatives, but Proposition \ref{us7prop8}
shows that $f_0$ is unique in the possibly larger class of weak
solutions of \eq{us3eq6}, which need only be once weakly
differentiable. So this is a stronger result than we actually need.

Finally, Theorem \ref{us7thm1} follows from Propositions
\ref{us7prop7} and \ref{us7prop8}, and Theorem \ref{us7thm2}
in the same way as Theorem~\ref{us6thm2}.

\end{document}